\documentclass[10pt,a4paper]{article}

\usepackage{graphicx}
\usepackage{epsfig}
\usepackage{amssymb}
\usepackage{amsthm}
\usepackage{amsmath}
\usepackage[left=2.5cm,right=2.5cm,bottom=3.5cm,top=3.5cm]{geometry}
\usepackage{float}
\usepackage{epstopdf}
\usepackage{multirow}
\usepackage{color}
\usepackage{hhline}
\usepackage{booktabs}
\usepackage{bm}
\usepackage{bbm}
\usepackage{newtxmath}
\usepackage{accents}
\usepackage{dsfont}

\newtheorem{theorem}{Theorem}[section]

\newtheorem{remark}{Remark}[section]

\newtheorem*{weakproblem}{Weak problem}

\definecolor{nverde}{RGB}{0,61,0} 
\definecolor{cr1}{RGB}{200,0,0}
\definecolor{cr2}{RGB}{0,0,200}
\definecolor{cr12}{RGB}{100,0,100}

\newcommand{\halb}{\frac{1}{2}}
\newcommand{\xx}{\mathbf{x}}
\newcommand{\dt}{\Delta t}
\newcommand{\nn}{\mathbf{n}}

\DeclareMathOperator{\lapla}{\Delta}
\DeclareMathOperator{\dive}{\nabla \cdot}
\DeclareMathOperator{\gra}{\nabla}
\DeclareMathOperator{\grae}{\nabla}

\DeclareMathOperator{\dS}{\mathrm{dS}}
\DeclareMathOperator{\dV}{\mathrm{dV}}
\newcommand{\abs}[1]{\left| #1 \right|}

\newcommand{\Q}{\mathbf{Q}}
\newcommand{\QI}{\mathbf{Q}_{I}}
\newcommand{\QE}{\mathbf{Q}_{E}}
\newcommand{\Flux}{\mathbf{\mathcal{F}}}
\newcommand{\ww}{\mathbf{w}}
\newcommand{\vel}{\mathbf{u}}
\newcommand{\ke}{k}
\newcommand{\btau}{\boldsymbol{\tau}}
\newcommand{\btaus}{\btau^{\ast}}

\newcommand{\press}{p}
\newcommand{\temp}{\theta}
\newcommand{\hf}{\mathbf{q}}

\newcommand{\pdstar}[1]{\tilde{\press}^{#1}}

\newcommand{\wstar}{\ww^{\ast}}
\newcommand{\wsstar}{\ww^{\ast\ast}}

\newcommand{\pstar}{\press^{\ast}}
\newcommand{\psstar}{\ke^{\ast\ast}}

\newcommand{\vstar}{\vel^{\ast}}

\newcommand{\ent}{H}  
\newcommand{\Id}{\mathbf{I}}

\newcommand{\cv}{P}
\newcommand{\cvs}{h} 
\newcommand{\cvv}[1]{\left| \cv_{#1} \right|}
\newcommand{\cvb}{\partial \cv}

\newcommand{\Vop}{\vmathbb{V}}  
\newcommand{\Cop}{\vmathbb{C}}  

\newcommand{\HONE}{H^1}
\newcommand{\HONEO}{\HONE_{0}}

\newcommand{\PS}[1]{\mathbbm{P}_{#1}}
\newcommand{\CS}[1]{C^{#1}}

\newcommand{\vsh}[1]{v_{#1}}
\newcommand{\tf}{\mathit{z}}
\newcommand{\NMB}{N}
\newcommand{\NE}{\NMB_{P,e}}   

\newcommand{\NDOF}{N^{\textrm{dof}}_{\cv}}
\newcommand{\dof}{\textrm{dof}}
\newcommand{\dofv}[2]{\hat{\vel}_{#1}^{#2}}

\newcommand{\doftau}[2]{\hat{\btau}_{#1}^{#2}}
\newcommand{\dofp}[2]{\hat{\press}_{#1}^{#2}}
\newcommand{\dofpp}[2]{\hat{\mathbf{\press}}_{#1}^{#2}}

\newcommand{\kk}{\kappa} 
\newcommand{\proj}{\Pi^\nabla_{\cv,\kk}}
\newcommand{\projL}{\Pi^0_{\cv,\kk}}
\newcommand{\projLb}{\mathbf{\Pi}^0_{\cv,\kk}}
\newcommand{\projLm}{\Pi^0_{\cv,\kk-1}}
\newcommand{\Rey}{\text{Re}}

\newcommand{\M}{\text{M}}

\allowdisplaybreaks

\renewcommand{\arraystretch}{1.2}	

\begin{document}

\thispagestyle{plain}
\begin{center}
 	\textbf{ \Large{An all Mach number semi-implicit hybrid Finite Volume/ Virtual Element method for compressible viscous flows on Voronoi meshes} }
 	
 	\vspace{0.5cm}
 	{Walter Boscheri$^{(a,b)}$\footnote{Corresponding author.\\ \hspace*{0.25cm} Email addresses: walter.boscheri@univ-smb.fr (W. Boscheri), saray.busto.ulloa@usc.es (S. Busto), michael.dumbser@unitn.it (M. Dumbser)}, Saray Busto$^{(c,d)}$, Michael Dumbser$^{(e)}$}
 	
 	\vspace{0.2cm}
 	{\small
 		\textit{$^{(a)}$ Laboratoire de Mathématiques UMR 5127 CNRS, Universit{\'e} Savoie Mont Blanc, 73376 Le Bourget du Lac, France}
 		
 		\textit{$^{(b)}$ Department of Mathematics and Computer Science, University of Ferrara, 44121 Ferrara, Italy}
 		
 		\textit{$^{(c)}$ Department of Applied Mathematics, Universidade de Santiago de Compostela, 15782 Santiago de Compostela, Spain}
 		
 		\textit{$^{(d)}$ Galician Center for Mathematical Research and Technology, CITMAga, 15782 Santiago de Compostela, Spain }
 		
 		\textit{$^{(e)}$ Department of Civil, Environmental and Mechanical Engineering, University of Trento, 38123 Trento, Italy}
 	}
\end{center}

\begin{abstract}
We present a novel high order semi-implicit hybrid finite volume/virtual element numerical scheme for the solution of compressible flows on Voronoi tessellations. 
The method relies on the flux splitting of the compressible Navier-Stokes equations into three sub-systems: a convective sub-system solved explicitly using a finite volume (FV) scheme, and the viscous and pressure sub-systems which are discretized implicitly at the aid of a virtual element method (VEM).  
Consequently, the time step restriction of the overall algorithm depends only on the mean flow velocity and not on the fast pressure waves nor on the viscous eigenvalues. 
As such, the proposed methodology is well suited for the solution of low Mach number flows at all Reynolds numbers. 
Moreover, the scheme is proven to be globally energy conserving so that shock capturing properties are retrieved in high Mach number flows. To reach high order of accuracy in time and space, an IMEX Runge-Kutta time stepping strategy is employed together with high order spatial reconstructions in terms of CWENO polynomials and virtual element space basis functions.
The chosen discretization techniques allow the use of general polygonal grids, a useful tool when dealing with complex domain configurations. 
The new scheme is carefully validated in both the incompressible limit and the high Mach number regime through a large set of classical benchmarks for fluid dynamics, assessing robustness and accuracy.
\end{abstract}

\vspace{0.2cm}

\vspace{0.4cm}

\section{Introduction}\label{sec:intro}
Fluids are involved in numerous natural phenomena and industrial processes including, e.g., geophysical applications, as weather forecasting and pollution, flood studies in biomedicine, optimization of energy devices or aerodynamics design. As a consequence, understanding their behavior is of great importance for the development of the society. 
From the mathematical point of view, fluids are typically modeled by the Euler equations, which are directly derived from the fundamental principles of mass, momentum and energy conservation. 
Then, under appropriate assumptions for the viscous stress tensor and the heat flux \cite{Ber05}, the compressible Navier-Stokes equations are derived, which constitute the most extended model for the simulation of viscous flows. 
Even if these systems have been put forward more than 200 years ago, the development of accurate schemes for their solution is still a major field of research in applied mathematics and engineering, especially in the context of multiple space and time scales that might arise in real world applications.

Traditionally, numerical schemes for fluid flows have been divided into two big families: pressure-based and density-based solvers. The election between these two kinds of numerical methods was initially motivated by the compressibility properties of the flow. Generally, fluids are considered to be a compressible medium where even strong shock waves may arise and they have been commonly discretized using explicit finite volume (FV) or discontinuous Galerkin (DG) density-based solvers, see e.g. \cite{God59,LW60,HLL83,Toro,Munz94,BR97b,Dol04}. However, depending on the Mach number $\M={\left| \vel \right|}/{c}$, i.e. the ratio between the mean flow velocity and the sound speed in the medium, we can identify different flow regimes.
An asymptotic analysis of the equations as $M\rightarrow 0$ yields to the so-called incompressible regime with the well-known divergence-free condition for the velocity field \cite{KlaMaj,KlaMaj82}. Moreover, this analysis shows that the pressure in the incompressible limit can be decomposed into a constant plus a fluctuation that is governed by an elliptic Poisson equation. As such, the nature of the system changes with respect to the original hyperbolic-parabolic compressible Navier-Stokes equations and the incompressible system is often discretized using semi-implicit pressure-based solvers \cite{HW65,Cho67,Pat80,BCG89,Cas14,CB24}, or implicit finite element methods (FE) \cite{TH73,For81,HR82,HMM86,HR88,Guer06}. 

One important feature of the explicit density-based solvers is that they are built upon the integration of the governing equations on spatial control volumes, hence ensuring conservation properties by construction, and they are well suited for capturing strong discontinuities and shock waves typical of the high Mach regime. However, in the low Mach regime these schemes may present excessive numerical diffusivity due to an incorrect scaling with respect to the Mach number of the flow under consideration \cite{GM04,Dellacherie1}. To extend the density-based methodology to the low Mach regime, several corrections have been proposed in the bibliography, see e.g. \cite{Thornber,NRKCC16,AdamsLowMach,barsukow2021AllMach}. Still, the allowed time step depends on a CFL stability condition which accounts for the pressure waves. Consequently, these schemes are quite expensive in the incompressible limit. To circumvent this issue, fully implicit approaches and preconditioning techniques may be employed \cite{LG08}, but then the difficulty lies in the treatment of the nonlinear upwind discretization of the convective terms.
On the other hand, pressure-based solvers, which generally allow larger time steps due to the implicit treatment of the pressure sub-system, were typically based on a non-conservative form of the equations which limited their capability to deal with highly compressible flows, and should be applied only in the low Mach number regime \cite{CG84,Klein2001,MRKG03}. 

The first semi-implicit pressure-based schemes able to deal with both incompressible and high Mach number flows have been presented in \cite{HVW03,ParkMunz2005}. The seminal idea is the use of a conservative version of the equations which is then split into a convective sub-system, to be discretized explicitly, and a pressure sub-system, solved implicitly. 
Following also this idea, a novel flux-vector splitting strategy was then proposed in \cite{TV12} for the Euler equations obtaining an explicit sub-system whose CFL condition is independent of the sound speed. 
Since then, several authors have focused on the development of the so-called \textit{all Mach number solvers}, using different families of numerical methods as finite volume schemes \cite{CordierDegond,CGK13,DC16,RussoAllMach,AbateIolloPuppo,ABIR19,BDT2021,LPT23} or discontinuous Galerkin methods \cite{TD17,BTBD20,BT22,TB24}. Semi-implicit time marching algorithms fall into the more general framework of implicit-explicit (IMEX) time integrators \cite{AscRuuSpi,Hofer,PR_IMEX,BFR16,BP17}, that easily permit to account for multiple time scales coexisting in the flow.

The particular splitting of the equations into a hyperbolic system, where the momentum is unknown, and an elliptic problem, for the unknown pressure, has also recently motivated the development of hybrid Finite Volume/Finite Element methods \cite{Hybrid2}. This approach profits from the conservation and shock capturing properties of Godunov methods for the solution of the convective system. Meanwhile, implicit continuous finite element methods, well known for their ability to treat elliptic and parabolic problems, are employed to solve the pressure system. To perform the spatial discretization of complex geometries, unstructured face-based staggered grids generated from a primal triangular or tetrahedral mesh are considered. However, a shortcoming of this hybrid methodology is the difficulty related to its extension to arbitrary polygonal grids. More precisely, even if the finite volume approach used to treat the convective system can be directly conveyed to arbitrary elements, the discretization of the pressure system using continuous finite elements is not straightforward for general grids. A possibility may be the use of a primal grid made of polygonal elements for the convective FV stage, which is then combined with a staggered subgrid made of triangles for the pressure sub-system, hence finite elements can be employed \cite{BTC23}. 

A powerful alternative is given by the use of virtual element methods (VEM) for the discretization of the Poisson-type sub-system following the novel framework put forward in \cite{HybridFVVEMinc} for the incompressible Navier-Stokes and shallow water equations. Virtual element methods can be seen as a generalization of finite element methods where the basis functions are not known explicitly easing their extension to general polytopal grids. The method requires the definition of adequate projector operators onto polynomial spaces which directly employ the known degrees of freedom for the computation of the integrals involved in the variational formulation of the equations. 
This approach has been widely employed for the implicit solution of elliptic equations in solid mechanics as for linear elasticity, elastodynamics and fracture problems \cite{VBM13,GTP14,DLV20,NZNRW18,HAHWGA19,BCMS22,AMMMV21} 
and for the simulation of porus media \cite{BHCWS21,BBKMS22}. 
However, the use of VEM for time-dependent partial differential equations remains almost unexplored. 
In fluid mechanics, the firstly developed VEM schemes dealt with the steady state Navier-Stokes equations \cite{WWW21,CMM21,WMWH21,BDPD22,AVV23,ABBVV24}, while their use within a numerical method for the solution of time evolving flows has been recently presented in \cite{HybridFVVEMinc}. Only very recently, a space-time VEM-DG method has been devised in \cite{antonietti2023cvem} for the linear scalar dissipative wave equation, while in \cite{VEMDG_heat2024,VEMDG_heat2024_2} a genuinely space-time VEM scheme has been proposed for the one-dimensional heat conduction equation. Besides, an innovative methodology combining the fundamentals of VEM and DG schemes has been proposed in \cite{VEMDGVoronoi} for the incompressible and compressible Navier-Stokes equations. This VEM-DG approach employs the virtual element framework to construct the solution in each single cell and thus performs a local projection. Consequently, a nonconforming representation of the solution with discontinuous data across the element boundaries is obtained making this method specially well-suited for the study of fluid problems presenting discontinuities and shock waves.

Following \cite{HybridFVVEMinc}, in this work the VEM is extended to compressible flows not only for the discretization of the pressure system but also for the implicit computation of the viscous terms in the Navier-Stokes equations. The resulting pressure algebraic system is symmetric, thus permitting to use very efficient iterative solvers like the conjugate gradient method. Furthermore, the CFL condition depends linearly on the characteristic grid size and on the bulk velocity instead of presenting the standard quadratic dependency on the characteristic mesh size arising in the explicit solution of a combined transport-diffusion system. This feature is particularly useful when dealing with low Reynolds numbers since it reduces the computational cost of the overall algorithm \cite{FL94,RBJ95,FFL95,SF96,LBS07,HybridImplicit}. Our novel numerical method is an asymptotic preserving (AP) all Mach number flow solver for general polygonal grids able to deal with both low and high Reynolds number flows. The numerical solution is stored at cell centers, with no use of staggered meshes, and it is transferred via suitable $L_2$ projection operators from the FV to the VEM solution space and vice-versa.
Moreover, high order of accuracy in time is achieved by means of the class of implicit-explicit Runge-Kutta (IMEX-RK) schemes \cite{ARW95,BFR16,BQRX19},  which overtake the well-known decrease of accuracy order stemming from the use of flux-splitting methodologies, as demonstrated in \cite{BP2021,BT22}. High order in space is instead reached for the finite volume method relying on Central WENO (CWENO) schemes, which have been originally proposed in \cite{LPR99,LPR00}, and subsequently extended to deal with general polygonal meshes in \cite{CWENOGBK,FVBoltz}.

The rest of the paper is organized as follows. The compressible Navier-Stokes equations are recalled in Section~\ref{sec:goveq} while the novel hybrid FV/VEM method is introduced in Section~\ref{sec:numdisc}. First, we present the splitting of the equations into the convective, viscous and pressure sub-systems giving an overview of the overall methodology and describing the time integration strategy, Sections~\ref{sec.fluxsplitting}-\ref{sec.imex}. 
Then, the spatial discretization is introduced. The explicit finite volume scheme for the convective sub-system is described in Section~\ref{sec.convective}. The virtual element scheme developed for the discretization of the viscous and pressure sub-systems is detailed in Section~\ref{sec.viscouspressure}.
The proposed hybrid FV/VEM algorithm is validated in Section~\ref{sec:numericalresults} where different test cases going from the compressible regime to the incompressible limit are analyzed, for different Reynolds numbers. Finally, in Section~\ref{sec:conclusions}, we draft the conclusions and provide an outlook to future research lines.

\section{Mathematical model} \label{sec:goveq}
Let us consider a computational domain $\Omega \subset \mathds{R}^d$ in $d=2$ space dimensions, with boundary $\partial \Omega \subset \mathds{R}^{d-1}$. The spatial position vector is $\xx=(x,y) \in \mathds{R}^2$, and  $t \in \mathds{R}^+$ denotes the time coordinate. The mathematical model is given by the compressible Navier-Stokes equations which describe the mass, momentum and energy conservation:
\begin{subequations}\label{eqn.navierstokes}
\begin{align}
	\frac{\partial \rho}{\partial t} + \dive \left( \rho\vel\right)=0,\label{eq:mass1}\\
	\frac{\partial \rho\vel}{\partial t} + \dive \left( \rho\vel\otimes \vel\right)  + \grae \press - \dive \btau = \mathbf{0},\label{eq:momentum1}\\
	\frac{\partial \rho E}{\partial t} + \dive \left[ \vel\left(\rho E + \press \right) \right]   - \dive \left(\btau\vel\right) +\dive \hf = 0,\label{eq:totenerg1}
\end{align}
\end{subequations}
with $\rho$ being the density, $\vel=(u,v)$ the velocity, $\press$ the pressure, and $E$ the total energy. Moreover, $\mu$ is the dynamic viscosity of the fluid and $\btau$ is the tensor of the viscous stresses, which, after introducing the identity matrix $\mathbf{I}$, writes 
\begin{equation}
	\btau = \mu \left(\gra \mathbf{u} + \gra \vel^{T} \right) -\frac{2}{3} \mu \left( \dive \vel  \right) \mathbf{I}, \label{eq:stresstensor}
\end{equation}
and $\hf$ corresponds to the heat flux,
\begin{equation}
	\hf = -\lambda \gra \temp, \label{eq:heatflux}
\end{equation}
with $\temp$ the temperature and $\lambda$ the thermal conductivity. 
To close the system an extra equation relating the pressure and the density is needed. In particular, we consider the
 ideal gas equation of state (EOS), that reads 
\begin{equation}
	\press = \rho R \temp, \label{eq:stateequation}
\end{equation}
with $R = c_p - c_v$ the specific gas constant being $c_{p}$ and  $c_{v}$ the heat capacity at constant pressure and at constant volume, respectively. Using \eqref{eq:stateequation}, the relation between the total energy $E$, the kinetic energy $k$, and the internal energy $e$, results 
\begin{equation}
	\rho E = \rho e +\rho \ke = \frac{1}{\gamma-1}\press + \halb\rho \left| \vel \right|^{2}, \label{eq:E.as.pk}
\end{equation}
where $\gamma={c_{\press}}/{c_{v}}$ denotes the ratio of specific heats.  
By substituting \eqref{eq:E.as.pk} in \eqref{eq:totenerg1}, the total energy density conservation equation can be rewritten in terms of the pressure and the kinetic energy as
\begin{equation}
	\frac{1}{\gamma-1} \frac{\partial  \press  }{\partial t}
	+ \frac{\partial \rho\ke }{\partial t}
	 + \dive \left( \rho \ent \vel + \rho \ke \vel \right)   
	 - \dive \left(\btau\vel\right) 
	 +\dive \hf 
	 = \mathbf{0}, \label{eq:presss}
\end{equation}
with the specific enthalpy
\begin{equation}
	\ent=\dfrac{\gamma\, \press}{\left( \gamma-1\right) \rho}.
	\label{eqn.ent}
\end{equation}

\section{Numerical method} \label{sec:numdisc}
To discretize the system \eqref{eq:mass1}-\eqref{eq:momentum1}-\eqref{eq:presss}, we propose a novel hybrid finite volume/virtual element method (FV/VEM) on general polygonal grids which extends the approach presented in \cite{HybridFVVEMinc} for incompressible and shallow water flows. First, we apply a time discretization to the governing equations which naturally leads to a splitting of the original system \cite{TD17,Hybrid2,BP2021,HybridNNT,BT22}. Then, the resulting convective sub-system is discretized in space using classical finite volumes while the Poisson-type sub-system obtained for the pressure as well as the viscous sub-system in the momentum equation are solved employing a virtual element method.

\subsection{Semi-discrete semi-implicit scheme}\label{sec.fluxsplitting}
Let $T= [0, t_f ]$ be the time interval, where $t_f \in \mathds{R}^+_0$ denotes the final time, and let $t \in T$. The time interval is discretized by means of a sequence of discrete points $t^n$ such that
\begin{equation}
	t^{n+1}=t^n + \dt,
\end{equation}
where $\dt=t^{n+1}-t^n$ denotes the time step. Let us represent with a super-index $n$ the approximation of the solution at time $t^{n}$, e.g.  $\rho^{n}=\rho\left( \xx,t^{n}\right)$. Performing a time discretization of \eqref{eq:mass1}-\eqref{eq:momentum1}-\eqref{eq:presss}, we get the the semi-discrete scheme
\begin{subequations}\label{eqn.navierstokes_sd}
	\begin{align}
		\frac{\rho^{n+1} -\rho^n}{\dt} = \dive \ww^n, \label{eqn.mass_sd}\\
		\frac{\ww^{n+1} -\ww^{n}}{\dt}  = 
		-\dive \left(\frac{1}{\rho^{n}}\ww^n\otimes\ww^{n}\right) 
		- \gra \press^{n+1} 
		+ \dive \btaus,\label{eqn.momentum_sd}\\
		\frac{1}{\dt  }\left[\frac{1}{ \left( \gamma-1\right) }\left( \press^{n+1} -\press^{n}\right)  
		+ \halb \vel^{n+1}\cdot\ww^{n+1} 
		- \rho^{n}\ke^{n}\right]  = 
		- \dive \left( \ke^{n} \ww^{n} \right) 
		- \dive \left(\ent^n \ww^{n+1}\right)
		+ \dive \left( \btaus\vstar\right)  -\dive \hf^{n}, \label{eqn.press_sd}		
	\end{align}
\end{subequations}
where $\ww$ denotes the linear momentum, i.e. $\ww = \rho\vel$. The new density is readily available from \eqref{eqn.mass_sd}, thus $\rho^{n+1}$ is updated explicitly. The momentum equation \eqref{eqn.momentum_sd} can be split into three contributions: (i) the convective equation, which does not depend on the pressure gradient nor on the viscous terms; (ii) the diffusion equation that incorporates the contribution of the viscous terms; (iii) the pressure equation that accounts for the pressure gradient. Hence, a first intermediate approximation for the linear momentum is given by the pure advection part, that is,
\begin{eqnarray}
	\wsstar  = \ww^{n}  - \dt \dive \left(\frac{1}{\rho^{n}}\ww^n\otimes\ww^{n}\right), \label{eqn.intermediatemomentum}
\end{eqnarray}
which is then summed to the viscous contribution, thus yielding
\begin{eqnarray}
	\wstar  = \wsstar +\dt \dive \btaus, \label{eqn.intermediatemomentumv}
\end{eqnarray}
and finally the new momentum is obtained by taking into account the pressure gradient as
\begin{eqnarray}
	\ww^{n+1}  = \wstar - \dt \gra \press^{n+1}. \label{eqn.velupdate}
\end{eqnarray}

The new pressure is still unknown, therefore let us gather the convective and viscous terms of the pressure equation \eqref{eqn.press_sd} in an auxiliary variable as
\begin{eqnarray}
	\pstar    = \rho^{n}\ke^{n} - \dt \dive \left(\ke^{n} \ww^{n} \right) + \dt \dive \left( \btaus\vstar\right)   - \dt \dive \hf^{n}, \label{eqn.pressure_td}
\end{eqnarray}
hence yielding
\begin{eqnarray}
	\frac{1}{ \gamma-1 }\press^{n+1} + \dt \dive \left(\ent^n \ww^{n+1}\right) 
	=\frac{1}{ \gamma-1 } \press^{n} + \pstar - 	\halb \vel^{n+1}\cdot\ww^{n+1}. \label{eqn.press_sdf}
\end{eqnarray}
Equations \eqref{eqn.velupdate} and \eqref{eqn.press_sdf} form a linear system that can be solved by formally inserting the momentum equation \eqref{eqn.velupdate} in the pressure equation \eqref{eqn.press_sdf}, hence leading to a wave equation for the scalar pressure field $\press^{n+1}$ as unknown:
\begin{gather}
	\frac{1}{ \gamma-1 }\press^{n+1} 	
	- \dt^2 \dive \left(\ent^n  \gra \press^{n+1} \right)
	= \frac{1}{ \gamma-1 } \press^{n} 
	+ \pstar - \halb \vel^{n+1}\cdot\ww^{n+1}	
	- \dt \dive \left(\ent^n  \wstar \right). \label{eqn.pressure_sd}
\end{gather}
The corresponding solution can be used within a correction stage in \eqref{eqn.velupdate} to update the intermediate velocity $\wstar$ with the contribution of the new pressure gradient. 

Let us notice that the third term in the right-hand side of \eqref{eqn.pressure_sd} involves a nonlinear contribution in the velocity at the new time step, which is still unknown and corresponds to the kinetic energy at the new time level. To avoid the solution of the resulting nonlinear system, a Picard iterative method can be used as proposed in \cite{DC16}, otherwise also a semi-implicit time linearization can be performed along the lines of \cite{BP2021}. In this work, we propose a different strategy. The pressure wave equation \eqref{eqn.pressure_sd} is solved firstly to find the \textit{pressure flux} $\pdstar{n+1}$ by using the intermediate momentum $\wstar$ for the computation of the new kinetic energy:
\begin{gather}
	\frac{1}{ \gamma-1 }\pdstar{n+1} 	
	- \dt^2 \dive \left(\ent^n  \gra \pdstar{n+1} \right)
	= \frac{1}{ \gamma-1 } \press^{n} 
	+ \pstar - \halb \frac{\wstar\cdot\wstar}{\rho^{n+1}}	
	- \dt \dive \left(\ent^n  \wstar \right) . \label{eqn.pressure_sd2}
\end{gather}
Then, we update the linear momentum as
\begin{eqnarray}
	\ww^{n+1}  = \wstar - \dt \gra \pdstar{n+1}. \label{eqn.velpicupdate}
\end{eqnarray}
Finally, the pressure wave equation \eqref{eqn.pressure_sd} is solved once again for the \textit{pressure state} $\press^{n+1}$ that guarantees thermodynamic consistency and conservation of the new energy $(\rho E)^{n+1}$, defined by \eqref{eq:E.as.pk}. Hence, we address
\begin{gather}
	\frac{1}{ \gamma-1 }\press^{n+1} 	
	- \dt^2 \dive \left(\ent^n  \gra \press^{n+1} \right)
	= \frac{1}{ \gamma-1 } \press^{n} 
	+ \pstar - \halb \frac{\ww^{n+1}\cdot\ww^{n+1}}{\rho^{n+1}}	
		- \dt \dive \left(\ent^n  \wstar \right), \label{eqn.pressure_sd3}
\end{gather}
with the new kinetic energy that is readily computed from the new momentum \eqref{eqn.velpicupdate}. The new energy $(\rho E)^{n+1}$ follows from its definition \eqref{eq:E.as.pk}, and we will prove that it is globally conserved by our scheme.

\begin{remark}[Asymptotic preserving property]
Let us remark that the semi-discrete scheme \eqref{eqn.velupdate}-\eqref{eqn.pressure_sd} is asymptotic preserving in the incompressible limit of the model by construction. Indeed, letting the Mach number to zero ($\M \to 0$), we have that the sound speed $c$ tends to infinity, that is
\begin{equation}
	c^2 =\gamma\frac{\press}{\rho} =\left( \gamma-1\right) \ent \rightarrow \infty,
\end{equation}
where we use the definition \eqref{eqn.ent} for the specific enthalpy. Further, assuming a constant density, then the enthalpy tends to a constant. Hence, retaining only the terms of order $c^2$ in \eqref{eqn.pressure_sd}, we get
\begin{equation}
	\dt \lapla \press^{n+1} =  \dive  \wstar,  \label{eqn.incpressure_sd}
\end{equation}
which corresponds to the elliptic pressure equation for the incompressible Navier-Stokes equations, see \cite{HybridFVVEMinc}.

Moreover, the asymptotic preserving property of the semi-discrete scheme for the viscous stiff limit is analogous to the proof of the asymptotic preserving property for the incompressible Navier-Stokes equations in \cite[Theorem~3]{HybridFVVEMinc}, therefore the semi-discrete scheme \eqref{eqn.velupdate} is asymptotic preserving for vanishing Reynolds numbers.	
\end{remark}

\subsection{Overall methodology}\label{sec.overall} 
The semi-implicit hybrid FV/VEM scheme is organized into three stages.
\begin{enumerate}
	\item Convection stage. An explicit finite volume method is employed to solve \eqref{eqn.intermediatemomentum} obtaining the first intermediate value for the momentum, i.e. $\wsstar$. Likewise, a first intermediate kinetic energy $\psstar$ is computed from \eqref{eqn.pressure_td} neglecting the viscous terms, that is
	 \begin{eqnarray}
	 	\psstar    = \rho^{n}\ke^{n} - \dt \dive \left(\ke^{n} \ww^{n} \right)   - \dt \dive \hf^{n}, 
	 	\label{eqn.pressure_tdsv}
	 \end{eqnarray}
    which accounts for the advection contribution of the kinetic energy and for the heat flux.
	
	\item Viscous stage. The viscous sub-system \eqref{eqn.intermediatemomentumv} is solved implicitly at the aid of a virtual element method (VEM), thus obtaining the intermediate momentum $\wstar$. This is then used to calculate the work of the viscous forces in \eqref{eqn.pressure_td}, hence obtaining the intermediate pressure 
	\begin{eqnarray}
		\pstar = \psstar + \dt \dive \left( \btaus\vstar\right). \label{eqn.pressure_tdwv}
	\end{eqnarray}
	
	\item Pressure stage. First, the sub-system \eqref{eqn.pressure_sd2} is solved for the \textit{pressure flux} $\pdstar{n+1}$ relying on the VEM. The new momentum $\ww^{n+1}$ can thus be updated by means of \eqref{eqn.velpicupdate}. Next, the new \textit{pressure state} $\press^{n+1}$ is obtained upon the solution of the sub-system \eqref{eqn.pressure_sd3}, that is discretized once again with the VEM approach.
	
\end{enumerate}

\subsection{High order extension in time: semi-implicit IMEX scheme}\label{sec.imex}
The semi-discrete scheme \eqref{eqn.velupdate}-\eqref{eqn.pressure_sd} leads to a first order scheme in time due to the split of the fluxes into explicit and implicit contributions. To attain high accuracy in time, a semi-implicit Implicit-Explicit (IMEX) Runge-Kutta methodology is employed. IMEX schemes have been widely used in the last years to achieve high order of accuracy when a splitting technique is adopted for the discretization of PDE systems, see e.g. \cite{PR_IMEX,BP17,BQRX22,OBB23,BTC23}. As shown in the previous section, an implicit discretization is chosen for the potentially stiff terms, namely the pressure terms in \eqref{eqn.momentum_sd}-\eqref{eqn.pressure_sd} and the viscous terms in \eqref{eqn.momentum_sd}, while an explicit treatment is adopted for the convective terms. Therefore, the framework of the semi-implicit IMEX schemes presented in \cite{BFR16} fits our time discretization. IMEX techniques belong to the Method-of-Lines (MOL) integrators, and they are based on a multi-step time stepping method characterized by two Butcher tableaux: one related to the explicit scheme and a second one for the implicit scheme, namely
\begin{center}
	\begin{tabular}{c | c}
		$\tilde{c}$ & $\tilde{A}$ \\ \hline
		& $\tilde{b}$
	\end{tabular}
	\hspace{1cm}
	\begin{tabular}{c | c}
		$c$ & $A$ \\ \hline
		& $b$
	\end{tabular}
\end{center}
where $\tilde{A}=\left(\tilde{a}_{ij}\right)\in\mathbb{M}_{s \times s}$ is a lower triangular matrix with null diagonal elements, $A=\left(a_{ij}\right)\in\mathbb{M}_{s\times s}$ is a lower triangular matrix, $\tilde{b},\, b \in \mathbb{M}_{s \times 1}$ are the weight vectors and $s$ indicates the number of implicit Runge-Kutta stages.

Following \cite{BFR16}, we rewrite the governing equations \eqref{eqn.navierstokes} as an autonomous system of the form
\begin{equation}
	\frac{\partial \Q}{\partial t} = \mathcal{H}\left(\QE(t),\QI(t)\right), \label{eqn.autonomous}
\end{equation} 
where $\Q=(\rho,\rho \vel,\rho e + \rho k)^\top$ is the vector of unknowns and $\mathcal{H}$ represents the spatial discretization of the convective, diffusive and pressure fluxes of \eqref{eqn.navierstokes}. The first argument of $\mathcal{H}$, denoted with $\QE$, is discretized explicitly, and the second argument, referred to as $\QI$, is taken implicitly, according to the flux splitting introduced in Section~\ref{sec.fluxsplitting}, thus obtaining a partitioned system:
\begin{equation}
	\left\{\begin{aligned}
		\frac{\partial \QE}{\partial t} &=  \mathcal{H}\left(\QE, \QI \right), \\[0.7pt]
		\frac{\partial \QI}{\partial t} &=  \mathcal{H}\left(\QE, \QI \right).
	\end{aligned} \right.
	\label{eqn.Hauto}
\end{equation}
Even if it may seem that the number of unknowns has been doubled, since $\mathcal{H}$ does not present an explicit time dependency, only one set of stage fluxes needs to be computed \cite{BFR16}.
More precisely, the stage fluxes for the semi-implicit IMEX scheme are calculated as
\begin{equation}
		\Q_E^i = \Q^n + \dt \sum \limits_{j=1}^{i-1} \tilde{a}_{ij} k_j, \quad
		\tilde{\Q}_I^i = \Q^n + \dt \sum \limits_{j=1}^{i-1} a_{ij} k_j, \quad 
		k_i = \mathcal{H} \left( \Q_E^i, \tilde{\Q}_I^i + \dt \, a_{ii} \, k_i \right), \quad 1 \leq i \leq s, \label{eq.IMEXstages} 
\end{equation}
and the final solution reads
\begin{equation}
	\Q^{n+1} = \Q^n + \dt \sum \limits_{i=1}^s b_i k_i.
	\label{eqn.QRKfinal}
\end{equation}

In this work, stiffly accurate schemes are used, so that the final solution $\Q^{n+1}$ coincides with the last stage value computed from \eqref{eq.IMEXstages} for $i=s$, and asymptotic preserving properties are proven to hold \cite{PR_IMEX}. More specifically, we have selected the LSDIRK2(2,2,2) and the SA DIRK(3,4,3) IMEX schemes of second and third order of accuracy, respectively.
Let us recall that the triplet $(s,\tilde{s},p)$ next to the name of the IMEX scheme indicates the number of stages of the implicit method, $s$, the number of stages of the explicit method, $\tilde{s}$, and the order of the resulting scheme, $p$. Moreover, DIRK stands for Diagonally Implicit Runge-Kutta schemes, LS indicates L-Stability and SA refers to Stiffly Accurate. 
A comprehensive description of the stages of these two IMEX schemes is given in \ref{app.imex}.

\subsection{General unstructured mesh}\label{sec.grid}
The computational domain $\Omega$ is discretized employing Voronoi meshes. Our numerical scheme does not require the orthogonality property of Voronoi meshes, so it can be applied to any type of polygonal tessellation. Consequently, we assume that $\Omega$ is paved with a general unstructured mesh which counts a total number $N$ of non-overlapping polygonal control volumes $\cv_{i}$ with $i\in\left\lbrace 1,\dots, N\right\rbrace$. The tessellation of $\Omega$ is then given by
\begin{equation}
	\mathcal{T}_{\Omega}=\bigcup\limits_{i=1}^{N} \cv_{i}.
\end{equation}
The boundary of $\cv_{i}$ is denoted by $\partial \cv_{i}$ and is defined by the outward pointing normal vector $\nn$, and $\left|\cv_{i} \right|$ identifies the cell volume. The barycenter of the cell is computed as
\begin{equation}
	\xx_{\cv_i}=\frac{1}{\left|\cv_i\right|} \int\limits_{P_i} \xx \dV.
\end{equation}
Let $e_k$ and $\xx_k=\left(x_k,y_k\right)$ represent the edges and vertices of $\cv_{i}$, respectively, with $k\in\left\lbrace 1,\dots, N_{\cv_{i},e}\right\rbrace$. The number of vertices (and edges) of the cell is labeled with $N_{\cv_{i},e}$, and the outward pointing normal vector of edge $e$ writes $\nn_{e}$. The characteristic mesh size of each cell is evaluated as 
\begin{equation}
	\cvs_i=\frac{2\left|\cv_{i} \right|}{\sum\limits_{e=1}^{N_{\cv_{i},e}}\left|\cvb_{i,e}\right|},
	\label{eqn.h}
\end{equation}
with $\left|\cvb_{i,e}\right|$ being the length of edge $e$ of cell $P_i$. 

Assuming the existence of a constant $\varrho>0$ common for all elements, we require each cell of the computational mesh to verify the following regularity assumptions:
\begin{itemize}
	\item $\cv_{i}$ is star-shaped with respect to a disk with radius
	$r>\varrho \max \limits_{i=1,\dots,N} \cvs_{i}$. Hence, the elements are simply connected subsets of $\mathbb{R}^{2}$ with a finite number of vertices and edges.
	\item Each edge $e \in \cvb_{i}$ verifies $|e| \geq\varrho\cvs_{i}$, so that the number of edges of any element is limited over the whole mesh.
\end{itemize}

\begin{remark}[Numerical integration]
In the numerical scheme that will be presented in the sequel, Gauss-Lobatto formulae \cite{stroud} are used to compute integrals over the boundary $\cvb_{i}$, while for volume integrals we adopt the efficient numerical integration proposed in \cite{SOMMARIVA2009886,SOMMARIVA2020} for arbitrary shaped polygonal cells.
\end{remark}

\subsection{Convective sub-system: Finite Volume scheme}\label{sec.convective}
The contribution of the convective terms is computed using an explicit finite volume method on the polygonal grid. Data are stored as integral cell averages within each control volume, that is
\begin{equation}
	q_{i}  = \frac{1}{\cvv{i}}\int\limits_{\cv_{i}} q\left(\xx\right)  \dV,
	\label{eqn.cellAv}
\end{equation}
for a generic quantity $q(\xx,t)$. Let us consider only the explicit contributions in the semi-discrete scheme \eqref{eqn.mass_sd}-\eqref{eqn.press_sd}. Recalling the definitions of the intermediate momentum \eqref{eqn.intermediatemomentum} and kinetic energy \eqref{eqn.pressure_tdsv}, integration over the control volume $\cv_{i}$ with subsequent application of Gauss theorem yields
\begin{subequations}
	\label{eqn.conective_fv}
	\begin{align}
	\rho_i^{n+1} &= \rho_i^n -  \frac{\dt}{\cvv{i}} \sum\limits_{e=1}^
	{N_{\cv_{i},e}} \int\limits_{\cvb_{i},e} \Flux_{\rho} \cdot \nn_{e} \dS, \label{eqn.rho_fv} \\
	\wsstar_{i}  &= \ww^{n}_{i}  - \frac{\dt}{\cvv{i}} \sum\limits_{e=1}^
	{N_{\cv_{i},e}} \int\limits_{\cvb_{i},e} \Flux_{\ww} \cdot \nn_{e} \dS, \label{eqn.mom_fv} \\ 
	\psstar_i &= \rho^{n}\ke^{n} - \frac{\dt}{\cvv{i}} \sum\limits_{e=1}^
	{N_{\cv_{i},e}} \int\limits_{\cvb_{i},e} \Flux_{\ke} \cdot \nn_{e} \dS. \label{eqn.ek_fv}
\end{align}
\end{subequations}
The flux contribution across each edge $e$ of cell $\cvb_{i,e}$ is computed using the Rusanov numerical flux function: 
\begin{subequations}
		\label{eqn.rusanov}
	\begin{align}
	\mathcal{F}_{\rho} \cdot \nn_e &= \frac{1}{2} \left( \ww_e^{+,n} + \ww_e^{-,n} \right) \cdot \nn_{e} - \frac{1}{2} \left|s_{\max}\right|_e \left( \rho_e^{+,n} - \rho_e^{-,n} \right), \\
	\mathcal{F}_{\ww} \cdot \nn_e &= \frac{1}{2} \left( \frac{\ww_e^{+,n}\otimes\ww_e^{+,n} }{\rho_e^{+,n}}+ \frac{\ww_e^{-,n}\otimes\ww_e^{-,n}}{\rho_e^{-,n}} \right) \cdot \nn_{e} - \frac{1}{2} \left|s_{\max}\right|_e \left( \ww_e^{+,n} - \ww_e^{-,n} \right), \\
	\mathcal{F}_{\ke} \cdot \nn_e &= \frac{1}{2} \left( \ke_e^{+,n} \ww_e^{+,n} + \hf_e^{+,n} + \ke_e^{-,n} \ww_e^{-,n} + \hf_e^{-,n} \right) \cdot \nn_{e} - \frac{1}{2} \left|s_{\max}\right|_e \left( \rho_e^{+,n} \ke_e^{+,n}  -  \rho_e^{-,n} \ke_e^{-,n} \right),
	\end{align}
\end{subequations}
where the superscripts $(+,-)$ denote the right and left states to the edge, respectively. Since the numerical dissipation coefficient corresponds to the maximum absolute eigenvalue of the convective sub-system associated to the left and right states, 
\begin{equation}
	\left|s_{\max}\right|_e:= \max \left\lbrace \left| \vel^{+,n}_{e}\cdot\nn_{e} \right|,\left| \vel^{-,n}_{e}\cdot\nn_{e} \right| \right\rbrace,\label{eqn.alpharusanov}
\end{equation}
it does not depend on the sound speed. Consequently, the resulting scheme is well suited for the solution of low Mach number flows because the numerical dissipation is vanishing in the incompressible limit of the model.

Regarding the computation of the left and right states in \eqref{eqn.rusanov}, we may simply consider the cell average values defined by \eqref{eqn.cellAv}. However, this approach would only lead to a first order accurate scheme in space. To increase the order of accuracy, a polynomial reconstruction of the numerical solution can be performed and used in the computation of the numerical fluxes \eqref{eqn.rusanov}. In particular, we consider a CWENO reconstruction on general polygonal meshes following the procedure detailed in \cite{CWENOGBK,FVBoltz,BTC23}.

\subsubsection{High order extension in space: CWENO reconstruction}\label{sec.cweno}
Let us consider a generic variable $q(\xx)$ and a generic control volume $\cv$. Central WENO (CWENO) schemes \cite{LPR99,LPR00} are based on the computation of a reconstruction polynomial $w(\xx)$ of arbitrary degree $\kk$ at each cell $\cv$ of the quantity $q(\xx)$, starting from the known cell averages \eqref{eqn.cellAv}. The reconstruction polynomial is expressed as 
\begin{equation}
	\label{eqn.recPoly}
	w(\xx) = \sum \limits_{\ell=1}^{n_{\kk}} \beta_{\ell}(\xx) \, \hat{w}_{\ell},
\end{equation} 
where $\hat{w}_{\ell}$ denote the degrees of freedom and $\beta_{\ell}(\xx)$ represent the conservative Taylor basis functions that, according to \cite{FVBoltz}, write
\begin{equation}
	\label{eqn.Voronoi_modal}
	\beta_{\ell}(\xx)_{|_{\cv}} = m_{\bm{\kappa}} - \frac{\eta}{\cvv{\,}}\int\limits_{\cv} m_{\bm{\kappa}} \, \dV,
\end{equation}
where $m_{\bm{\kappa}}=\left(\frac{\xx-\xx_{\cv}}{\cvs_{\cv}}\right)^{\bm{\kappa}}$ are the scaled monomials of degree 
$\left|\bm{\kappa}\right|=\kappa_1+\kappa_2$ with $\xx^{\kappa}=\left(x^{\kappa_1},y^{\kappa_2}\right)$.
The coefficient $\eta$, 
\begin{equation}
	\eta=\left\lbrace \begin{array}{ll}
		0 & |\bm{\kappa}|=0\\
		1 & |\bm{\kappa}|>0
	\end{array} \right. ,
\end{equation}
ensures the conservation property of the basis functions, i.e. they verify 
\begin{equation}
	\label{eqn.modal_cons}
	\frac{1}{|\cv|}\int\limits_{\cv} \sum_{\ell=1}^{n_{\kk}}\beta_{\ell}(\xx) \, \dV = 1,
\end{equation}
thus $\hat{w}_{1}$ corresponds to the cell averaged value of $q(\xx)$ at cell $\cv$ given by \eqref{eqn.cellAv}. The number of degrees of freedom $n_{\kk}$ depends on the polynomial degree, and for $d=2$ it is explicitly given by
\begin{equation}
	n_{\kk} = \frac{(\kk+1)(\kk+2)}{2}.
	\label{eqn.nk}
\end{equation}
Once the reconstructed polynomials \eqref{eqn.recPoly} are obtained in each cell, the high order extrapolated states at each edge are evaluated and substituted into the numerical flux function \eqref{eqn.rusanov} leading to a numerical scheme of order $\kk+1$.  
For further details related to the CWENO reconstruction procedure on arbitrary polygonal grids we refer to Appendix~B in \cite{BTC23}.

\subsection{Viscous and pressure sub-systems: Virtual Element Method}\label{sec.viscouspressure}
The viscous and pressure stages are carried out relying on a virtual element method. To this end, the variational formulations of the sub-systems \eqref{eqn.intermediatemomentumv}, \eqref{eqn.pressure_sd2} and \eqref{eqn.pressure_sd3}
are first obtained. Then, we introduce the virtual element space and the related projector operators. Finally, the  virtual element discretization of the weak problems is presented.

\subsubsection{Variational formulation of the viscous sub-system} \label{ssec.weak}
Multiplication of \eqref{eqn.intermediatemomentumv} by a test function $\tf\in \HONEO(\Omega)=\left\lbrace\tf\in \HONE\left(\Omega\right) \mid \int\limits_{\Omega} \tf \dV = 0 \right\rbrace$, integration over the computational domain and use of Green's formula for the viscous term yield the following weak problem:
\begin{weakproblem} Find $\vstar\in \HONEO(\Omega)$ such that
\begin{equation}
	\int\limits_{\Omega} \rho^{n+1} \, \vstar \tf \dV + \, \dt \int\limits_{\Omega} \btaus \cdot \gra \tf \dV= \int\limits_{\Omega} \wsstar \tf \dV + \, \dt \int\limits_{\partial\Omega} \btaus \cdot \nn \, \tf \dS, \qquad \forall \tf\in \HONEO(\Omega). \label{eqn.wfvsicous}
\end{equation}
\end{weakproblem}

\subsubsection{Variational formulation of the pressure sub-system}
Multiplication of \eqref{eqn.pressure_sd2} by a test function $\tf\in \HONEO(\Omega)$ and integration over the computational domain leads to
\begin{gather}
	\frac{1}{ \gamma-1 }\int\limits_{\Omega}^{} \pdstar{n+1} \tf	\dV
	- \dt^2 \int\limits_{\Omega}^{}\dive \left(\ent^n  \gra \pdstar{n+1} \right) \tf \dV 
	= \frac{1}{ \gamma-1 } \int\limits_{\Omega}^{}\press^{n} \tf \dV
	+ \int\limits_{\Omega}^{} \pstar \tf \dV \notag\\
	- \halb \int\limits_{\Omega}^{} \frac{1}{\rho^{n+1}}\wstar\cdot \wstar	\tf \dV
	- \dt \int\limits_{\Omega}^{}  \dive \left(\ent^n  \wstar \right) \tf \dV. \label{eqn.pressure_vf1}
\end{gather}
Using Green's formulas for the enthalpy dependent terms, and the definition of the new momentum $\ww^{n+1}$ given by \eqref{eqn.velupdate}, yield the weak problem for the \textit{pressure flux} $\pdstar{n+1} $:
\begin{weakproblem} Find $\pdstar{n+1}\in H^{1}_{0}(\Omega)$ such that
	\begin{gather}
		\frac{1}{ \gamma-1 }\int\limits_{\Omega}^{} \pdstar{n+1} \tf	\dV
		+ \dt^2 \int\limits_{\Omega}^{} \left(\ent^n  \gra \pdstar{n+1} \right) \cdot \gra \tf \dV 
		=  	\frac{1}{ \gamma-1 } \int\limits_{\Omega}^{}\press^{n} \tf \dV+ \int\limits_{\Omega}^{} \pstar \tf \dV\notag\\	 
		- \halb \int\limits_{\Omega}^{} \frac{1}{\rho^{n+1}}\wstar\cdot \wstar	\tf \dV
		+ \dt \int\limits_{\Omega}^{}   \left(\ent^n  \wstar \right)\cdot \gra \tf \dV
		- \dt \int\limits_{\partial\Omega}^{}   \left(\ent^n  \ww^{n+1} \right)\cdot\nn\, \tf \dS, \qquad \forall \tf\in \HONEO(\Omega). \label{eqn.pressure_vf3}
	\end{gather}
\end{weakproblem}
\noindent Similarly, the weak problem for the \textit{pressure state} $\press^{n+1}$ given by \eqref{eqn.pressure_sd3} results as follows:
\begin{weakproblem} Find $\press^{n+1}\in \HONEO(\Omega)$ such that
	\begin{gather}
		\frac{1}{ \gamma-1 }\int\limits_{\Omega}^{} \press^{n+1} \tf	\dV
		+ \dt^2 \int\limits_{\Omega}^{} \left(\ent^n  \gra \press^{n+1} \right) \cdot \gra \tf \dV 
		=  	\frac{1}{ \gamma-1 } \int\limits_{\Omega}^{}\press^{n} \tf \dV 
		+ \int\limits_{\Omega}^{} \pstar \tf \dV	 
		\notag\\- \halb \int\limits_{\Omega}^{} \frac{1}{\rho^{n+1}}\ww^{n+1}\cdot \ww^{n+1}	\tf \dV
		+ \dt \int\limits_{\Omega}^{}   \left(\ent^n  \wstar \right)\cdot \gra \tf \dV
		- \dt \int\limits_{\partial\Omega}^{}   \left(\ent^n  \ww^{n+1} \right)\cdot\nn\, \tf \dS, \qquad \forall \tf\in \HONEO(\Omega). \label{eqn.pressure_vf4}
	\end{gather}
\end{weakproblem}
Since \eqref{eqn.pressure_vf3} and \eqref{eqn.pressure_vf4} involve the same integral operators and exhibit the same mathematical structure, in what follows we simply introduce the discretization for \eqref{eqn.pressure_vf4} which corresponds to the weak problem of \eqref{eqn.pressure_sd}. Indeed, replacing $\ww^{n+1}$ by $\wstar$ in the third term of the right hand side of the above weak problem leads to the discretization of \eqref{eqn.pressure_vf3}.

\subsubsection{The virtual element space} \label{ssec.lvem}
The conforming virtual element space of order $\kk$ is defined on a polygonal element $\cv$ as
\begin{equation}\label{eq:Vh:def}
	V^h_{\kk}(\cv) = \left\lbrace
	\vsh{}\in\HONE(\cv)\,\mid \, \lapla v \in\PS{\kk-2}(\cv) ,
	\vsh{\mid \cvb}\in\CS{0}(\cvb), \
	\vsh{\mid e}\in\PS{\kk}(e)\,\,\forall e \in \cvb \, \right\rbrace,
\end{equation}
with $\PS{\kk}(\cv)=\left<\{ m_{\bm{\kappa}} \mid 0 \leq \left|\bm{\kappa}\right| \leq \kk \}\right>$ representing the space of polynomials of degree lower or equal to $\kk$ that can be generated by the basis of scaled monomials
\begin{equation}
\mathcal{M}_{\kk}(\cv) =	\{ m_{\bm{\kappa}} \mid 0 \leq \left|\bm{\kappa}\right| \leq \kk \}. \label{eqn.Pkbasis}
\end{equation} 
Consequently,  $\PS{\kk}(\cv)$ is a subset of $V^h_{\kk}$  of dimension $n_{\kk} := \textrm{dim}\,  \PS{\kk}(\cv)$, with $n_{\kk}$ given by \eqref{eqn.nk}. As presented in \cite{BBCMMR13}, each function $v\in V^h_{\kk}(\cv)$ is uniquely determined by the following degrees of freedom:
\begin{itemize}
	\item the $\NE$ values of $\vsh{}$ at the vertices of $\cv$;
	\item the $(\kk-1)\NE$ values of $\vsh{}$ on the $(\kk+1)$ Gauss-Lobatto quadrature points on each edge $e$;  
	\item the moments up to degree $\kk-2$ of $\vsh{}$ in $\cv$, defined as
	\begin{equation}
		\label{eqn.dof_mom}
		\frac{1}{\cvv{}} \int\limits_{\cv} \vsh{} m_{\alpha} \, \dV, \qquad \alpha \in\left\lbrace 1,\dots,n_{\kk-2}\right\rbrace.
	\end{equation}
\end{itemize}
Hence, each $\vsh{} \in V^h_{\kk}(\cv)$ has a total number of degrees of freedom
\begin{equation}\label{eq:dim_Vh} 
	\NDOF := \textrm{dim} \, V^h_{\kk}(\cv) =  \kk \, \NE + \frac{(\kk-1)\,\kk}{2},
\end{equation}
and can be written using the Lagrange interpolation as 
\begin{align}\label{eq:vh_repr} 
	\vsh{} = \sum_{i=1}^{\NDOF} \varphi_i \hat{\vsh{}}_{i} ,
\end{align}
where $\hat{\vsh{}}_{i}:=\dof_i(\vsh{})$ 
denotes the value of the $i$-th degree of freedom and $\{\varphi_{i}\}_{i=1}^{\NDOF}$ is a basis of $V^h_{\kk}(\cv)$ verifying
\begin{equation}
	\dof_i(\varphi_j) = \delta_{ij}, \qquad i,j\in\left\lbrace 1, \ldots, \NDOF\right\rbrace.
\end{equation}
Finally,  gathering together the elemental spaces $V^h_{\kk}(\cv_{i})$ for all $\cv_{i}\in \mathcal{T}_{\Omega}$, we get the global conforming virtual element space associated to the tessellation:
\begin{equation}\label{eq:Vh_global:def}
	V^h_{\kk} = \left\lbrace
	\vsh{} \in\HONE(\cv)\,\mid\, \vsh{\mid\cv}\in V^h_{\kk}(\cv)\,\,\forall\cv\in\mathcal{T}_{\Omega} \, \right\rbrace.
\end{equation}

\subsubsection{Elliptic projection operator}
A key feature of virtual element spaces is that the explicit expression of the virtual basis functions $\varphi_i$ is not known but only their values at the degrees of freedom are available. Thus, to discretize the weak problems, we make use of an elliptic projector operator mapping functions from the virtual element space to the corresponding polynomial space, namely $\proj:V^h_{\kk}(\cv) \rightarrow \PS{\kk}(\cv)$. Following \cite{BBCMMR13}, the orthogonality condition
\begin{equation}\label{eq:vem_proj} 
	\int\limits_{\cv} \nabla p_{\kk} \cdot \nabla (\proj\vsh{} - \vsh{}) \, \dV = 0, \quad \forall p_{\kk} \in \PS{\kk}(\cv),
\end{equation}
allows the definition of this projector operator up to a constant which is then determined at the aid of an additional operator $P_0: V^h_{\kk}(\cv) \rightarrow \PS{0}(\cv)$ verifying
\begin{align}\label{eq:vem_proj_const} 
	P_0(\proj\vsh{} - \vsh{}) = 0.
\end{align}
Among the different possibilities, according to \cite{vem2013} we select the $P_0$ operator to be given by
\begin{equation}\label{eq:vem_proj_0} 
	P_0 \vsh{} := \left\lbrace 
	\begin{array}{lr}
		\displaystyle \frac{1}{\NE} \sum\limits_{i=1}^{\NE} \vsh{}(\xx_i) & \textrm{if} \enskip \kk = 1, \\[10pt]
		\displaystyle \frac{1}{\abs{\cv}}\int\limits_{\cv} \vsh \, \dV & \textrm{if} \enskip \kk \geq 2. 
	\end{array}	\right. 
\end{equation}
Consequently, the projection $\proj\vsh{}$ can be computed simply using the degrees of freedom $\hat{\vsh{}}_{i}$. Furthermore, since $ p_\kk \in \PS{\kk}(\cv)$ and $\proj\vsh{} \in \PS{\kk}(\cv)$, we can express \eqref{eq:vem_proj} as an algebraic system in terms of the basis of monomials $\mathcal{M}_\kk(\cv)$ given by \eqref{eqn.Pkbasis} as
\begin{align}\label{eq:vem_proj2} 
	\sum_{\beta=1}^{n_{\kk}} s^\beta \int\limits_{\cv} \nabla m_\alpha \cdot \nabla m_\beta \, \dV = \int\limits_{\cv} \nabla m_\alpha \cdot \nabla \vsh \, \dV, \qquad  \alpha \in\left\lbrace 1,\dots,n_{\kk}\right\rbrace,
\end{align}
where we have taken into account that
\begin{equation}\label{eqn.projinpk}
	\proj\vsh{}= \sum_{\beta=1}^{n_{\kk}} s^\beta m_\beta .
\end{equation}
Then, applying integration by parts to the right hand side of \eqref{eq:vem_proj2}, we get
\begin{align}\label{eq:vem_proj3} 
	\sum_{\beta=1}^{n_{\kk}} s^\beta \int\limits_{\cv} \nabla m_\alpha \cdot \nabla m_\beta \, \dV = -\int\limits_{\cv}\Delta m_\alpha \vsh \, \dV + \int\limits_{\cvb} \frac{\partial m_\alpha}{\partial n} \vsh \, \dS, \qquad  \alpha \in\left\lbrace 1,\dots,n_{\kk}\right\rbrace .
\end{align}
The left hand-side of this system involves the integration of known polynomials over $\cv$ while the right hand side accounts of a first integral that can be explicitly computed employing the internal degrees of freedom of $\vsh{}$ and a second integrand that can be exactly calculated using the Gauss-Lobatto quadrature points along the edges defining $\cvb$, which coincides with the known boundary degrees of freedom of the virtual basis.
The above system is finally supplemented with the following condition that directly comes from \eqref{eq:vem_proj_0}:
\begin{equation}\label{eq:vem_proj_aux} 
	\sum_{\beta=1}^{n_{\kk}} s^\beta P_0 m_\beta = P_0 \vsh.
\end{equation}
Gathering \eqref{eq:vem_proj3} and \eqref{eq:vem_proj_aux}, we have a linear system in the unknowns $s^{\beta}$. Once the solution is computed, $s^{\beta}$ is substituted into \eqref{eqn.projinpk} to evaluate the projection of all basis functions $\varphi_{i}$ as
\begin{equation} 
	\proj\varphi_i = \sum_{\alpha=1}^{n_{\kk}} s^\alpha_i m_\alpha, \quad i \in\left\lbrace 1,\dots,\NDOF\right\rbrace.
\end{equation}
Defining the matrices $\mathbf{G}\in \mathbb{M}_{n_{\kk} \times n_{\kk}}$,  $\mathbf{B}\in \mathbb{M}_{n_{\kk} \times \NDOF}$ with non-zero entries
\begin{align}\label{gb_matrices} 
	(\mathbf{G})_{\alpha \beta} &= P_0 m_\beta, \quad &  \alpha = 1, &\enskip \beta \in\left\lbrace 1,\dots,n_{\kk} \right\rbrace, \\
	(\mathbf{G})_{\alpha \beta} &= \int\limits_{\cv}\nabla m_\alpha \cdot \nabla m_\beta \, \dV, \quad & \alpha \geq 2, & \enskip \beta \in\left\lbrace 1,\dots,n_{\kk} \right\rbrace, \\
	(\mathbf{B})_{\alpha  i} &= P_0 \varphi_{i}, \quad &  \enskip \alpha = 1, & \enskip i  \in\left\lbrace 1,\dots,\NDOF \right\rbrace, \\
	(\mathbf{B})_{\alpha  i} &= \int\limits_{\cv} \nabla m_\alpha \cdot \nabla \varphi_i \, \dV, \quad & \alpha \geq 2, & \enskip i \in\left\lbrace 1,\dots,\NDOF \right\rbrace,
\end{align}
system \eqref{eq:vem_proj3}-\eqref{eq:vem_proj_aux} compactly writes
\begin{align}\label{eq:proj_system} 
	\mathbf{G} \accentset{\star}{\mathbf{\Pi}}^{\nabla}_k = \mathbf{B},
\end{align}
with $\accentset{\star}{\mathbf{\Pi}}^{\nabla}_k \in\mathbb{M}_{n_{\kk} \times \NDOF}$ the matrix representation of the projector operator $\proj$ in the basis $\mathcal{M}_k(\cv)$, i.e. $(\accentset{\star}{\mathbf{\Pi}}^{\nabla}_k)_{\alpha i} = s_i^\alpha$. 
Moreover, introducing the following matrix for a change of basis  
\begin{align}\label{D_matrix} 
	(\mathbf{D})_{i \alpha} = \dof_i (m_\alpha), \qquad i \in\left\lbrace 1,\dots,\NDOF \right\rbrace, \quad \alpha \in\left\lbrace 1,\dots,n_{\kk} \right\rbrace,
\end{align}
we can also represent the operator $\proj$ in the canonical basis of $V^{h}_{\kk}(\cv)$ as 
\begin{align}\label{Pi_nabla} 
	\mathbf{\Pi}^\nabla_{\kk} = \mathbf{D}  \accentset{\star}{\mathbf{\Pi}}^{\nabla}_{\kk} = \mathbf{D}\mathbf{G}^{-1}\mathbf{B},\qquad \mathbf{\Pi}^\nabla_{\kk}\in\mathbb{M}_{\NDOF \times \NDOF}.
\end{align}

\subsubsection{$L_{2}$ projection operator}
Similarly to what has been done for the elliptic projector operator, we also define the $L_2$ projection operator $\projL: V^h_{\kk}(\cv) \rightarrow \PS{\kk}(\cv)$ employing the orthogonality condition
\begin{align}\label{eq:vem_projL1} 
	\int\limits_{\cv} p_{\kk} (\projL\vsh{} - \vsh{})  \dV = 0, \qquad \forall p_{\kk} \in \PS{\kk}(\cv).
\end{align}
The projector $\projL \vsh{}$ can be computed by considering the degrees of freedom $\hat{\vsh{}}_{i}$ and, since $\projL\vsh{} \in \PS{\kk}(\cv)$, it can be represented in the basis $\mathcal{M}_k(\cv)$ as
\begin{align}\label{eq:vem_projL2} 
	\projL\vsh{} = \sum_{\beta=1}^{n_{\kk}} r^\beta m_\beta.
\end{align}
Using the above definition in \eqref{eq:vem_projL1}, a linear system for the $n_{\kk}$ unknowns $r^\beta$ is obtained:
\begin{align}\label{eq:vem_projL3} 
	\sum_{\beta=1}^{n_{\kk}} r^\beta \int\limits_{\cv} m_\alpha m_\beta \dV = \int\limits_{\cv} m_\alpha \vsh{} \ \dV, \quad \alpha \in\left\lbrace 1,\dots,n_{\kk}\right\rbrace. 
\end{align}
As for the elliptic projector system, the left-hand side term of \eqref{eq:vem_projL3} consists in the integral of known polynomials over $\cv$ and can be easily computed. Regarding the right-hand side, we know moments as degrees of freedom only for $m_\alpha \in \PS{\kk-2}(\cv)$ so we replace $\vsh{}$ with its elliptic projection $\proj \vsh{}$ for monomials of degree $\kk$ and $\kk-1$. Consequently, we get the system
\begin{align}\label{eq:proj_system2} 
	\mathbf{H} \accentset{\star}{\mathbf{\Pi}}^0_\kk = \mathbf{C},
\end{align}
with matrices
\begin{subequations}
	\label{eqn.HCmatrix}
	\begin{align} 
		(\mathbf{H})_{\alpha \beta} &= \int\limits_{\cv}m_\alpha m_\beta  \dV, &  \alpha,\beta \, \in\left\lbrace 1,\dots,n_{\kk} \right\rbrace, & \label{h_matrix} \\
		(\mathbf{C})_{\alpha  i} &= 
		\int\limits_{\cv}m_\alpha \varphi_i \dV, &  \alpha\in\left\lbrace 1,\dots,n_{\kk-2} \right\rbrace, \quad i \in\left\lbrace 1,\dots,\NDOF \right\rbrace, \\
		(\mathbf{C})_{\alpha  i} &= 
		\int\limits_{\cv}m_\alpha \proj \varphi_i \dV, &  \alpha\in\left\lbrace n_{\kk-2}+1,\dots,n_{\kk} \right\rbrace, \quad i \in\left\lbrace 1,\dots,\NDOF \right\rbrace,
	\end{align}
\end{subequations}
$\mathbf{H}\in\mathbb{M}_{n_{\kk} \times n_{\kk}}$, $\mathbf{C}\in\mathbb{M}_{n_{\kk} \times \NDOF}$ and $\accentset{\star}{\mathbf{\Pi}}^0_\kk\in\mathbb{M}_{n_{\kk} \times \NDOF}$ the matrix representation of $\projL$. Making use of the change of basis \eqref{D_matrix}, we can also express the $L_{2}$ operator in terms of the canonical basis of $V^h_{\kk}(\cv)$ as
\begin{align}\label{Pi_0} 
	\mathbf{{\Pi}^{0}_{\kk}} = \mathbf{D} \accentset{\star}{\mathbf{\Pi}}^0_{\kk} = \mathbf{D} \mathbf{H}^{-1} \mathbf{C}.
\end{align}

Finally, denoting by $\mathbf{H}'\in\mathbb{M}_{n_{\kk-1} \times n_{\kk-1}}$ the sub-matrix of $\mathbf{H}$ formed by its first $n_{\kk-1}$ rows and columns, and by $\mathbf{C}'\in\mathbb{M}_{n_{\kk-1} \times \NDOF}$ the sub-matrix containing the first $n_{\kk-1}$ rows of $\mathbf{C}$, we get the linear system
\begin{align}\label{eq:proj_system3} 
	\mathbf{H}' \accentset{\star}{\mathbf{\Pi}}^0_{\kk-1} = \mathbf{C}',
\end{align}
that defines the $L_2$ projection onto $\PS{\kk-1}(\cv)$. The related matrix representation of $\projLm$ in terms of the canonical basis of $V^h_{\kk-1}(\cv)$ reads
\begin{align}\label{eq:proj_system4} 
	\mathbf{\Pi}^0_{\kk-1} = \mathbf{D}'\accentset{\star}{\mathbf{\Pi}}^0_{\kk-1},
\end{align}
where 
$\mathbf{D}'\in\mathbb{M}_{\NDOF \times n_{\kk-1}}$ is given by the first $n_{\kk-1}$ columns of $\mathbf{D}$.

\subsubsection{Discretization of the viscous sub-system}
The discrete intermediate velocity $\vstar$ as well as the viscous stress tensor $\btaus$ in \eqref{eqn.wfvsicous} are approximated by means of an expansion in the local virtual element space $V^{h}_{\kk}\left(\cv\right)$ using the virtual basis functions within each control volume $\cv\in \mathcal{T}_{\Omega}$ as
\begin{equation}
	\vstar\left(\xx \right)_{\mid \xx\in\cv}  = \sum\limits_{i=1}^{\NDOF} \varphi_{i} \dofv{i}{\ast}, \qquad 
	\btaus \left(\xx \right)_{\mid \xx\in\cv}  = \sum\limits_{i=1}^{\NDOF} \varphi_{i} \doftau{i}{\ast}. \label{eqn.momvem}
\end{equation}
Choosing $\tf =\varphi_{j}$ and substituting \eqref{eqn.momvem} into the weak problem \eqref{eqn.wfvsicous}, we obtain the local discrete viscous sub-system for the element $\cv$:
\begin{equation}
	\int\limits_{\cv} \rho^{n+1} \, \varphi_{i} \varphi_{j} \dV \, \dofv{i}{\ast} + \dt \int\limits_{\cv} \varphi_{i} \cdot \gra \varphi_{j} \dV \, \doftau{i}{\ast} = \left( \mathbf{F}^{n}_{\vel,\cv}\right)_{j}  , \label{eqn.wfvsicous_disc}
\end{equation}
with the assumption $ \cvb \notin \partial \Omega$ and the right hand side term given by
\begin{equation}
	\left( \mathbf{F}^{n}_{\vel,\cv}\right)_{j}  = \int\limits_{\cv} \wsstar \varphi_{j} \dV.
\end{equation}
By introducing the following matrix definitions
\begin{align} 
	(\mathbf{M}_{\cv}^{\, \rho})_{i,j} &= \int\limits_{\cv} \rho^{n+1}\varphi_i \varphi_j  \, \dV,  \label{mass_mat_weigh}\\\
	(\mathbf{K}_{\vel,\cv})_{i,j} &= \int\limits_{\cv} \varphi_{i}\cdot \gra \varphi_{j} \dV, \label{stiffness} 
\end{align}
we can rewrite the local system \eqref{eqn.wfvsicous_disc} in matrix form:
\begin{equation}\label{eqn.discrweakform_visc} 
	\mathbf{M}_{\cv}^{\, \rho} \, \dofv{\cv}{\ast} + \Delta t \, \mathbf{K}_{\vel,\cv} \, \doftau{\cv}{\ast} = \mathbf{F}^{n}_{\vel,\cv},
\end{equation}
where the vectors $\dofv{\cv}{\ast}$ and $\doftau{\cv}{\ast}$ collect the degrees of freedom of the expansions \eqref{eqn.momvem} for cell $\cv$.

Since $\mathbf{M}_{\cv}^{\, \rho}$, $\mathbf{K}_{\vel,\cv}$ and $\mathbf{F}^{n}_{\vel,\cv}$ involve the integration of functions with unknown explicit expression in the interior of $\cv$, their approximation is done employing  the elliptic and $L_2$ projection operators introduced in the previous sections. Let us first introduce the basis functions expansion
\begin{align}
	\varphi_i = \projL \varphi_i + (\Id-\projL) \varphi_i.
\end{align}
Since the mass matrix \eqref{mass_mat_weigh} is weighted by the new density value over the cell, then it can be computed by means of the above expansion as the sum of two terms:
\begin{equation}
\int\limits_{\cv}\rho^{n+1} \varphi_i  \varphi_j  \dV = \int\limits_{\cv}\rho^{n+1} \projL\varphi_i  \projL \varphi_j  \dV + \,  \mathcal{S}_{\cv}\left( \, (\Id-\projL)\varphi_i, \,   (\Id-\projL) \varphi_j \, \right).
\end{equation}
The first term guarantees consistency and can be exactly computed. The second term $\mathcal{S}_{\cv}$, which ensures stability, is approximated using the so called $\textrm{dof}_i-\textrm{dof}_i$ stabilization \cite{BBCMMR13}:
\begin{align}\label{proj_M3} 
	\mathcal{S}_{\cv}\left( \, (\Id-\projL)\varphi_i, \,   (\Id-\projL) \varphi_j \, \right) := 
	\cvv{} \bar{\rho}_{\cv}^{n+1} \sum_{r=1}^{\NDOF} \dof_r((\Id-\projL)\varphi_i) \dof_r((\Id-\projL) \varphi_j),
\end{align}
where $\bar{\rho}_{\cv}^{n+1}$ is the averaged density in the cell that has already been computed by the explicit finite volume schemes at the convection stage.
Further alternatives for the approximation of this term can be found in \cite{VEM2,Mascotto2018}. 
Recalling the matrix representation of the $L_2$ projector \eqref{Pi_0}, the approximation of the mass matrix $\mathbf{M}_{\cv}^{\, \rho}$ is
expressed in matrix form as
\begin{align}\label{proj_M4} 
	\mathbf{M}_{\cv}^{\, \rho} = \left(\mathbf{C}^{\, \rho}\right) ^{\top}\left(\mathbf{H}^{\, \rho}\right)^{-1}\mathbf{C}^{\rho} +  \cvv{}\bar{\rho}_{\cv}^{n+1}(\mathbf{I}-\projLb)^\top(\mathbf{I}-\projLb)
\end{align}
where the matrix $\mathbf{H}^{\, \rho}$ is defined as 
\begin{equation}
		\left(\mathbf{H}^{\, \rho}\right)_{\alpha,\beta} := \int\limits_{\cv} \rho^{n+1} m_{\alpha} m_{\beta} \, \dV, \qquad \alpha,\beta \, \in\left\lbrace 1,\dots,n_{\kk} \right\rbrace,
\end{equation}
and the matrix $\mathbf{C}^{\, \rho}$ is readily obtained by its definition \eqref{eq:proj_system2} using $\mathbf{H}^{\, \rho}$ instead of $\mathbf{H}$.

Regarding the matrix \eqref{stiffness}, 
we use the $L_{2}$ projector operator $\projLm$ on polynomials of degree $\kk-1$ to approximate gradients of the virtual element basis functions. This is needed to avoid loss of optimal convergence, as studied in \cite{BBMR16}. More precisely, let us define the matrices  
\begin{align*}
	(\mathbf{E}^x)_{i\alpha} = \int\limits_{\cv} \varphi_{i,x}m_\alpha \, \dV, \quad\quad	(\mathbf{E}^y)_{i\alpha} = \int\limits_{\cv} \varphi_{i,y}m_\alpha \, \dV, \quad \quad \alpha \in \left\lbrace 1,\dots,n_{\kk-1}\right\rbrace,
\end{align*}
and let $\hat{\mathbf{H}}$ be the matrix formed by the first $n_{\kk-1}$ rows of $\mathbf{H}$ \eqref{h_matrix}. The matrix representation of the gradient projectors is then evaluated as
\begin{equation}
	\accentset{\star}{\mathbf{\Pi}}^{0,x}_{\kk-1}= \hat{\mathbf{H}}^{-1}\mathbf{E}^x, \qquad \accentset{\star}{\mathbf{\Pi}}^{0,y}_{\kk-1}= \hat{\mathbf{H}}^{-1}\mathbf{E}^y.
\end{equation}
Consequently, the term \eqref{stiffness} of the viscous sub-system can be divided into two contributions accounting for each spatial direction, that is
\begin{equation}
	\mathbf{K}_{\vel,\cv}^{x} = \left(\accentset{\star}{\mathbf{\Pi}}^{0,x}_{\kk-1}\right)^\top \mathbf{H} \accentset{\star}{\mathbf{\Pi}}^{0}_{\kk} , \qquad \mathbf{K}_{\vel,\cv}^{y} = \left(\accentset{\star}{\mathbf{\Pi}}^{0,y}_{\kk-1}\right)^\top \mathbf{H} \accentset{\star}{\mathbf{\Pi}}^{0}_{\kk}.
\end{equation}
Gathering the above expressions in \eqref{eqn.discrweakform_visc}, we get
\begin{equation}
	\mathbf{K}_{\vel,\cv}\doftau{\cv}{\ast}  =  \mathbf{K}_{\vel,\cv}^{x}\doftau{\cv}{\ast,x} + \mathbf{K}_{\vel,\cv}^{y}\doftau{\cv}{\ast,y},
\end{equation}
with $\doftau{\cv}{\ast}=\{ \doftau{\cv}{\ast,x},\doftau{\cv}{\ast,y}\}$ representing the degrees of freedom associated to the local virtual element expansion of the stress tensor along each spatial direction. To compute $\doftau{\cv}{\ast}$, we resort to an $L_2$ projection onto the local virtual element space:
\begin{equation}
	\doftau{\cv}{\ast} = \left( \mathbf{M}_{\cv} \right)^{-1} \int\limits_{\cv} \varphi \, \btaus  \, \dV,
\end{equation}
where the mass matrix $\mathbf{M}_{\cv}$ is simply evaluated as \eqref{proj_M4} by assuming unit density weight (see Eqn.\eqref{proj_M5} in the next section). The viscous stress is computed from the definition \eqref{eq:stresstensor} employing the virtual element basis functions for the approximation of the velocity gradients,
\begin{equation}
	\gra \vel_{\cv}^{\ast} =  \sum\limits_{i=1}^{\NDOF} \gra \varphi_{i} \, \dofv{i}{\ast} \approx \sum\limits_{i=1}^{\NDOF} \projLm \varphi_{i} \, \dofv{i}{\ast}.  
\end{equation}

Finally, the load term $\mathbf{F}^{n}_{\vel,\cv}$ is approximated making use again of the $L_2$ projector $\projL$ as
\begin{equation}
	\left( \mathbf{F}^{n}_{\vel,\cv}\right)_{j} = \int\limits_{\cv} \wsstar \, \projL \varphi_{j} \dV.
\end{equation}

The discrete viscous sub-system \eqref{eqn.discrweakform_visc} is solved using a matrix-free GMRES method \cite{SS86}. Once the intermediate velocity field $\vel^\ast$ is obtained, the intermediate momentum $\wstar$ is computed for each element from \eqref{eqn.intermediatemomentumv} as follows:
\begin{equation}
	\wstar_{\cv} = \wsstar_{\cv} + \dt \int\limits_{\cvb} \btaus_{\cv} \cdot \nn \, \dS,
	\label{eqn.wstar}
\end{equation}
with $\wsstar_{\cv}$ already evaluated with the explicit finite volume scheme \eqref{eqn.mom_fv}. In the same manner, the work of the viscous forces needed to obtain the intermediate pressure $\pstar$ in \eqref{eqn.pressure_tdwv} is calculated for each control volume as 
\begin{equation}
	\pstar_{\cv} = \psstar_{\cv} + \dt \int\limits_{\cvb} \left( \btaus_{\cv} \, \vstar_{\cv} \right) \cdot \nn \, \dS,
	\label{eqn.ekstar}
\end{equation}
where $\psstar_{\cv}$ is known from the convective stage \eqref{eqn.ek_fv}. Notice that the above boundary integrals can be efficiently calculated using Gauss-Lobatto quadrature rules, since the degrees of freedom of the local virtual element space for $\vel^\ast_{\cv}$, and consequently for $\btaus_{\cv}$, are readily available at quadrature points (see Section \ref{ssec.lvem}).

\subsubsection{Discretization of the pressure sub-system}
The discrete pressure is approximated in the local virtual element space $V^{h}_{\kk}\left(\cv\right)$ on each control volume $\cv\in \mathcal{T}_{\Omega}$ using the virtual basis functions as
\begin{gather}
	\press^{n+1}\left(\xx \right)_{\mid \xx\in\cv}  = \sum\limits_{i=1}^{\NDOF} \varphi_{i} \dofp{i}{n+1}. \label{eqn.pvem}
\end{gather}
Choosing $\tf =\varphi_{j}$ and substituting \eqref{eqn.pvem} into the weak problem \eqref{eqn.pressure_vf4}, we get the local discrete pressure sub-system assuming $ \cvb \notin \partial \Omega$:
\begin{gather}
	\frac{1}{ \gamma-1 }\int\limits_{\cv}^{} \varphi_{i}  \varphi_{j}	\dV \dofp{i}{n+1}
	+ \dt^2 \int\limits_{\cv}^{} \ent^n  \gra \varphi_{i}  \cdot \gra \varphi_{j} \dV  \dofp{i}{n+1}
	=  	\frac{1}{ \gamma-1 } \int\limits_{\cv}^{}\press^{n} \varphi_{j} \dV 
	+ \int\limits_{\cv}^{} \pstar \varphi_{j} \dV	 
	\notag\\- \halb \int\limits_{\cv}^{} \frac{1}{\rho^{n+1}}\ww^{n+1}\cdot \ww^{n+1}	\varphi_{j} \dV
	+ \dt \int\limits_{\cv}^{}   \left(\ent^n  \wstar \right)\cdot \gra \varphi_{j} \dV. 
	 \label{eqn.pressure_disc}
\end{gather}
System \eqref{eqn.pressure_disc} can be rewritten in matrix-vector form as
\begin{equation}\label{eqn.discrweakform_press} 
		\frac{1}{\gamma-1} \mathbf{M}_{\cv} \, \dofpp{}{n+1} + \Delta t^2 \mathbf{K}_{\cv} \, \dofpp{}{n+1}= \mathbf{F}^{n}_{\press,\cv},
\end{equation}
with the mass matrix, stiffness matrix and the right hand side term given by
\begin{eqnarray}
	&\mathbf{M}_{\cv} =&  \int\limits_{\cv}^{} \varphi_{i}  \varphi_{j}	\dV, \\
	&\mathbf{K}_{\cv} =& \int\limits_{\cv}^{} \ent^n  \gra \varphi_{i}  \cdot \gra \varphi_{j} \dV,
	\\
	&\mathbf{F}^{n}_{\press,\cv}  =&  \int\limits_{\cv}^{}\left( \frac{1}{ \gamma-1 } \press^{n}
	+  \pstar - \halb \frac{\ww^{n+1}\cdot\ww^{n+1}}{\rho^{n+1}}	\right) \varphi_{j} \dV
	+ \dt \int\limits_{\cv}^{}   h^n  \wstar \cdot \gra \varphi_{j} \dV
	. \label{eqn.rhspsys}
\end{eqnarray}

The above mass matrix is approximated making use of the projection operator $\projL$, as introduced in \eqref{proj_M4}, by assuming a constant unity density, that is
\begin{align}\label{proj_M5} 
	\mathbf{M}_{\cv} = \mathbf{C} ^{\top}\left(\mathbf{H}\right)^{-1}\mathbf{C} +  \cvv{}(\mathbf{I}-\projLb)^\top(\mathbf{I}-\projLb),
\end{align}
with the matrix definitions introduced in \eqref{eqn.HCmatrix}. The stiffness matrix $\mathbf{K}_{\cv}$ contains a space dependent coefficient, which is the enthalpy $\ent^n$. Notice that this matrix is formally equal to the stiffness matrix that is retrieved in the weak problem for the solution of the shallow water equations, as fully detailed in \cite{HybridFVVEMinc}. Indeed, the role of the total water depth is here replaced by the specific enthalpy, while keeping the mathematical structure identical. Thus, following \cite{HybridFVVEMinc}, the stiffness matrix is evaluated as
\begin{align}\label{K_consiststa} 
	\mathbf{K}^{n}_{\cv} = (\accentset{\star}{\mathbf{\Pi}}^{0,x}_{\kk-1})^\top \mathbf{H}^{\ent} \accentset{\star}{\mathbf{\Pi}}^{0,x}_{\kk-1} + (\accentset{\star}{\mathbf{\Pi}}^{0,y}_{\kk-1})^\top \mathbf{H}^{\ent} \accentset{\star}{\mathbf{\Pi}}^{0,y}_{\kk-1}
	+ \bar{\ent}^{n} (\mathbf{I} - \mathbf{\Pi}^{\nabla}_{k})^\top(\mathbf{I} - \mathbf{\Pi}^{\nabla}_{k}),
\end{align}
with $\bar{\ent}^{n}$ being the mean value of the enthalpy over the element $\cv$ and the enthalpy weighted  matrix $\mathbf{H}^{\ent}$ defined as
\begin{align}\label{H-h} 
	\left(\mathbf{H}^\ent\right)_{\alpha,\beta} := \int\limits_{\cv} \ent^n m_\alpha m_\beta \, \dV, \qquad \alpha,\beta\in\left\lbrace 1,\dots, n_{\kk-1}\right\rbrace.
\end{align}
The right hand side term \eqref{eqn.rhspsys} is determined by resorting to the $L_2$ projector $\projL$ as
\begin{equation}
	\left( \mathbf{F}^{n}_{\press,\cv}\right)_{j} = \int\limits_{\cv}^{}\left( \frac{1}{ \gamma-1 } \press^{n}
	+  \pstar - \halb \frac{\ww^{n+1}\cdot\ww^{n+1}}{\rho^{n+1}} \right)  \projL \varphi_{j} \dV
	+ \dt \int\limits_{\cv}^{}   \ent^n  \wstar \cdot \gra \projL \varphi_{j} \dV.
\end{equation}
The obtained implicit pressure sub-system is symmetric by construction, since both matrices $\mathbf{M}_{\cv}$ and $\mathbf{K}^{n}_{\cv}$ are symmetric. Therefore, it is solved using an efficient matrix-free conjugate gradient algorithm, hence providing the vectors of degrees of freedom $\dofpp{\cv}{n+1}$ for the new pressure field at each control volume.

\subsubsection{High order interpolation between finite volume and virtual element spaces}
Let us note that the convection stage, presented in Section~\ref{sec.convective}, computes cell averaged values for the intermediate velocity and pressure that need to be transferred to the virtual element space before the VEM is applied for the solution of the viscous and pressure weak problems formulated in Section~\ref{ssec.weak}. 

Therefore, the intermediate solution of the explicit convective stage undergoes the CWENO reconstruction (see Section~\ref{sec.cweno}) which provides the associated high order polynomials in terms of the Taylor modal basis functions \eqref{eqn.recPoly}-\eqref{eqn.Voronoi_modal}. Next, using the virtual element mass matrix \eqref{proj_M5}, we construct the operator 
\begin{equation}
	\Vop_{\cv} = \left( \mathbf{M}_{\cv} \right)^{-1} \, \int \limits_{\cv}  \, \projL \varphi_{i} \, m_{\alpha} \, \dV, \qquad i\in\left\lbrace 1,\dots, \NDOF\right\rbrace, \quad \alpha\in\left\lbrace 1,\dots, n_{\kk}\right\rbrace,\label{eqn.fv2vem}
\end{equation} 
that accounts for the mapping from the finite volume to the virtual element space. 
Similarly, after the implicit pressure stage, the solution is projected to the finite volume space using the operator
\begin{equation}
	\Cop_{\cv} = \left( \int \limits_{\cv} m_{\alpha} m_{\beta} \, \dV \right)^{-1} \, \int \limits_{\cv}  m_{\alpha} \, \projL \varphi_{i} \, \dV, \qquad i\in\left\lbrace 1,\dots, \NDOF\right\rbrace, \quad \alpha,\beta \in\left\lbrace 1,\dots, n_{\kk}\right\rbrace. \label{eqn.vem2fv}
\end{equation}
Obviously, the above operators verify the consistency property 
\begin{equation}
	\Cop_{\cv} \, \Vop_{\cv} = \mathbf{I}_{[n_{\kk} \times n_{\kk}]}, \qquad \Vop_{\cv} \, \Cop_{\cv} = \mathbf{I}_{[\NDOF \times \NDOF]},
\end{equation} 
with $\mathbf{I}$ the identity operator.

\subsection{Global conservation of the total energy}
We recall that the discrete pressure sub-system is solved twice, namely we first solve the weak problem \eqref{eqn.pressure_vf3} for the \textit{pressure flux}, then the new momentum is updated with \eqref{eqn.velpicupdate}, that is 
\begin{equation}
	\ww_{\cv}^{n+1} = \wstar_{\cv} - \dt \int\limits_{\cvb} \tilde{\mathbf{p}}_{\cv}^{n+1} \cdot \nn \, \dS,
	\label{eqn.wnew}
\end{equation}
and finally the new \textit{pressure state} is obtained as the solution of the weak problem \eqref{eqn.pressure_vf4}. Similarly to the computation of the viscous stress contribution \eqref{eqn.wstar}, even for obtaining the new momentum we exploit Gauss-Lobatto quadrature formulae for the evaluation of the boundary integral in \eqref{eqn.wnew}. The new kinetic energy is then calculated for each cell $\cv$  as
\begin{equation}
	\left( \rho^{n+1}\ke^{n+1} \right)_{\cv}= \frac{1}{\left|\cv\right|} \int\limits_{\cv} \halb \frac{\ww^{n+1}_{\cv}\cdot\ww^{n+1}_{\cv}}{\rho^{n+1}_{\cv}} \dV,
\end{equation}
which will be needed at the next time step  within the convective stage \eqref{eqn.pressure_td}.

\begin{theorem}\label{th.energy}
	Assuming impermeable boundary conditions $\int\limits_{\partial \Omega} \ww \cdot \nn = 0$, the semi-discrete scheme \eqref{eqn.pressure_vf4} with the intermediate results for pressure and kinetic energy, given by \eqref{eqn.pressure_tdwv} and \eqref{eqn.pressure_tdsv}, respectively, is globally energy conserving in the sense that
	\begin{equation}
		\int\limits_{\Omega}^{} \frac{\rho^{n+1} E^{n+1} - \rho^{n} E^{n}}{\Delta t} \, \dV = 0.
	\end{equation}
\end{theorem}

\begin{proof}
Without loss of generality, we assume $\tf=1 \in \HONE\left(\Omega\right)$ in \eqref{eqn.pressure_vf4}, hence obtaining
\begin{gather}
	\int\limits_{\Omega}^{} \frac{\press^{n+1}}{ \gamma-1 }	\dV +  \int\limits_{\Omega}^{} \halb\frac{\ww^{n+1}\cdot \ww^{n+1}}{\rho^{n+1}} \dV
	=  	\int\limits_{\Omega}^{} \frac{\press^{n}}{ \gamma-1 }  \dV 
	+ \int\limits_{\Omega}^{} \pstar \dV, \label{eqn.pressure_th_1}
\end{gather}
where the boundary term vanishes because of the impermeability condition. Using the definition of the intermediate pressure field $\pstar$ and the convective contribution of the kinetic energy $\psstar$, given by \eqref{eqn.pressure_tdwv} and \eqref{eqn.pressure_tdsv}, respectively, from the above expression we get
\begin{gather}
	\int\limits_{\Omega}^{} \frac{\press^{n+1}}{ \gamma-1 } 	\dV + \int\limits_{\Omega}^{} \halb \frac{\ww^{n+1}\cdot \ww^{n+1}}{\rho^{n+1}}\dV
	=  \int\limits_{\Omega}^{} \frac{\press^{n}}{ \gamma-1 }  \dV 
	+ \int\limits_{\Omega}^{} \rho^{n}\ke^{n} - \dt \dive \left(\ke^{n} \ww^{n} + \hf^{n} \right) + \dt \dive \left( \btaus\vstar\right) \dV. \label{eqn.pressure_th_2}
\end{gather}
In the above right hand side term, the fluxes $\left( \ke^{n} \ww^{n} + \hf^{n} \right)$ are computed by means of the finite volume scheme \eqref{eqn.ek_fv} with the conservative numerical fluxes \eqref{eqn.rusanov} for each cell $\cv$. Likewise, the divergence of the work of the stress tensor is discretized in conservative flux form according to \eqref{eqn.ekstar} for all the control volumes. Therefore, the integral over the entire computational domain of these terms vanishes due to the telescopic property of finite volume schemes combined with impermeable boundary conditions. What is left of the above equation then writes
\begin{gather}
	\int\limits_{\Omega}^{} \frac{\press^{n+1}}{ \gamma-1 } 	\dV + \int\limits_{\Omega}^{} \halb \frac{\ww^{n+1}\cdot \ww^{n+1}}{\rho^{n+1}}\dV
	=  \int\limits_{\Omega}^{} \frac{\press^{n}}{ \gamma-1 }  \dV 
	+ \int\limits_{\Omega}^{} \rho^{n}\ke^{n} \dV, \label{eqn.pressure_th_3}
\end{gather}
that, recalling the definition of the total energy \eqref{eq:E.as.pk}, yields
\begin{equation}
	\int\limits_{\Omega}^{} \rho^{n+1} E^{n+1} \, dV = \int\limits_{\Omega}^{} \rho^{n} E^{n} \, dV.
\end{equation}
Consequently, the total energy is globally conserved.
\end{proof}

\section{Numerical results}\label{sec:numericalresults}
In what follows, we assess the novel semi-implicit hybrid FV/VEM (SI-FVVEM) methodology for all Mach number flows using a set of classical test cases ranging form high Mach number flows to the incompressible limit. As previously described, the flux splitting of our novel scheme leads to an explicit sub-system for transport equations and two implicit sub-systems for the viscous and pressure terms. As a consequence, the time step of the overall method is restricted only by the CFL condition related to the bulk velocity $|\vel|$ and not to the speed of pressure waves nor the viscosity coefficient. Unless stated the contrary, all test cases have been therefore run computing the time step from
\begin{equation}
	\dt = \min\limits_{P_i} \left\lbrace \dt_{i} \right\rbrace, \qquad \dt_{i} = \mathrm{CFL} \, \frac{h_{i} }{|\vel_{i}| + c_{\alpha}},
\end{equation}
with $\mathrm{CFL}= 0.5$ and the characteristic size $h_i$ of the Voronoi element $P_i$ given by \eqref{eqn.h}. The coefficient $c_{\alpha}$ accounts for an artificial viscosity in case of strong shock waves, and it is not considered if not otherwise stated, i.e. $c_{\alpha}=0$. The ratio of specific heats is set to $\gamma=1.4$ and the gas constant is $R=1$. The  SI units are considered in all test cases.

\subsection{Numerical convergence study}
To assess the accuracy of the novel SI-FVVEM schemes, we study the isentropic vortex problem proposed in \cite{HuShuTri}. The known analytical solution, defined in the computational domain $\Omega=\left[0,10\right]^2$, reads
\begin{gather}
	\vel\left(\xx,t\right) = \frac{\varepsilon}{2\pi}e^{\frac{1-r^2}{2}} \left(\begin{array}{c} 5-y \\ x-5 \end{array}\right), 
	\quad 
	\rho\left(\xx,t\right) = (1+\delta \temp)^{\frac{1}{\gamma-1}}, 
	\quad \press\left(\xx,t\right) = (\press_{0}+\delta \temp)^{\frac{\gamma}{\gamma-1}},
	\notag\\ \delta \temp = -\frac{(\gamma-1)\varepsilon^2}{8\gamma\pi^2}e^{1-r^2},
	\quad r^2 = (x-5)^2+(y-5)^2, \label{ShuVort}
\end{gather} 
with $\varepsilon=5$ the vortex strength and $p_0$ the mean pressure defined according to the desired Mach number. The simulations are run up to time $t_f=0.1$ on a set of refined Voronoi grids with periodic boundary conditions everywhere. In Tables~\ref{tab.shuconvergence}-\ref{tab.shuconvergenceo3}, we observe that the expected order of accuracy is reached for both the second and third order schemes. 
The asymptotic preserving property is analyzed by considering a set of Mach numbers ranging in $\M\in\left\lbrace 10^{-6}, 10^{-4}, 10^{-2}, 1\right\rbrace $.

\begin{table}
	\begin{center}
		\caption{Isentropic vortex. $L^2$ errors and convergence rates obtained for the isentropic vortex with $\M\in\left\lbrace 10^{-6}, 10^{-4}, 10^{-2}, 1\right\rbrace $. The results are obtained using the second order scheme in space and time. $\cvs(\Omega)=1/N_x$ denotes the characteristic size of the elements generated by considering $N_x$ points in the $x-$direction.}
			\begin{tabular}{ccccccccc}
				$\cvs(\Omega)$ & ${L_2}(\rho)$ & $\mathcal{O}(\rho)$ & ${L_2}(u)$ & $\mathcal{O}(u)$ & ${L_2}(v)$ & $\mathcal{O}(v)$ & ${L_2}(\press)$ & $\mathcal{O}(\press)$\\
				\hline
				\multicolumn{9}{c}{$M=1$}\\
				\hline
				$\frac{1}{8}$ & $1.4227\cdot 10^{-1} $ & -& $2.9804\cdot 10^{-1} $ & - & $2.9307\cdot 10^{-1} $ & -& $1.5776\cdot 10^{-1}$ & -\\
				$\frac{1}{15}$& $3.4143\cdot 10^{-2} $ & 2.27 & $6.6575\cdot 10^{-2} $ &  2.38 & $6.7252\cdot 10^{-2} $ & 2.34 & $4.3657\cdot 10^{-2}$ &2.04\\
				$\frac{1}{30}$& $6.6055\cdot 10^{-3}   $ & 2.37 & $1.3143\cdot 10^{-2}   $ &  2.34 & $1.3256\cdot 10^{-2}   $ & 2.34 & $1.0845\cdot 10^{-2}  $ &2.01\\
				$\frac{1}{45}$& $2.6854\cdot 10^{-3} $ & 2.22 & $5.3015\cdot 10^{-3} $ &  2.24 & $5.3787\cdot 10^{-3} $ & 2.22 & $4.7830\cdot 10^{-3}$ &2.02\\				
				\hline
				\multicolumn{9}{c}{$M=10^{-2}$}\\
				\hline
				$\frac{1}{8}$ & $ 1.4229\cdot 10^{-1}$ & -    & $3.0822\cdot 10^{-1} $ & - & $3.0459\cdot 10^{-1} $ & -& $3.7373\cdot 10^{-3} $ & -\\
				$\frac{1}{15}$& $ 3.4289\cdot 10^{-2}$ & 2.26 & $7.3083\cdot 10^{-2} $ & 2.29  & $7.1447\cdot 10^{-2} $ & 2.31 & $1.1765\cdot 10^{-3} $ &1.84\\
				$\frac{1}{30}$& $ 6.7729\cdot 10^{-3}  $ & 2.34 & $1.5757\cdot 10^{-2}   $ & 2.21  & $1.5381\cdot 10^{-2}   $ & 2.22 & $2.8033\cdot 10^{-4}   $ &2.07\\
				$\frac{1}{45}$& $ 2.7514\cdot 10^{-3}$ & 2.22 & $6.8921\cdot 10^{-3} $ & 2.04  & $6.6339\cdot 10^{-3} $ & 2.07 & $1.1257\cdot 10^{-4} $ &2.25\\				
				\hline
				\multicolumn{9}{c}{$M=10^{-4}$}\\
				\hline
				$\frac{1}{8}$ & $1.4229\cdot 10^{-1} $ & -& $3.9657\cdot 10^{-1} $ & - & $3.9397\cdot 10^{-1}$ & -& $1.4430\cdot 10^{-4} $ & -\\
				$\frac{1}{15}$& $3.4402\cdot 10^{-2} $ & 2.26 & $9.1457\cdot 10^{-2} $ & 2.33  & $9.0785\cdot 10^{-2}$ & 2.33 & $2.4874\cdot 10^{-5} $ &2.80 \\
				$\frac{1}{30}$& $6.8437\cdot 10^{-3}   $ & 2.33 & $  1.8174\cdot 10^{-2} $ & 2.33  & $1.7310\cdot 10^{-2}  $ & 2.39 & $5.6645\cdot 10^{-6} $   &2.13 \\
				$\frac{1}{45}$& $2.7673\cdot 10^{-3} $ & 2.23 & $7.4450\cdot 10^{-3} $ & 2.20  & $7.3476\cdot 10^{-3}$ & 2.11 & $1.4336\cdot 10^{-6} $ &3.39 \\				
				\hline
				\multicolumn{9}{c}{$M=10^{-6}$}\\
				\hline
				$\frac{1}{8}$ & $1.4231\cdot 10^{-1}$ & -    & $4.0283\cdot 10^{-1}$ & -    & $4.0068\cdot 10^{-1}$ & -    & $1.3129\cdot 10^{-6}$ & -\\
				$\frac{1}{15}$& $3.4412\cdot 10^{-2}$ & 2.26 & $9.1890\cdot 10^{-2}$ & 2.35 & $9.1225\cdot 10^{-2}$ & 2.35 & $2.2917\cdot 10^{-7}$ & 2.78 \\
				$\frac{1}{30}$& $6.8530\cdot 10^{-3}  $ & 2.33 & $ 1.8260\cdot 10^{-2} $ & 2.33 & $  1.7392\cdot 10^{-2}$ & 2.39 & $  5.2820\cdot 10^{-8}$ & 2.12 \\
				$\frac{1}{45}$& $2.7706\cdot 10^{-3}$ & 2.23 & $7.4470\cdot 10^{-3}$ & 2.21 & $7.3621\cdot 10^{-3}$ & 2.12 & $1.8408\cdot 10^{-8}$ & 2.60 \\
				\hline
			\end{tabular}
		\label{tab.shuconvergence}
	\end{center}
\end{table}

\begin{table}
	\begin{center}
	\caption{Isentropic vortex. $L^2$ errors and convergence rates obtained for the isentropic vortex  with $\M\in\left\lbrace 10^{-6}, 10^{-4}, 10^{-2}, 1\right\rbrace $. The results computed using the third order scheme in space and time. $\cvs(\Omega)=1/N_x$ denotes the characteristic size of the elements generated by considering $N_x$ points in the $x-$direction.}
			\begin{tabular}{ccccccccc}
				$\cvs(\Omega)$ & ${L_2}(\rho)$ & $\mathcal{O}(\rho)$ & ${L_2}(u)$ & $\mathcal{O}(u)$ & ${L_2}(v)$ & $\mathcal{O}(v)$ & ${L_2}(\press)$ & $\mathcal{O}(\press)$\\
				\hline
				\multicolumn{9}{c}{$M=1$}\\
				\hline
				$\frac{1}{8}$ & $1.5274\cdot 10^{-1} $ & -& $3.2950\cdot 10^{-1} $& - & $3.1903\cdot 10^{-1} $ & -& $5.5250\cdot 10^{-2} $ & -\\
				$\frac{1}{15}$& $3.6517\cdot 10^{-2} $ & 2.28 & $6.7172\cdot 10^{-2} $& 2.53 & $6.6709\cdot 10^{-2} $ & 2.49 & $1.1205\cdot 10^{-2} $ & 2.54 \\
				$\frac{1}{30}$& $5.6375\cdot 10^{-3} $ & 2.70 & $1.0535\cdot 10^{-2} $& 2.67 & $1.0841\cdot 10^{-2} $ & 2.62 & $2.3378\cdot 10^{-3} $ & 2.26 \\
				$\frac{1}{45}$& $1.7193\cdot 10^{-3} $ & 2.93 & $3.3469\cdot 10^{-3} $& 2.83 & $3.4571\cdot 10^{-3} $ & 2.82 & $9.3600\cdot 10^{-4} $ & 2.26 \\				
				\hline
				\multicolumn{9}{c}{$M=10^{-2}$}\\
				\hline
				$\frac{1}{8}$ & $1.5277\cdot 10^{-1} $ & -& $3.6229\cdot 10^{-1} $ & - & $3.5371\cdot 10^{-1} $ & -& $ 1.0030\cdot 10^{-2} $ & -\\
				$\frac{1}{15}$& $3.6506\cdot 10^{-2} $ & 2.28 & $7.8451\cdot 10^{-2} $ & 2.43 & $7.7056\cdot 10^{-2} $ & 2.42 & $ 1.8207\cdot 10^{-3} $ & 2.71 \\
				$\frac{1}{30}$& $5.6480\cdot 10^{-3} $ & 2.69 & $1.3186\cdot 10^{-2} $ & 2.57 & $1.3097\cdot 10^{-2} $ & 2.56 & $ 2.6715\cdot 10^{-4} $ & 2.77 \\
				$\frac{1}{45}$& $1.7172\cdot 10^{-3} $ & 2.94 & $4.1989\cdot 10^{-3} $ & 2.82 & $4.3512\cdot 10^{-3} $ & 2.72 & $ 8.5675\cdot 10^{-5} $ & 2.80 \\				
				\hline
				\multicolumn{9}{c}{$M=10^{-4}$}\\
				\hline
				$\frac{1}{8}$ & $1.5260\cdot 10^{-1} $ & -& $4.4288\cdot 10^{-1} $ & - & $4.4010\cdot 10^{-1} $ & -& $1.7587\cdot 10^{-4} $ & -\\
				$\frac{1}{15}$& $3.6444\cdot 10^{-2} $ & 2.28 & $9.1760\cdot 10^{-2} $ & 2.50 & $9.0320\cdot 10^{-2} $ & 2.52 & $2.1214\cdot 10^{-5} $ &3.36 \\
				$\frac{1}{30}$& $5.6444\cdot 10^{-3} $ & 2.69 & $1.4374\cdot 10^{-2} $ & 2.67 & $1.4227\cdot 10^{-2} $ & 2.67 & $4.1840\cdot 10^{-6} $ &2.34 \\
				$\frac{1}{45}$& $1.7147\cdot 10^{-3} $ & 2.94 & $4.4069\cdot 10^{-3} $ & 2.92 & $4.5710\cdot 10^{-3} $ & 2.80 & $8.9859\cdot 10^{-7} $ &3.79 \\	
				\hline
				\multicolumn{9}{c}{$M=10^{-6}$}\\
				\hline
				$\frac{1}{8}$ & $1.5260\cdot 10^{-1} $ & -& $4.4589\cdot 10^{-1} $ & - & $4.4308\cdot 10^{-1} $ & -& $1.6441\cdot 10^{-6} $ & -\\
				$\frac{1}{15}$& $3.6441\cdot 10^{-2} $ &2.28  & $9.1938\cdot 10^{-2} $ & 2.51  & $9.0491\cdot 10^{-2} $ & 2.53 & $2.1137\cdot 10^{-7} $ &3.26\\
				$\frac{1}{30}$& $5.6443\cdot 10^{-3} $ &2.69  & $1.4400\cdot 10^{-2} $ & 2.67  & $1.4251\cdot 10^{-2} $ & 2.67 & $4.2159\cdot 10^{-8} $ &2.33\\
				$\frac{1}{45}$& $1.7147\cdot 10^{-3} $ &2.94  & $4.4118\cdot 10^{-3} $ & 2.92  & $4.5736\cdot 10^{-3} $ & 2.80 & $1.3446\cdot 10^{-8} $ &2.82\\
				\hline
			\end{tabular}
		\label{tab.shuconvergenceo3}
	\end{center}
\end{table}

\subsection{Riemann problems}
The shock-capturing and conservation properties of our novel schemes are analyzed using a set of Riemann problems. In particular, we consider the computational domain $\Omega=[-0.5,0.5]\times[-0.05,0.05]$  and the initial conditions are given by
\begin{gather}
	\rho\left(\xx,0\right) = \left\lbrace \begin{array}{ll}
		\rho_{L} &  x \le  x_{0},\\
		\rho_{R} &   x>  x_{0},
	\end{array}\right.\quad
	u\left(\xx,0\right) =\left\lbrace \begin{array}{ll}
		u_{L}  &  x \le  x_{0},\\
		u_{R} &  x >  x_{0},
	\end{array}\right. \quad
	v\left(\xx,0\right)  =0, \quad
	\press\left(\xx,0\right)  = \left\lbrace \begin{array}{ll}
		\press_{L} &  x \le  x_{0},\\
		\press_{R} &  x >  x_{0},
	\end{array}\right. 
\end{gather}
with the left and right states, the final simulation time and the initial location of the discontinuity specified in Table~\ref{tab.RP}. We also report the number of points $N_x$ in the $x-$direction, employed to generate the Voronoi grid.
The simulations are run setting periodic boundary conditions in $y-$direction while in the $x-$direction we impose Dirichlet boundary conditions. We use a SI-FVVEM scheme with third and first order of accuracy in space and time, respectively.

\begin{table}
	\caption{Riemann problems. Left and right initial states, initial position of the discontinuity $x_{0}$, final time $t_{f}$, and characteristic mesh size $h=1/N_x$, for each Riemann problem.}
	\begin{center}
		\begin{tabular}{cccccccccc}
			Test &  $\rho_{L}$ &  $\rho_{R}$  &  $u_{L}$ &  $u_{R}$ &  $\press_{L}$ &  $\press_{R}$ & $x_{0}$ & $t_{f}$ & $h$ \\ \hline
			RP1 & $ 1 $ & $ 0.125 $ & $ 0 $ & $ 0 $  & $ 1 $ & $ 0.1 $& $ 0 $ &  $ 0.2 $  & $1/200$\\ 
			RP2 & $ 1 $ & $ 1 $  & $ -1 $ & $ 1 $ & $ 0.4 $ & $ 0.4 $& $ 0 $ &  $ 0.15 $  & $1/200 $\\  
			RP3 & $ 0.445 $ & $ 0.5 $ & $ 0.698 $ & $ 0 $ & $ 3.528 $ & $ 0.571 $ & $ 0 $ &  $ 0.14 $  & $1/200 $\\
			RP4 & $ 1 $ & $ 1 $ & $ -19.59745 $ & $ -19.59745 $ & $ 1000.0 $ & $ 0.01 $ & $ 0.3 $ &  $ 0.012 $ & $1/200$\\
			\hline
		\end{tabular} \label{tab.RP}
	\end{center}
\end{table}

\begin{figure}[!htbp]
	\begin{center}
		\begin{tabular}{ccc}
			\includegraphics[trim= 5 5 5 5, clip,width=0.33\textwidth]{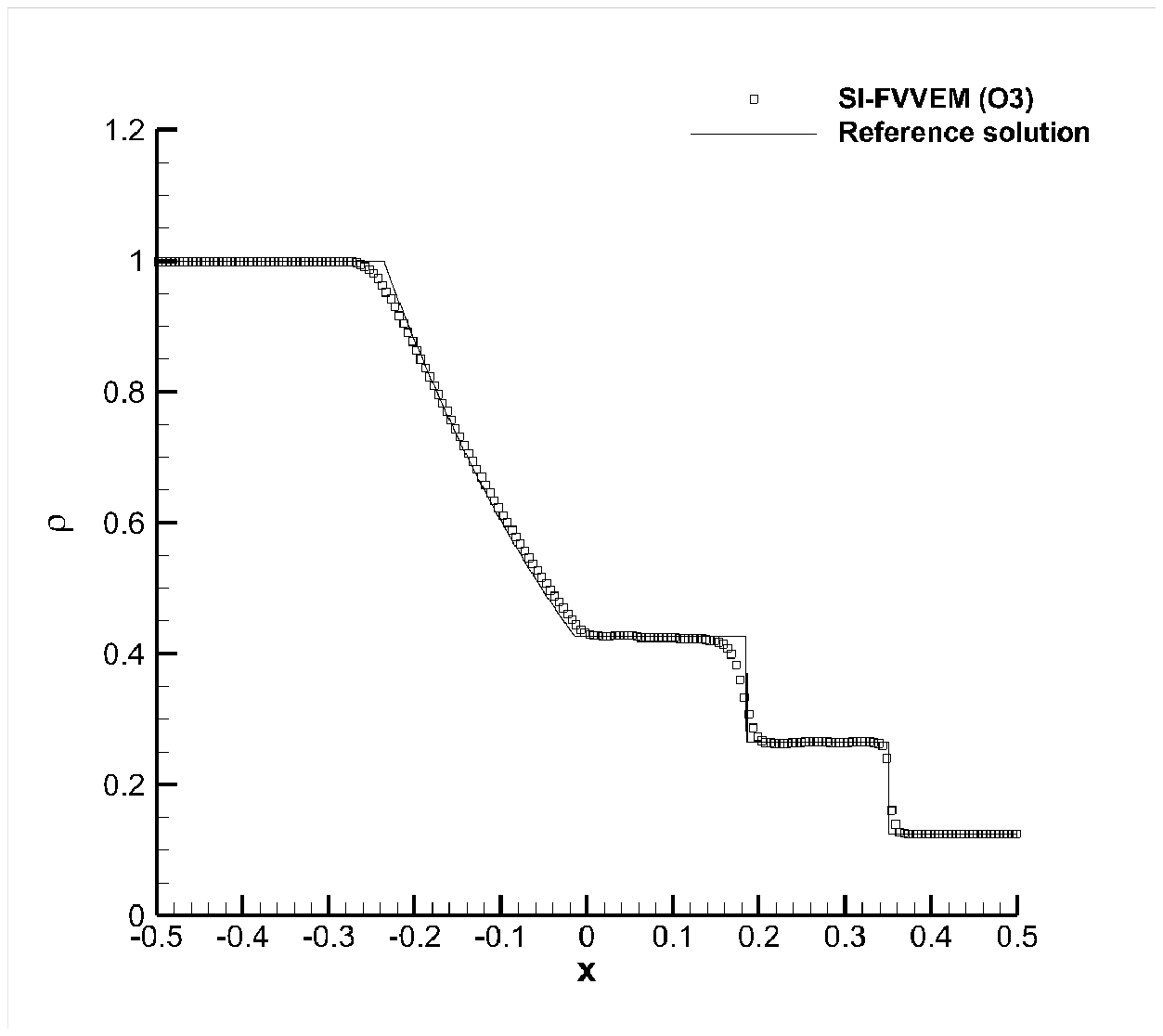}  &
			\includegraphics[trim= 5 5 5 5, clip,width=0.33\textwidth]{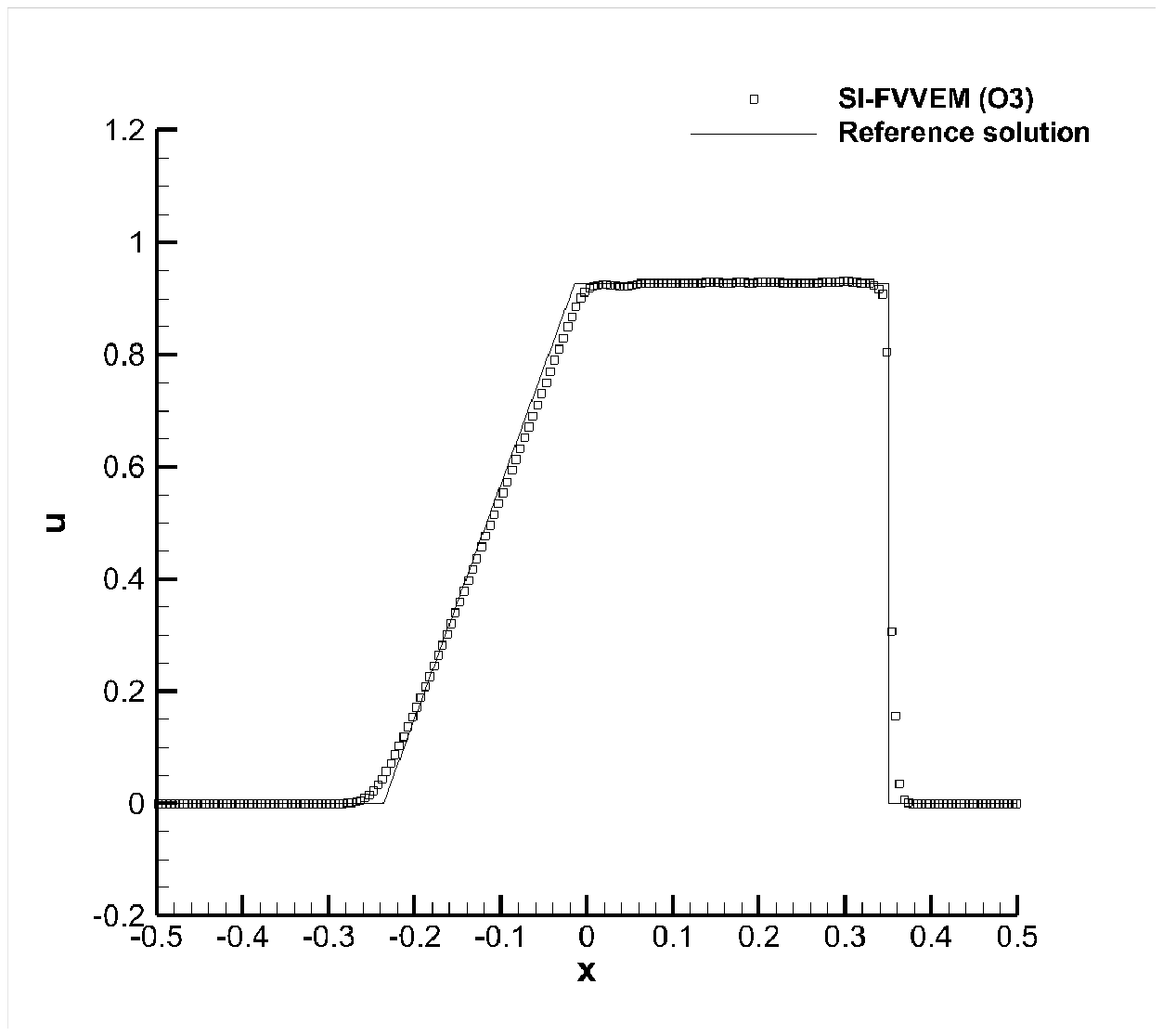} &
			\includegraphics[trim= 5 5 5 5, clip,width=0.33\textwidth]{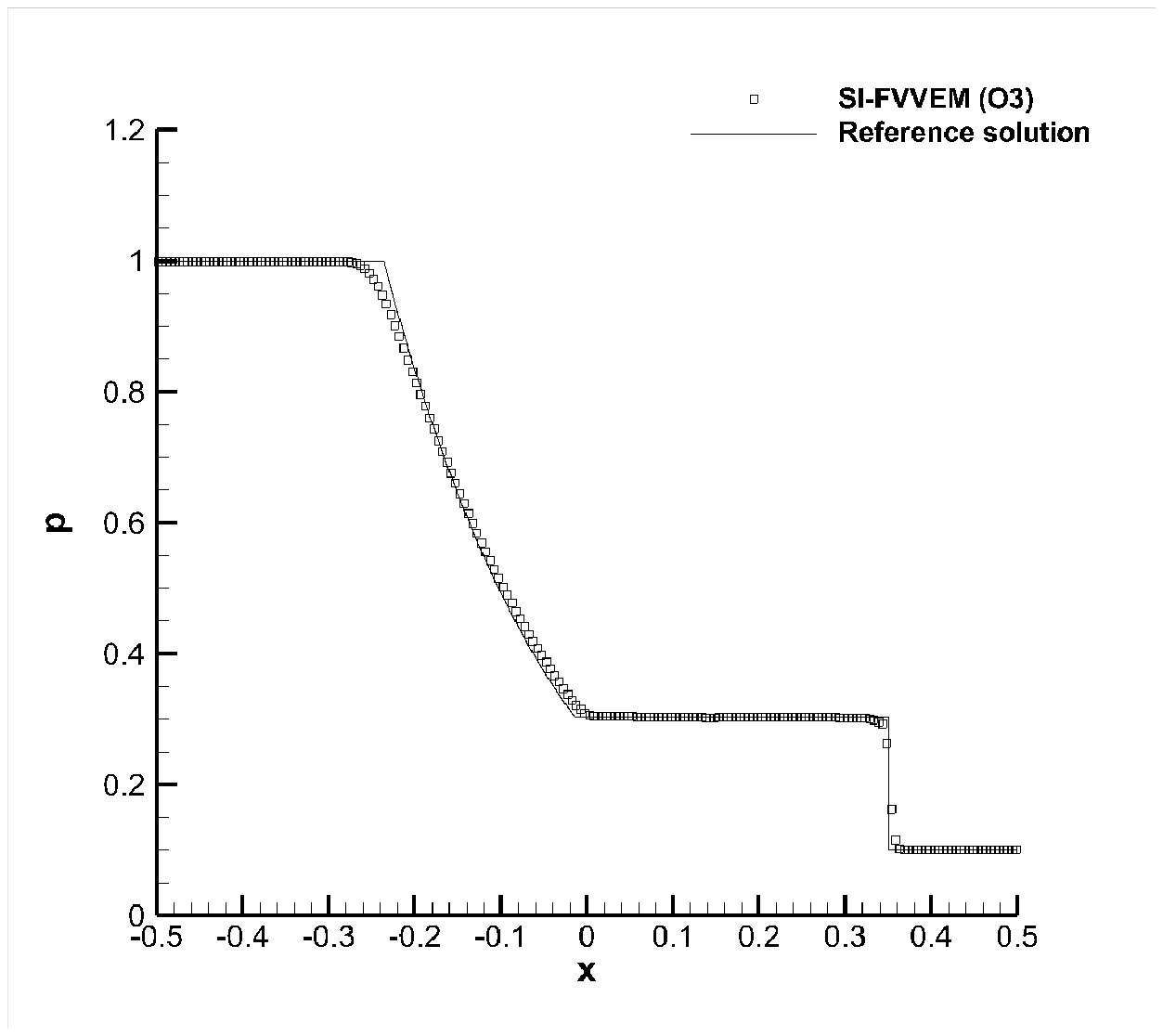}  \\
			\includegraphics[width=0.33\textwidth]{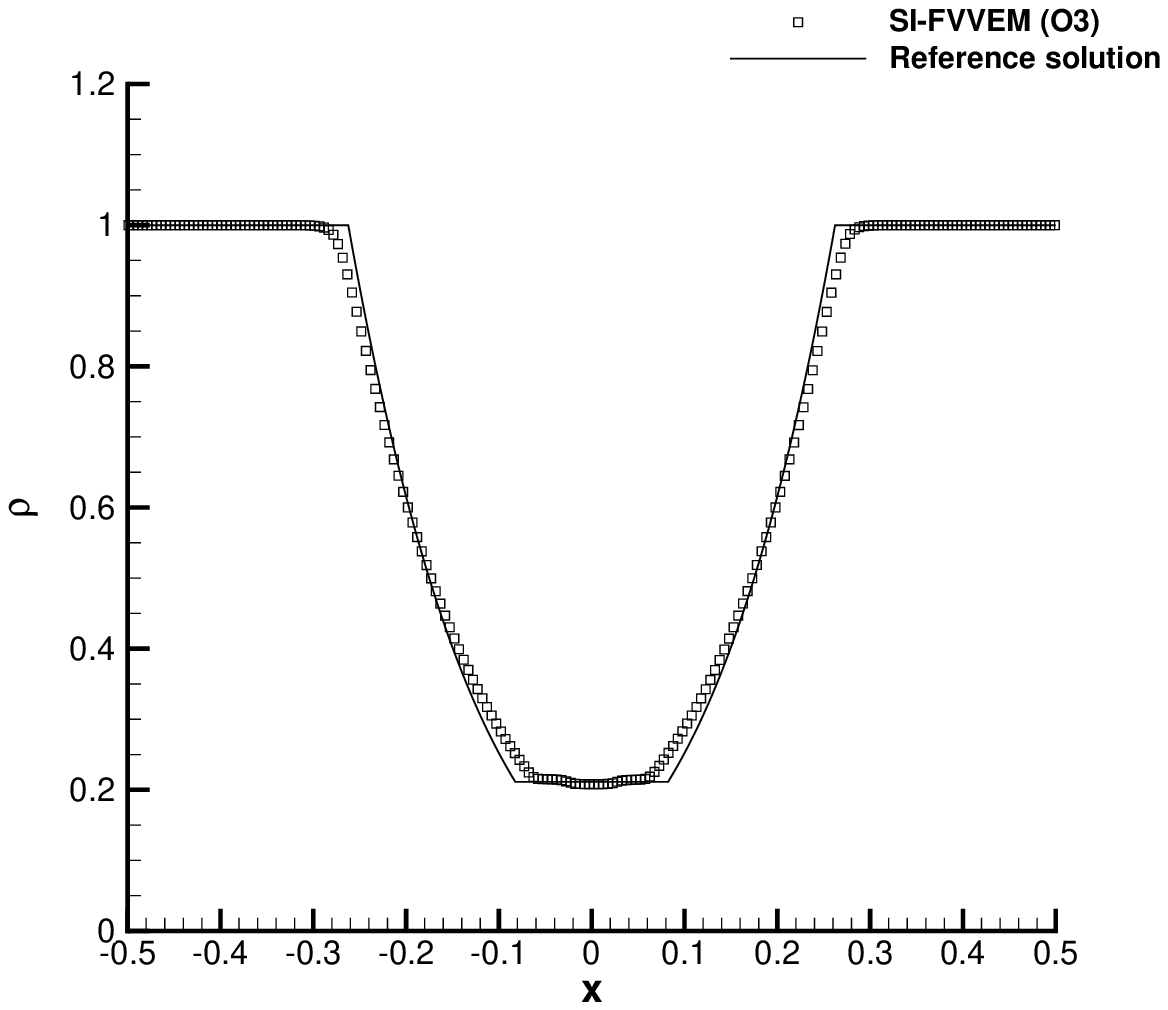}  &
			\includegraphics[width=0.33\textwidth]{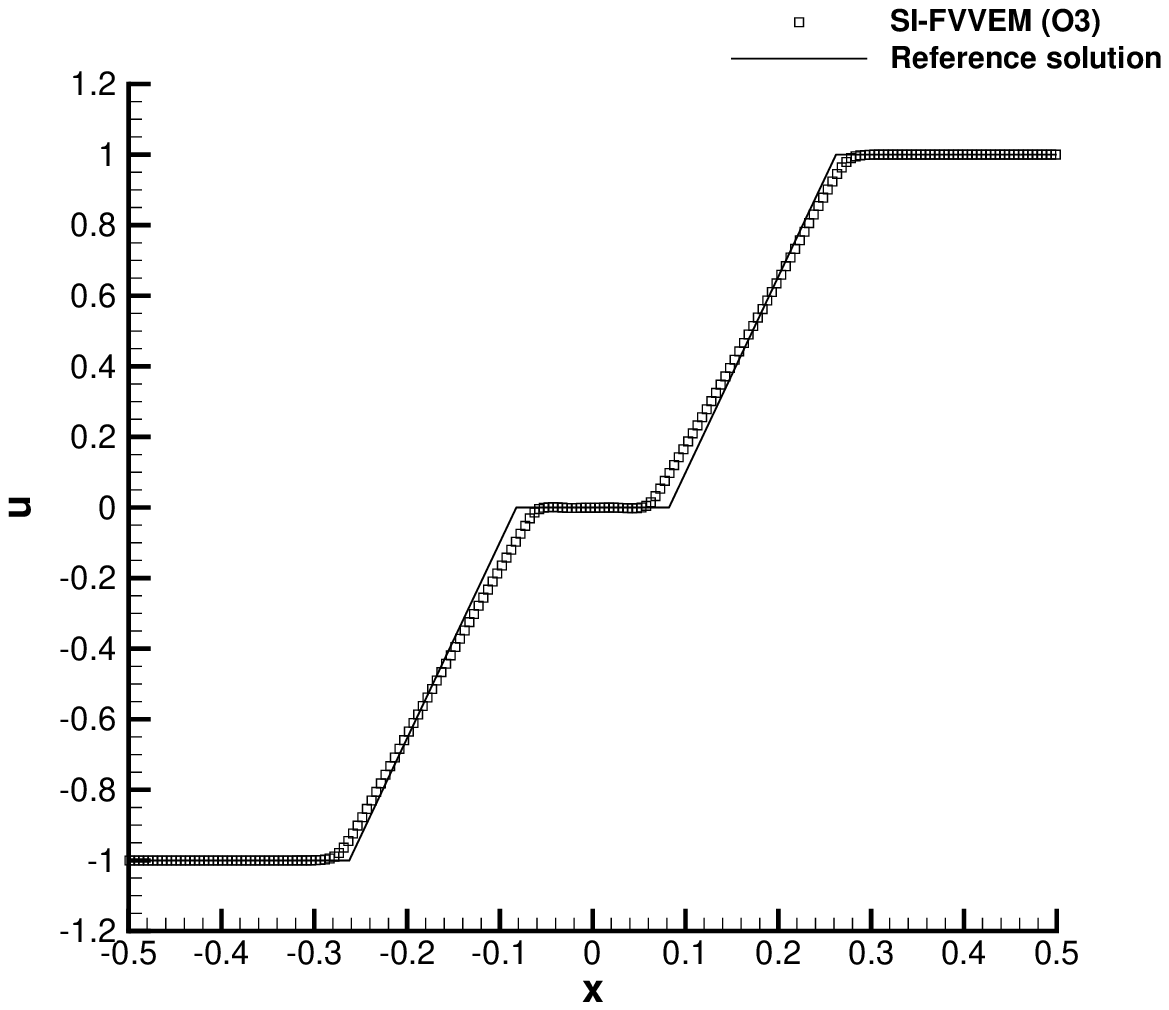} &
			\includegraphics[width=0.33\textwidth]{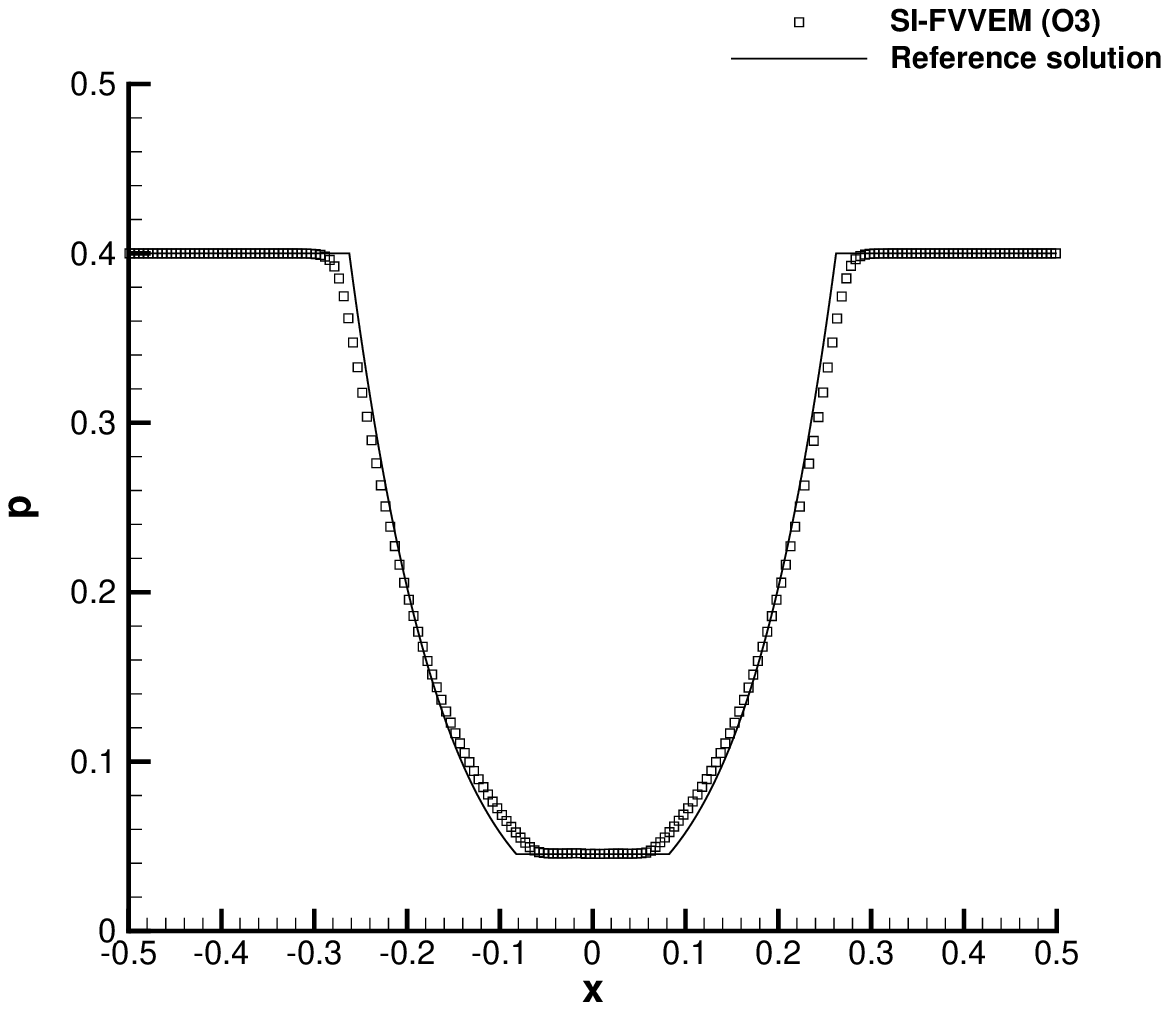}  
		\end{tabular}
		\caption{Riemann problems. RP1 at time $t_f=0.2$ (top) and RP2 at time $t_f=0.15$. Left: density $\rho$. Center: horizontal velocity component $u$. Right: pressure $p$.}
		\label{fig.RP12}
	\end{center}
\end{figure}

The first Riemann problem RP1, corresponds to the classical Sod shock tube test put forward in \cite{Sod78}. The numerical solution reported in Figure~\ref{fig.RP12} shows that the shock, contact and rarefaction waves are properly captured.
A good agreement between the numerical and exact solutions is also observed for the double rarefaction test RP2 in Figure~\ref{fig.RP12}.

\begin{figure}[!htbp]
	\begin{center}
		\begin{tabular}{ccc}&&\\
			\includegraphics[width=0.33\textwidth]{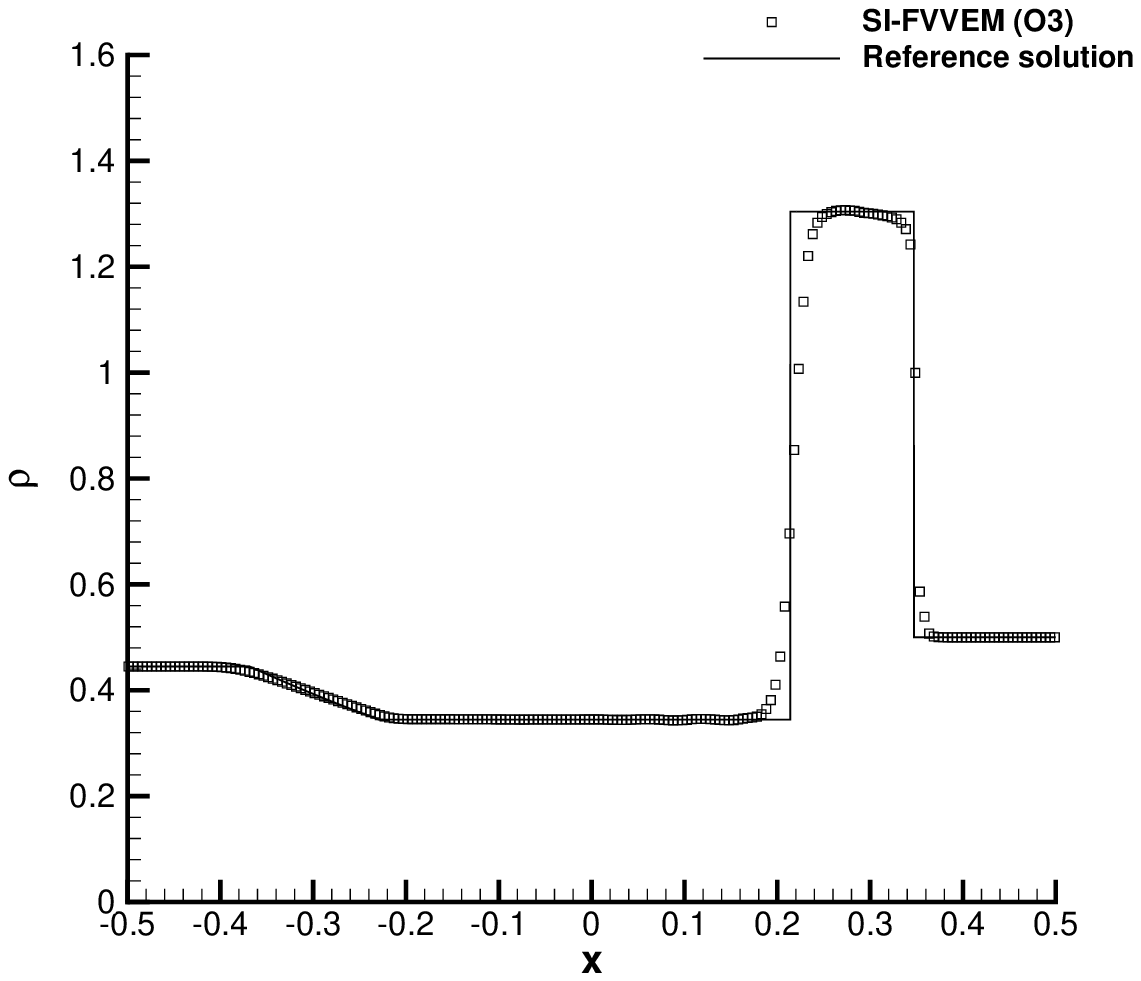}  &
			\includegraphics[width=0.33\textwidth]{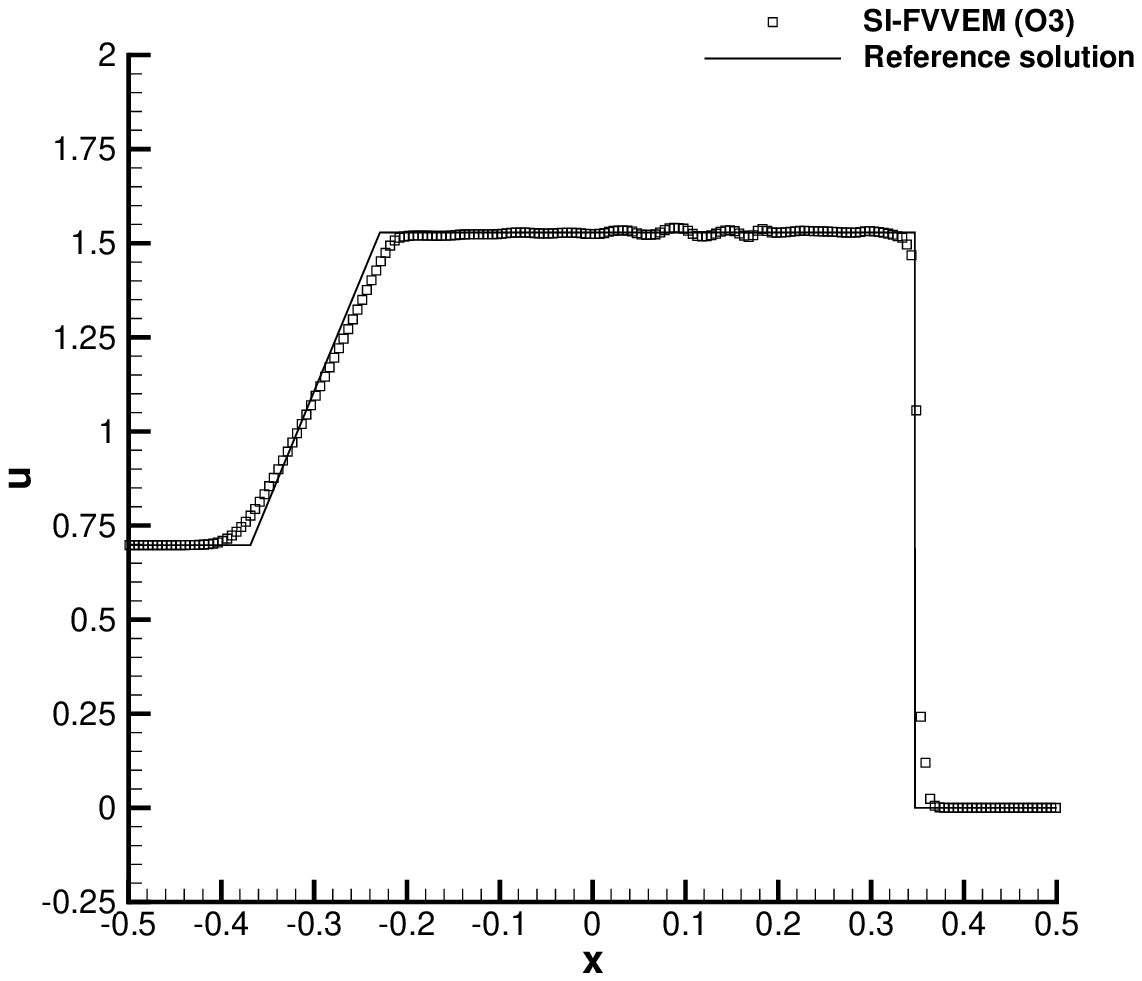} &
			\includegraphics[width=0.33\textwidth]{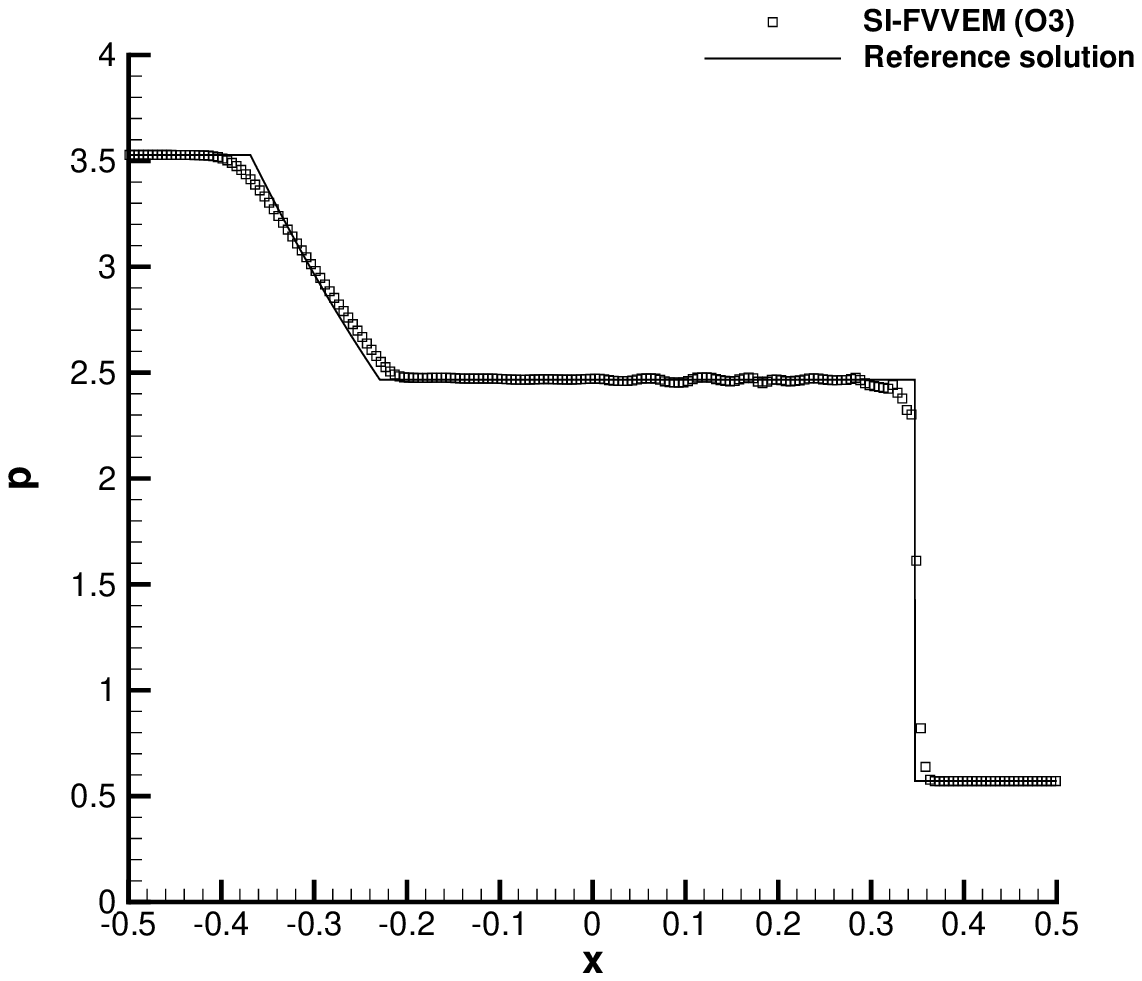}  \\
			\includegraphics[width=0.33\textwidth]{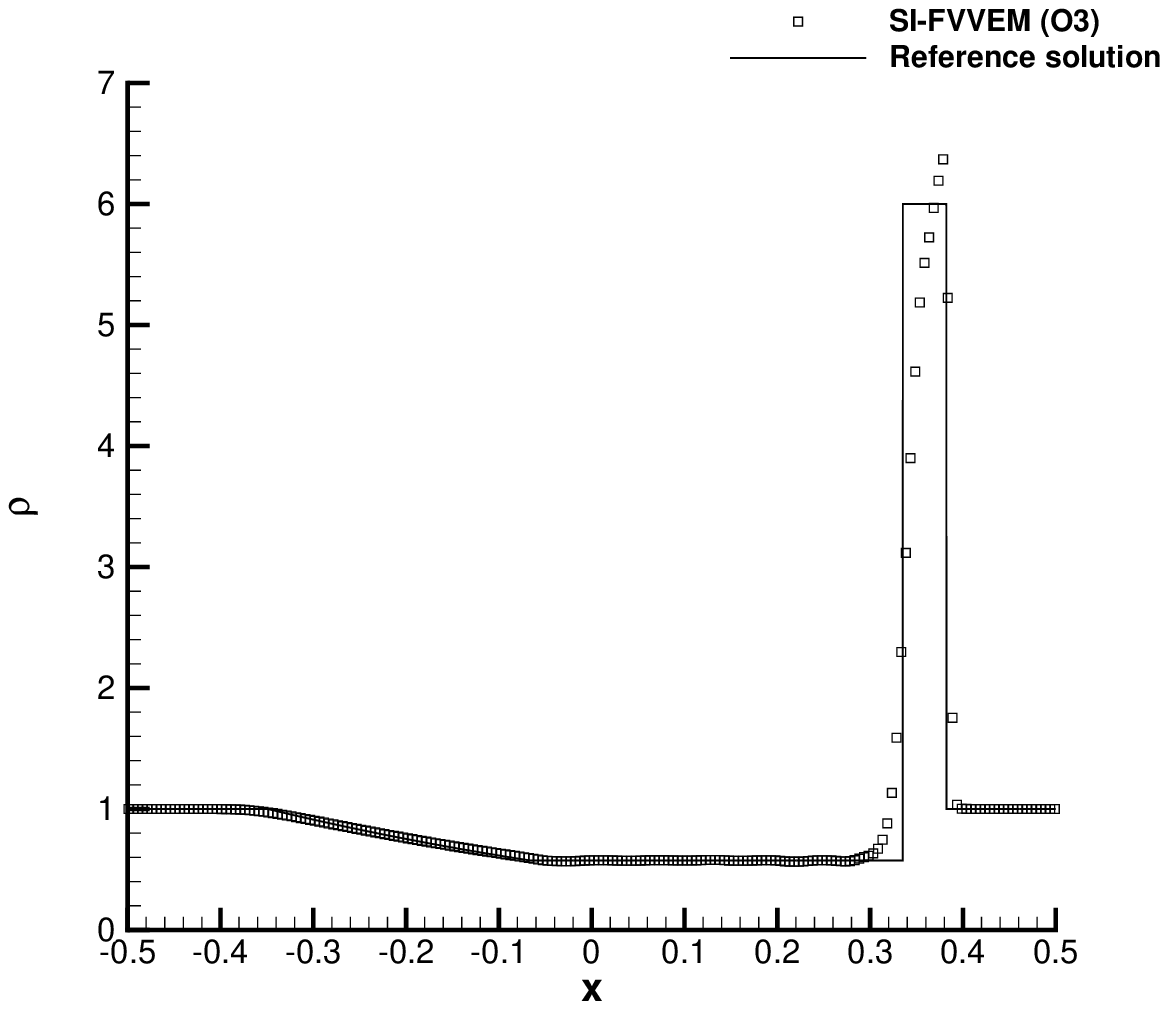}  &
			\includegraphics[width=0.33\textwidth]{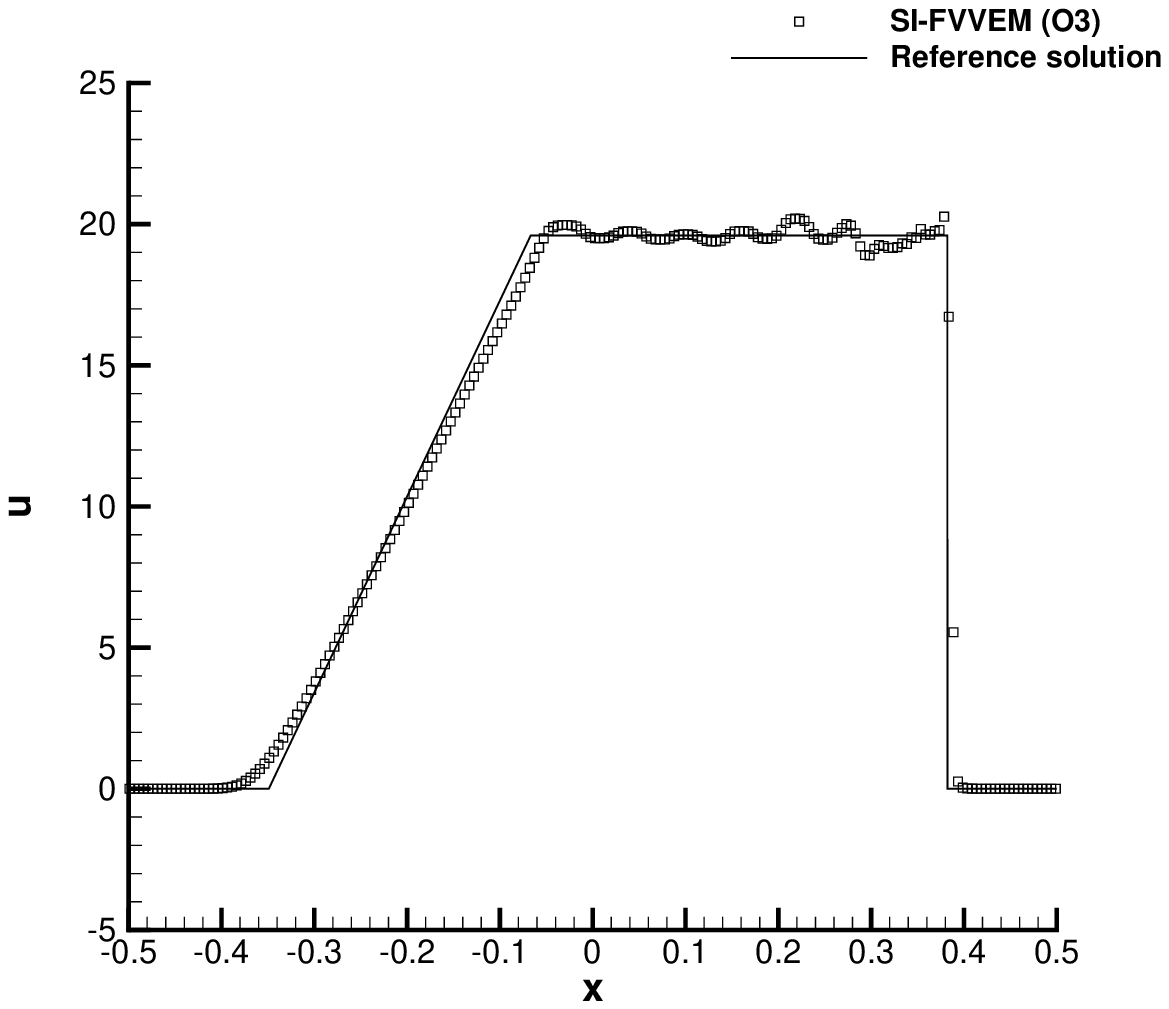} &
			\includegraphics[width=0.33\textwidth]{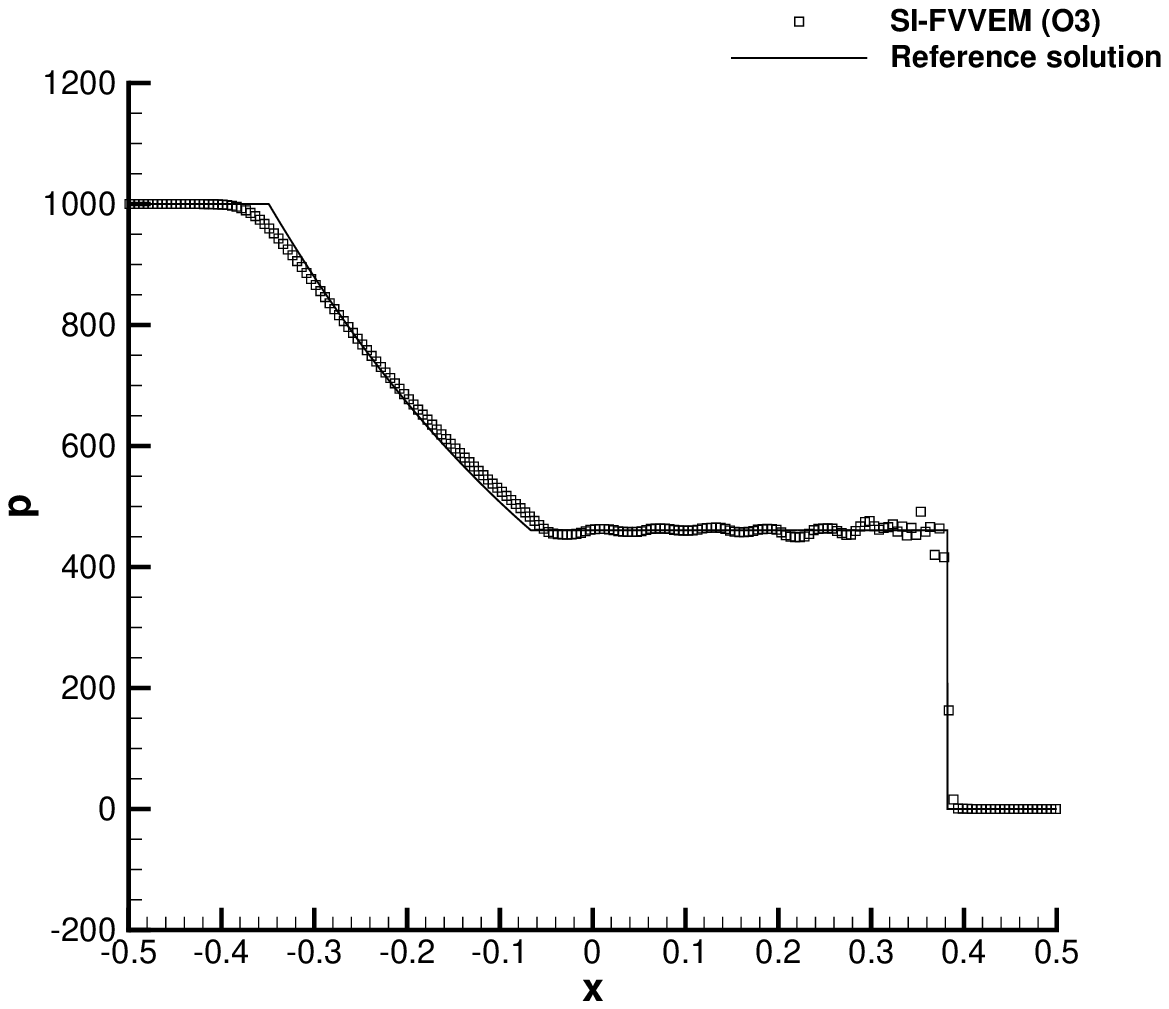}  \\
		\end{tabular}
		\caption{Riemann problems. RP3 at time $t_f=0.14$ (top) and RP4 at time $t_f=0.012$. Left: density $\rho$. Center: horizontal velocity component $u$. Right: pressure $p$.}
		\label{fig.RP34}
	\end{center}
\end{figure}

The third and fourth test cases correspond to the Lax shock tube problem and a severe test case proposed in \cite{Toro} as a modification of the left half of the blast problem presented in \cite{WC84}. The comparison of the numerical results obtained with the SI-FVVEM scheme against the corresponding exact solutions is provided in Figure \ref{fig.RP34}. A good match is observed for RP3. Regarding RP4, also the location of the right traveling shock wave as well as the left rarefaction wave correspond well with the known analytical solution. Moreover, the stationary discontinuity originated by the large pressure jump in the initial conditions is also properly captured. Even thought, small oscillations appear in the plateaux of the velocity and pressure fields that may be related to the high order space discretization used with no additional numerical dissipation.

\begin{figure}[!htbp]
	\begin{center}
		\begin{tabular}{cc}
			\includegraphics[width=0.47\textwidth]{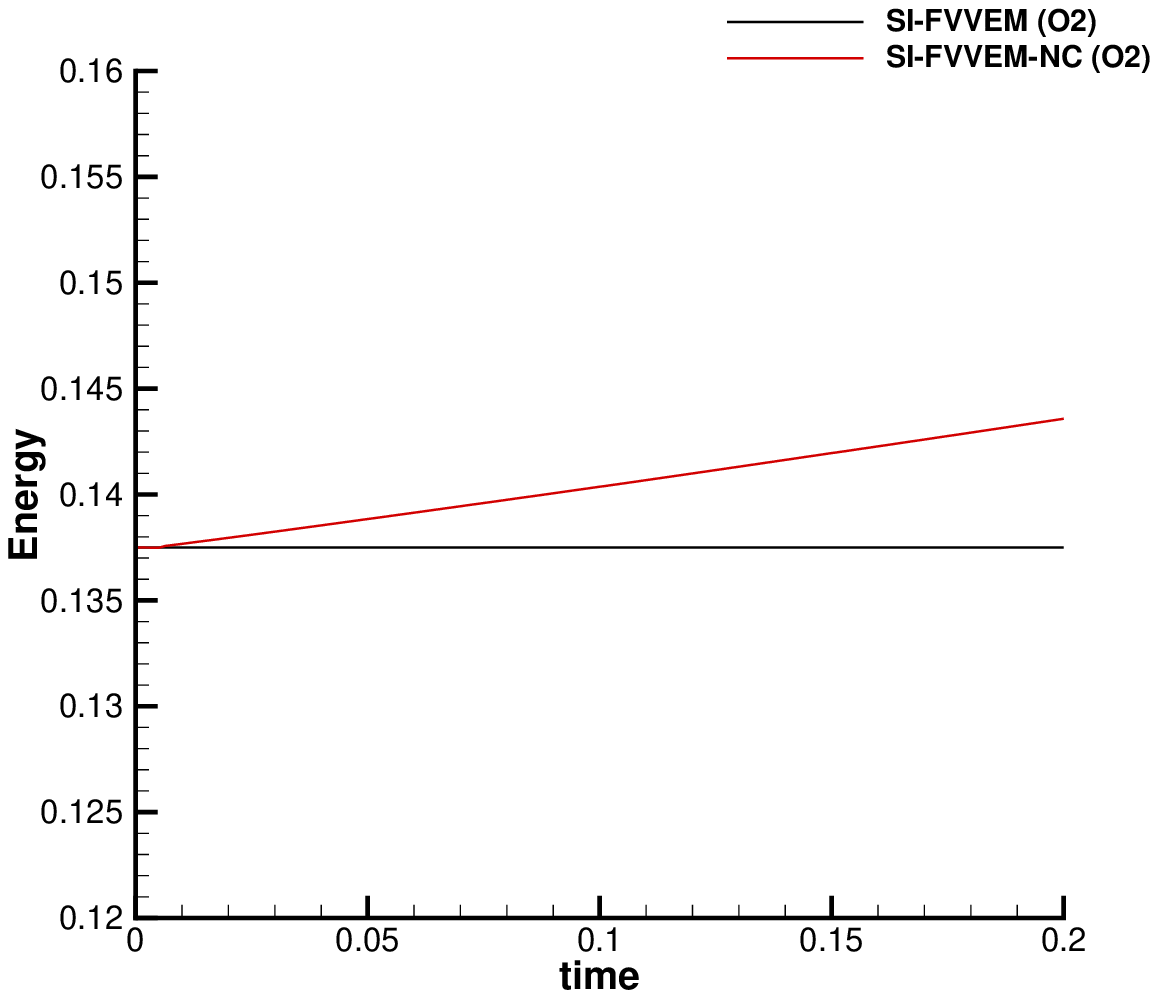}  &
			\includegraphics[width=0.47\textwidth]{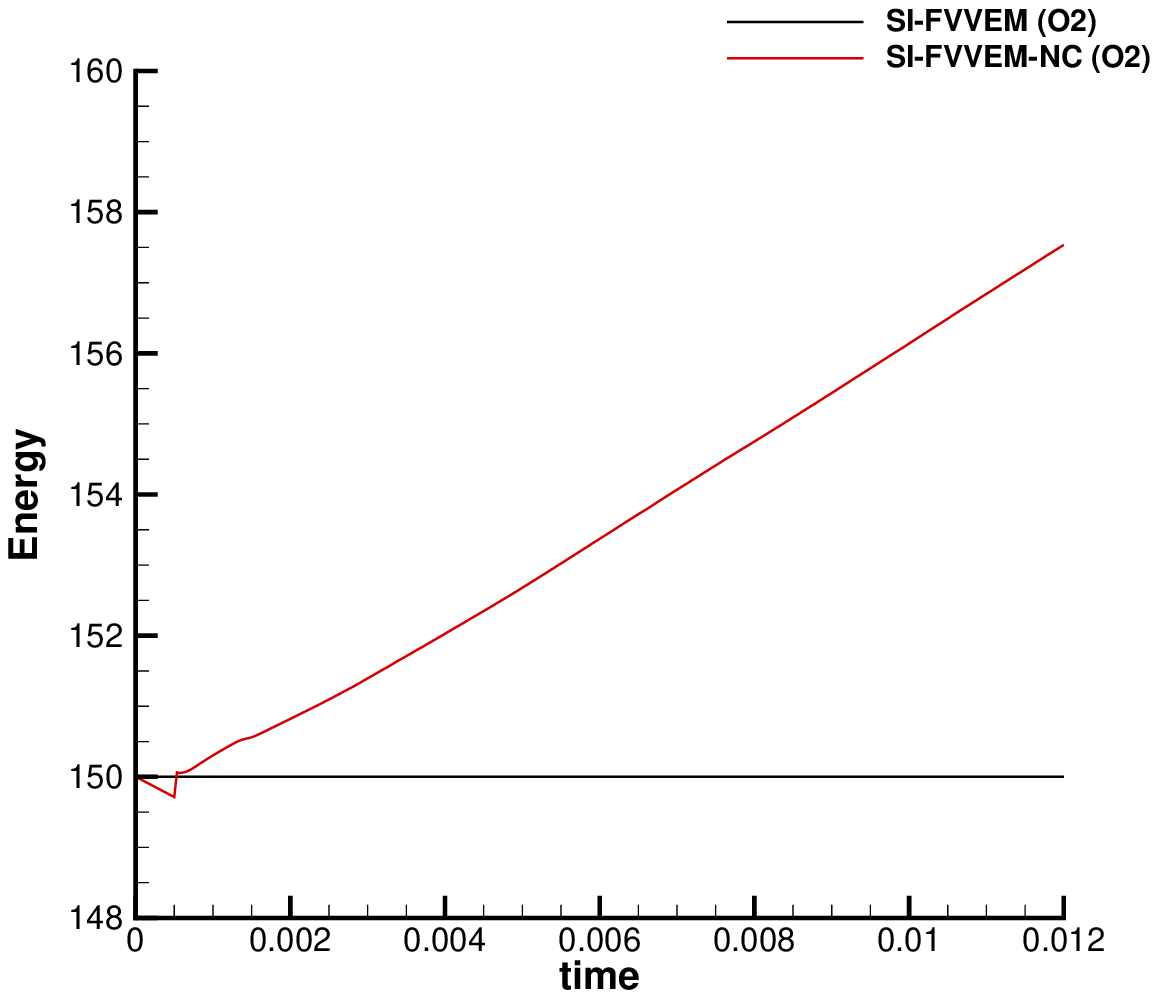}  \\
		\end{tabular}
		\caption{Riemann problems. Time evolution of the total energy for the globally energy conservative (black line) and the non-conservative (red line) SI-FVVEM schemes for RP1 (left) and RP4 (right).}
		\label{fig.RPenergy}
	\end{center}
\end{figure}

To illustrate the global energy conservation property of the proposed scheme proven in Theorem \ref{th.energy}, we report in Figure~\ref{fig.RPenergy} the evolution of the global total energy along the simulation of RP1 and RP4. We observe that the new SI-FVVEM method is able to globally preserve the energy up to machine precision. On the contrary, avoiding the two step procedure, i.e. not computing the solution of \eqref{eqn.pressure_vf4}, would lead to a non globally conservative scheme where the energy linearly increases in time.

\subsection{Circular explosion}
\begin{figure}[!htbp]
	\begin{center}
		\begin{tabular}{cc}  
			\includegraphics[trim= 5 5 5 5, clip, width=0.47\textwidth]{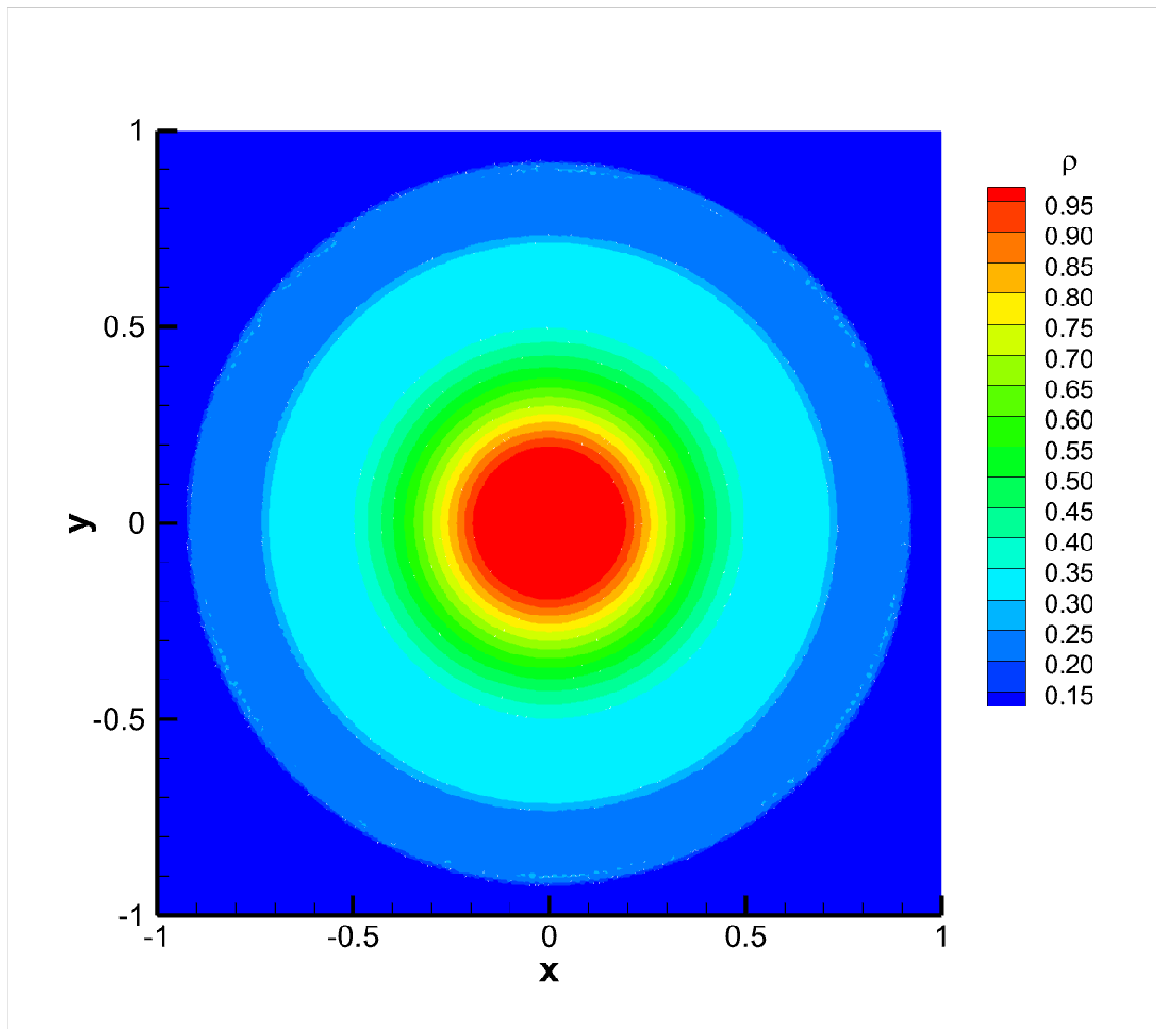}&
			\includegraphics[width=0.47\textwidth]{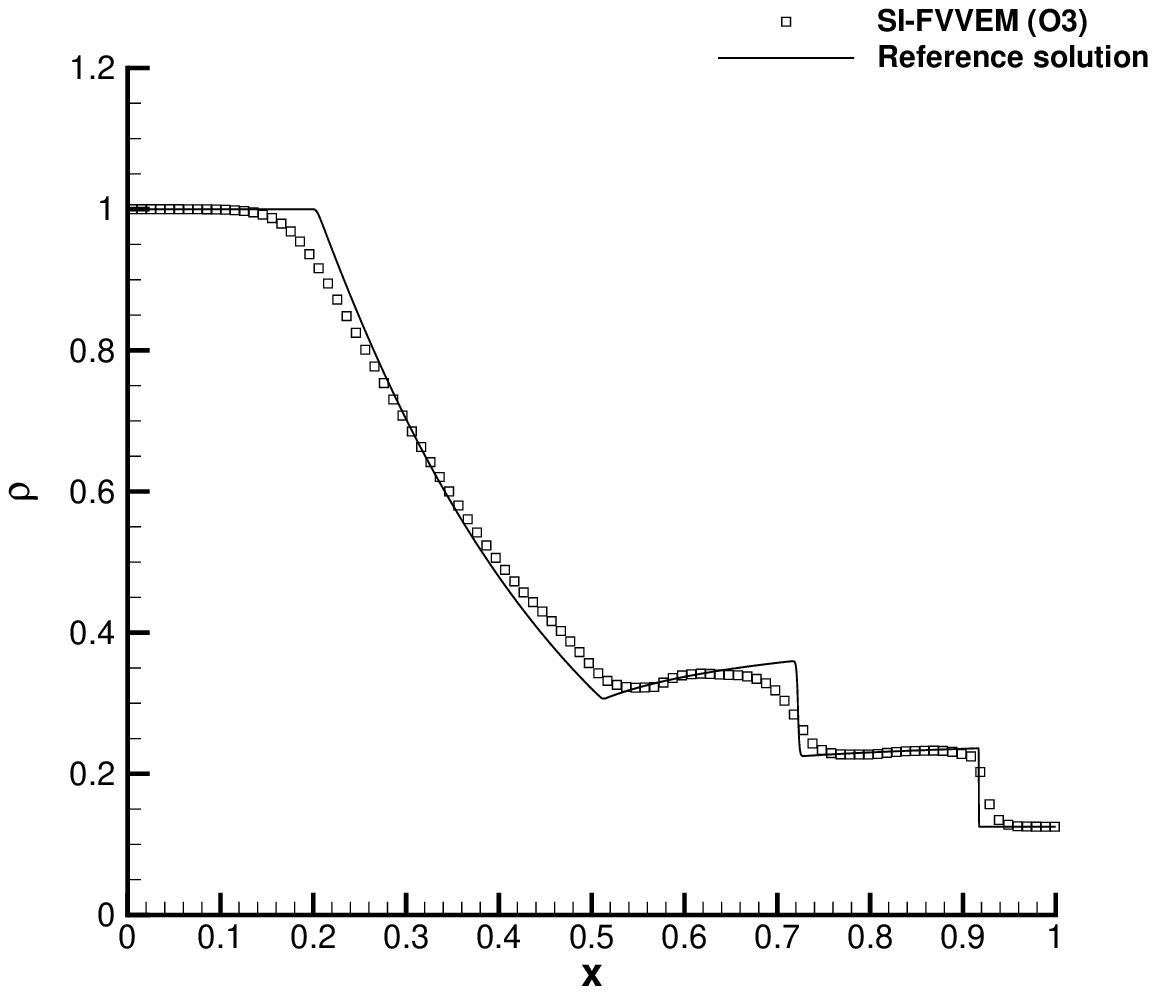} \\				
			\includegraphics[width=0.47\textwidth]{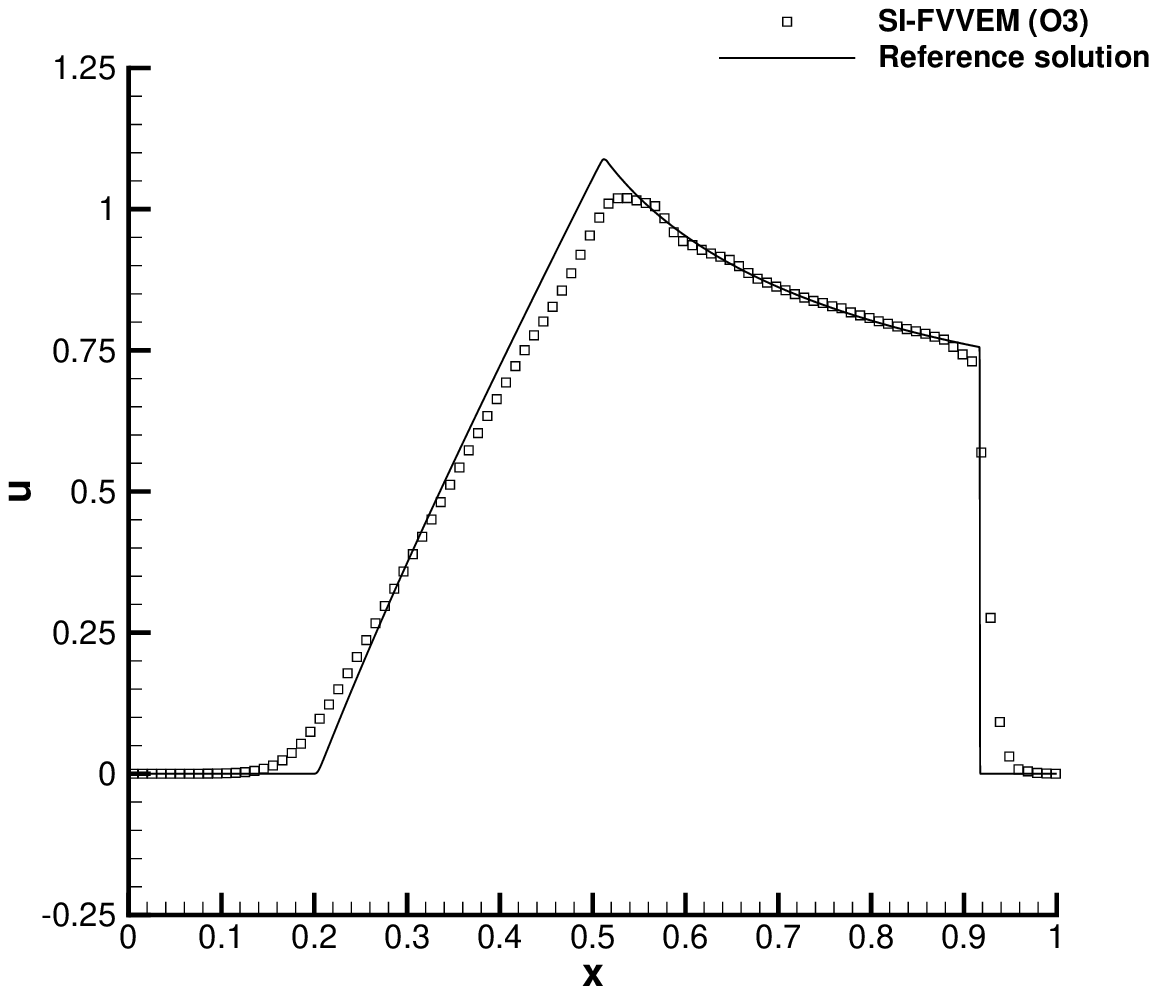} & 
			\includegraphics[width=0.47\textwidth]{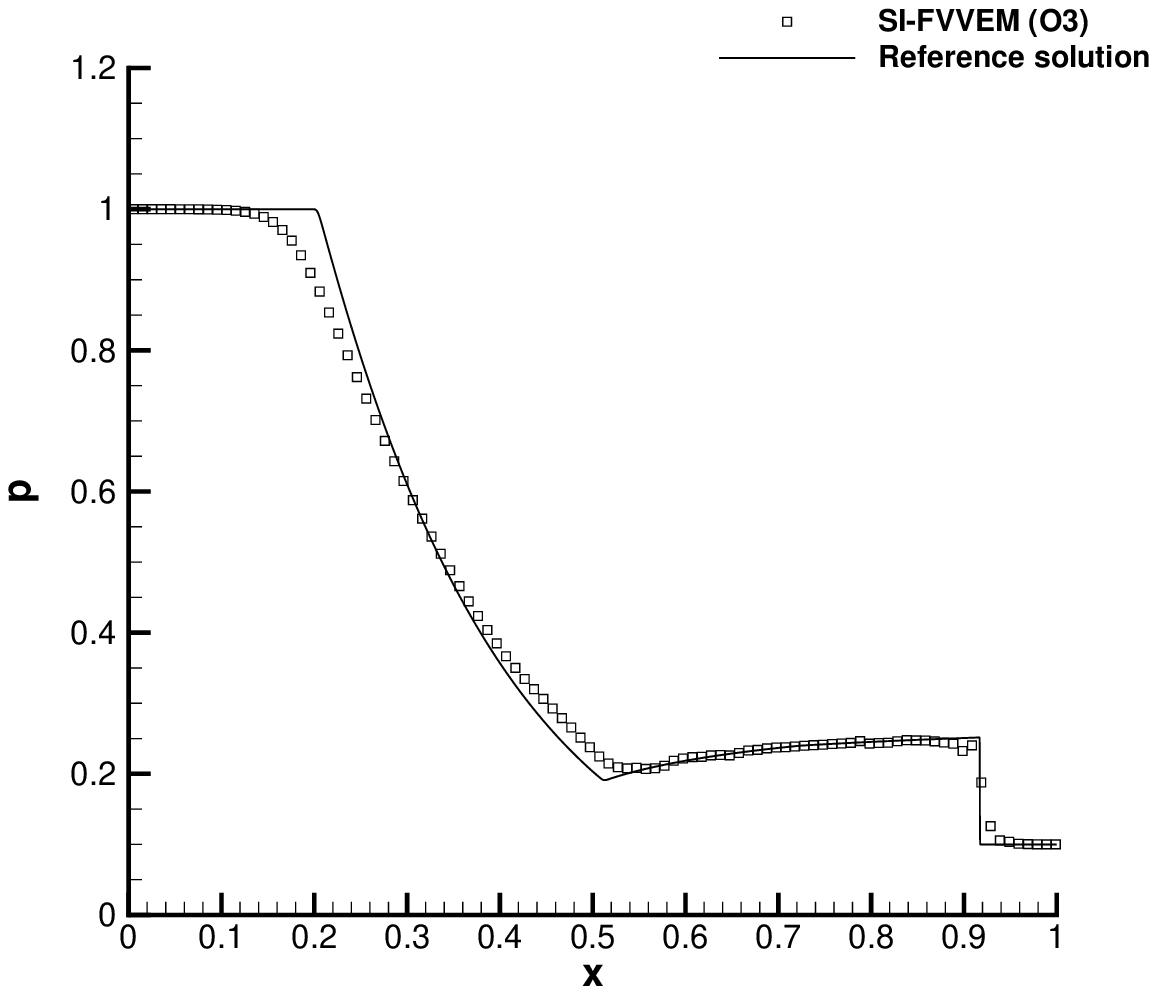} \\
		\end{tabular} 
		\caption{Circular explosion problem at time $t_f=0.25$. Top left: two-dimensional view of the density distribution. From top-left to bottom-right: density $\rho$, horizontal velocity $u$ and pressure $p$ distribution compared against the reference solution extracted with a one-dimensional cut of 200 equidistant points along the $x-$direction at $y=0$. }
		\label{fig.EP2D}
	\end{center}
\end{figure}

We consider a classical circular explosion problem with a radial initial condition based on the Sod shock tube benchmark:
\begin{equation}
	\rho\left(\xx,0\right) = \left\lbrace \begin{array}{lr}
		1 & \mathrm{if}\; r\leq 0.5\\
		0.125 & \mathrm{if}\; r> 0.5\\
	\end{array} \right., \qquad
	\press\left(\xx,0\right) = \left\lbrace \begin{array}{lr}
		1 & \mathrm{if}\; r\leq 0.5\\
		0.1 & \mathrm{if}\; r> 0.5\\
	\end{array} \right., \qquad \vel\left(\xx,0\right) = 0,
\end{equation}
with the generic radial position given by $r =\sqrt{x^2+y^2}$.
The computational domain $\Omega = [-1,1]^2$ is discretized using a Voronoi mesh of $N=25648$ 
elements and the initial condition is adopted to set Dirichlet boundary conditions everywhere. The results obtained with the third order SI-FVVEM scheme at $t_f=0.25$ are shown in Figure~\ref{fig.EP2D}. We observe an excellent agreement with the reference solution obtained using a second order TVD scheme to solve the corresponding one-dimensional PDE in the radial direction obtained from the compressible Euler equations with geometrical source terms \cite{Toro,Hybrid2}.

\subsection{Double Mach reflection problem}
\begin{figure}[!htbp]
	\begin{center}
		\begin{tabular}{cc} 
			\multicolumn{2}{c}{\includegraphics[trim= 2 2 2 2, clip,width=0.9\textwidth]{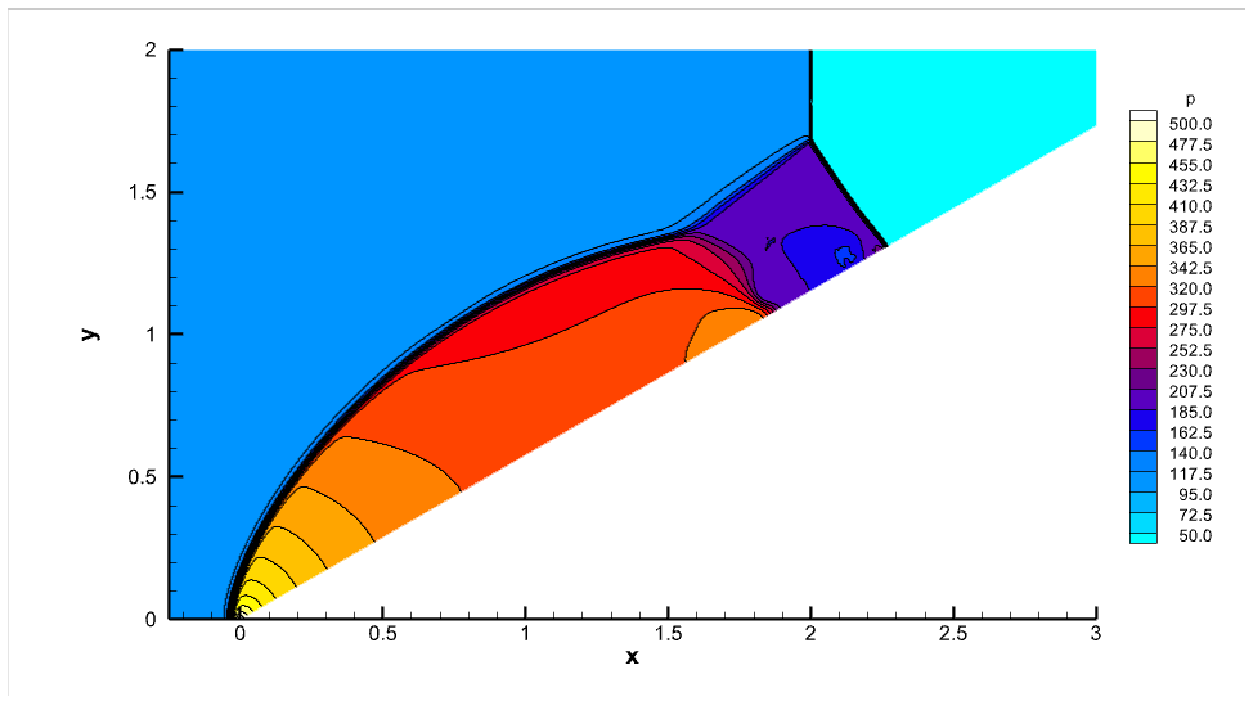}} \\
			\includegraphics[trim= 5 5 5 5, clip,width=0.47\textwidth]{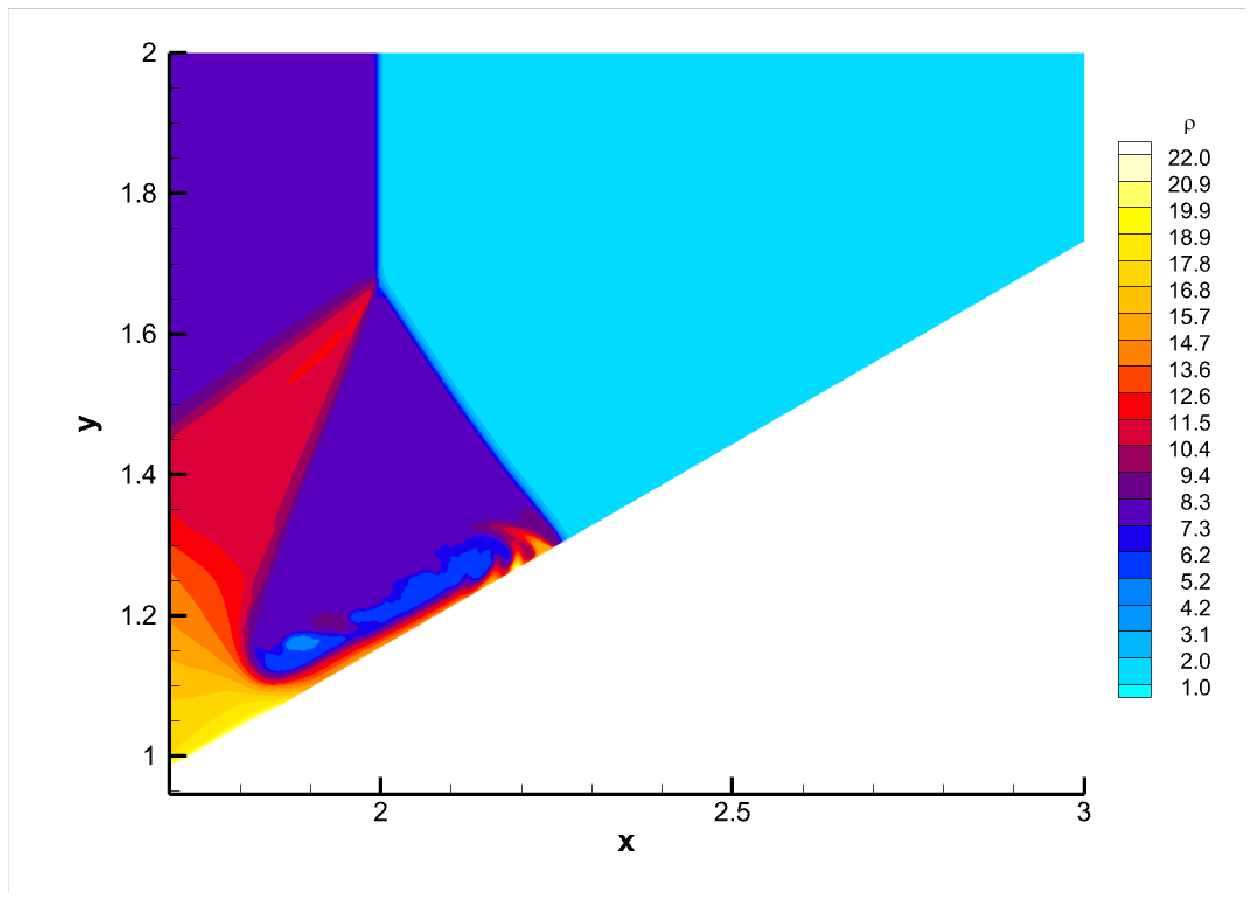} & \includegraphics[trim= 5 5 5 5, clip,width=0.47\textwidth]{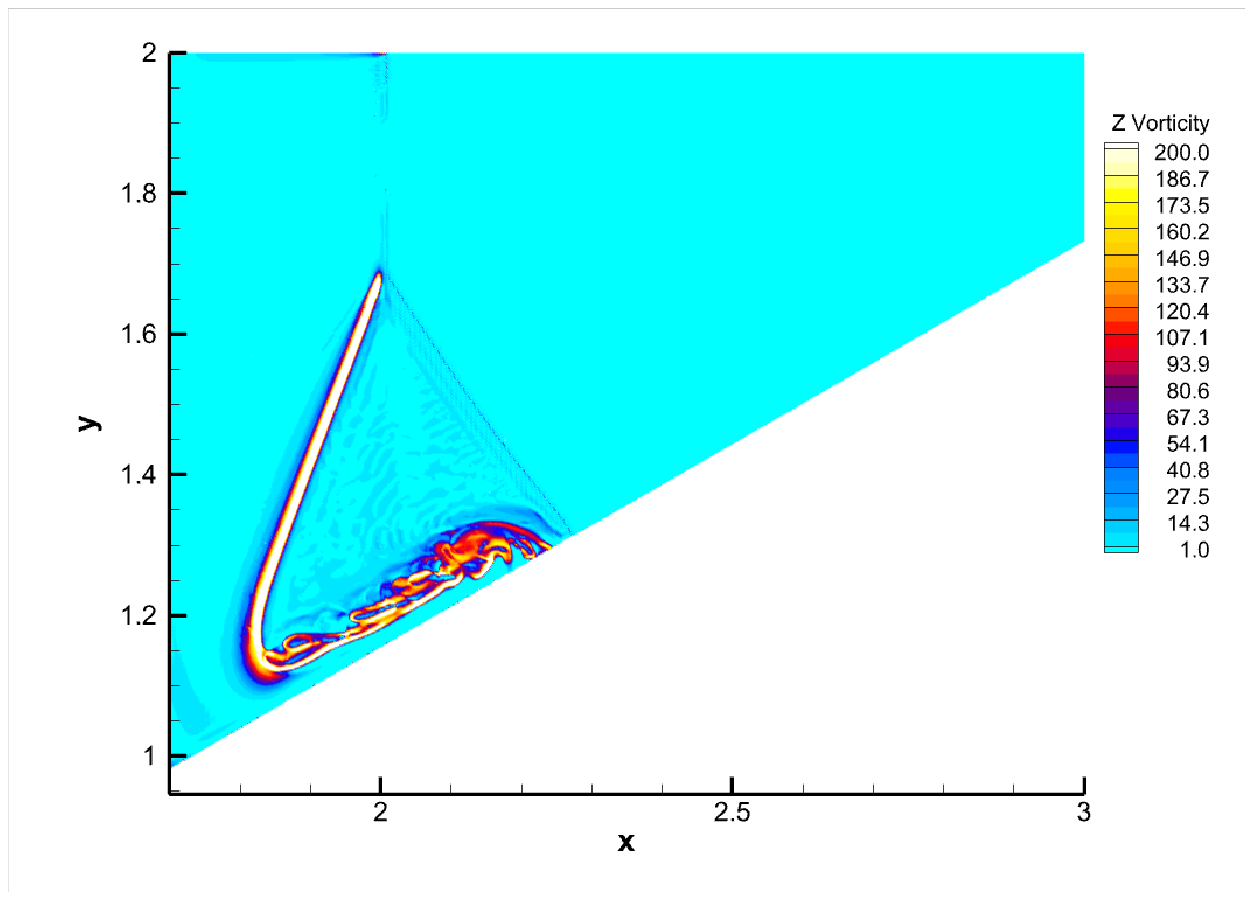}
		\end{tabular} 
		\caption{Double Mach reflection problem at time $t_f=0.2$. Top: 21 equidistant contour lines in the range $[50,500]$ for pressure. Bottom: zoom on the shock front with density (left) and vorticity (right) distribution.} 
		\label{fig.DMR}
	\end{center}
\end{figure}

To demonstrate the effectiveness of the new SI-FVVEM schemes in the high Mach number regime, we run this challenging test case. A very strong shock wave, which is initially located at $x_0=0$, moves along the $x-$direction of the computational domain, where a ramp with angle $\alpha=\frac{\pi}{6}$ is located. The computational domain is the box $\Omega=[-0.25,3]\times[0,2]$, from which the angle $\alpha$ is subtracted. It is discretized with $N=1429782$ polygonal control volumes of characteristic mesh size $h=1/400$. Inflow boundaries are set on the left side, outflow conditions on the right side, and wall boundary conditions are imposed elsewhere. The shock Mach number is $\text{M}_s=10$ and the initial condition is given in terms of primitive variables according to \cite{WC84}:
\begin{equation}
	(\rho(\xx,0),u(\xx,0),v(\xx,0),p(\xx,0)) = \left\{ \begin{array}{lll} \left( 8.0, 8.25, 0, 116.5 \right), & \textnormal{ if } & x \leq x_0, \\ 
		\left( 1.4,0,0,1.0 \right), & \textnormal{ if } & x > x_0  . 
	\end{array}  \right.
	\label{eqn.DMR_IC}
\end{equation}
The simulation stops at time $t_f=0.2$, and it is run with MPI parallelization on 94 CPUs. The computation is run using the second order SI-FVVEM scheme with an artificial viscosity coefficient of $c_{\alpha}=10^{-3}$. The results are depicted in Figure~\ref{fig.DMR}. The flow field agrees well with the numerical reference solutions shown in \cite{WC84,BB19}. Thanks to the energy conservation property of the new schemes, the shock wave is properly resolved and located in the correct position at $x=2$. 

\subsection{Taylor-Green vortex}
\begin{figure}[!htbp]
	\begin{center}
		\begin{tabular}{cc} 
			\includegraphics[trim= 5 5 5 5, clip, width=0.47\textwidth]{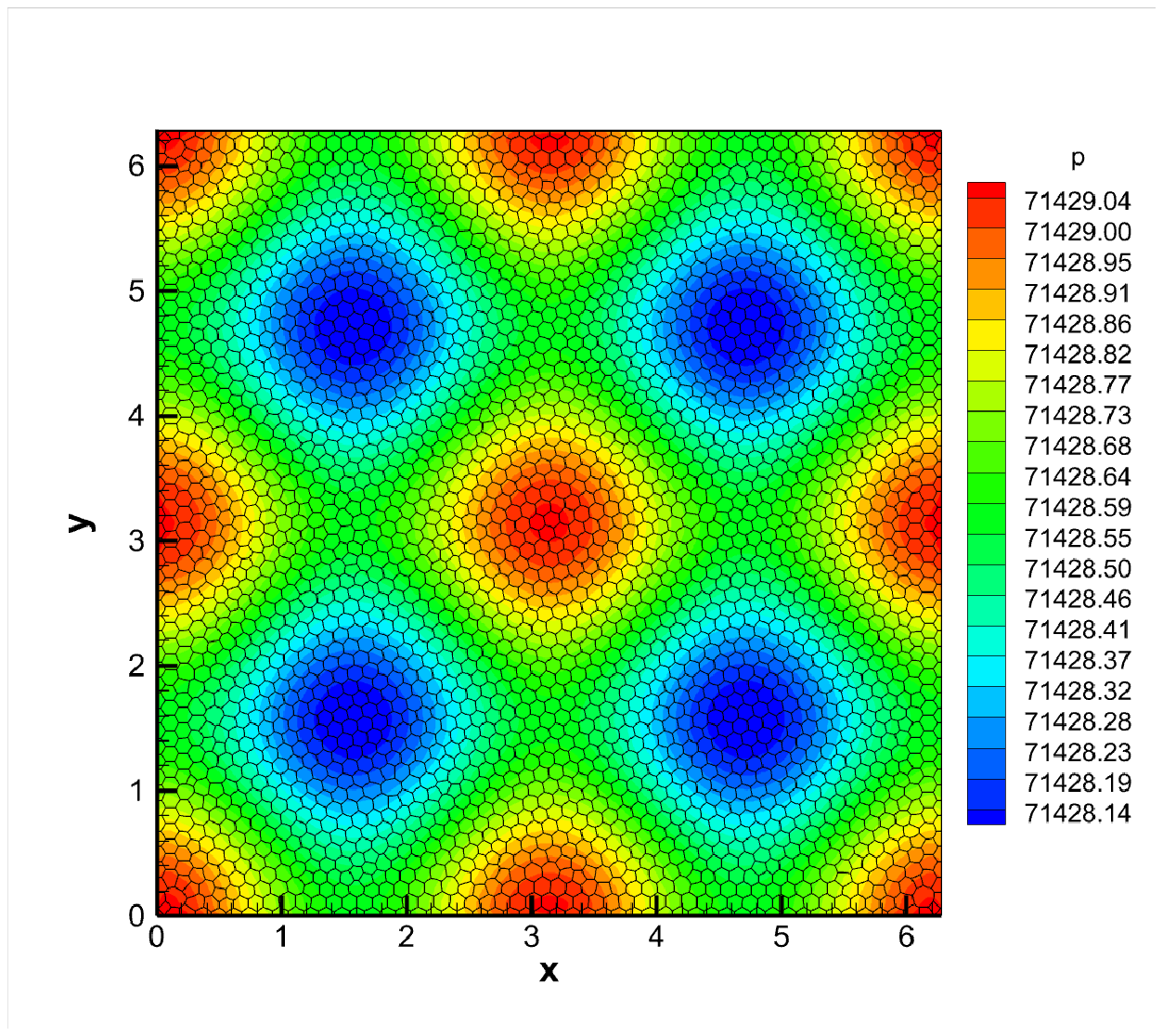} &  			
			\includegraphics[trim= 5 5 5 5, clip, width=0.47\textwidth]{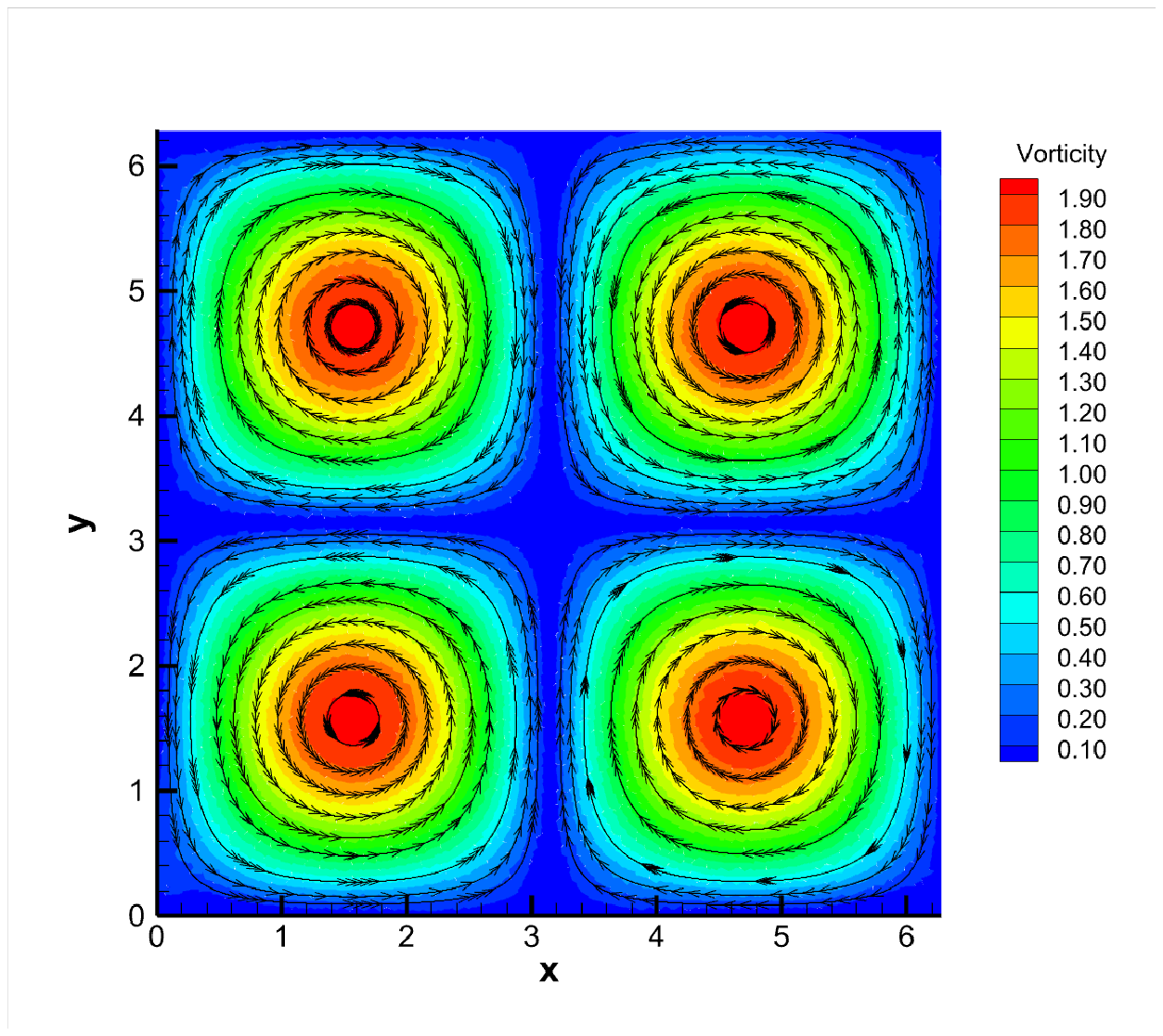} \\
			\includegraphics[width=0.47\textwidth]{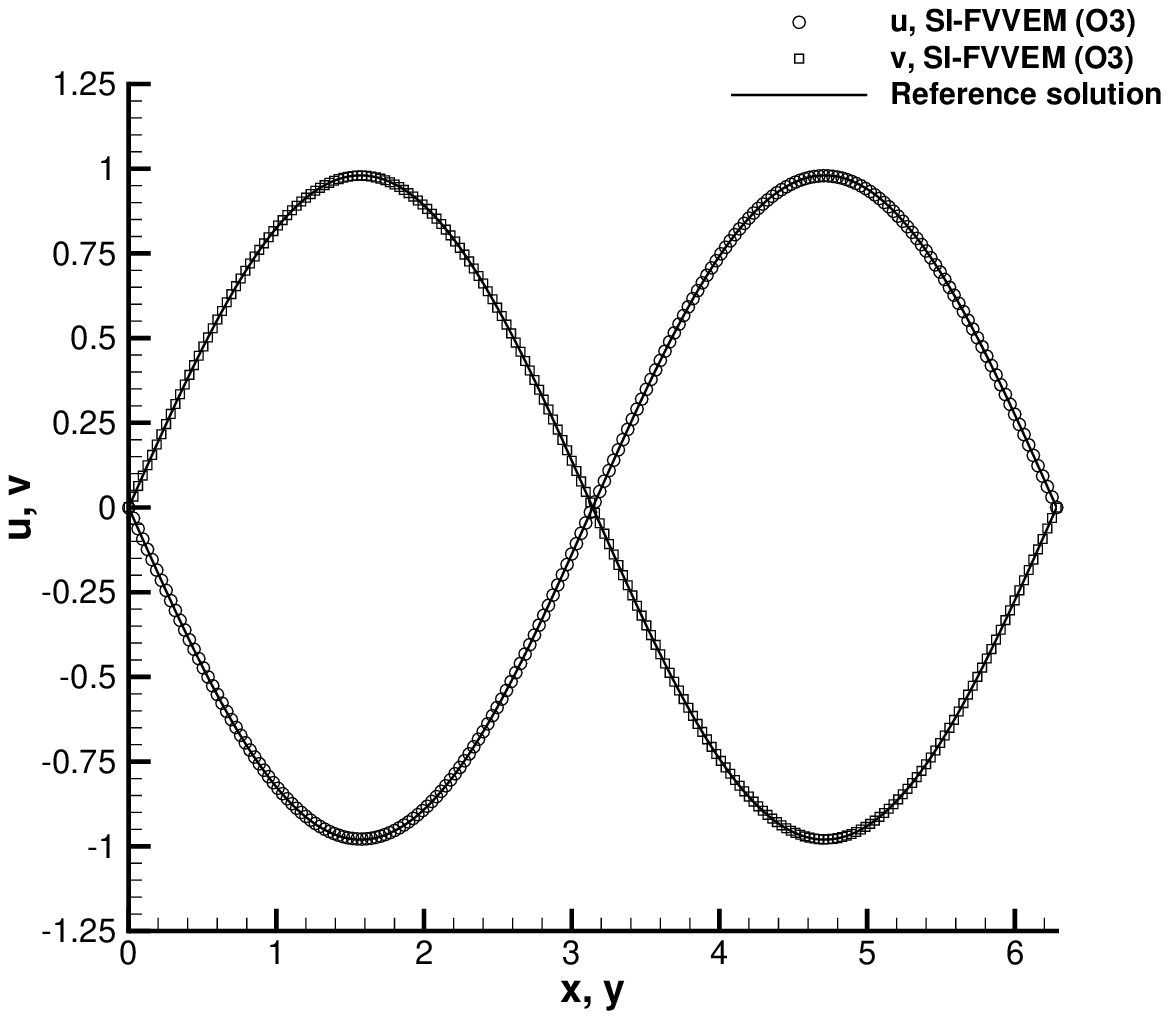} &  			
			\includegraphics[width=0.47\textwidth]{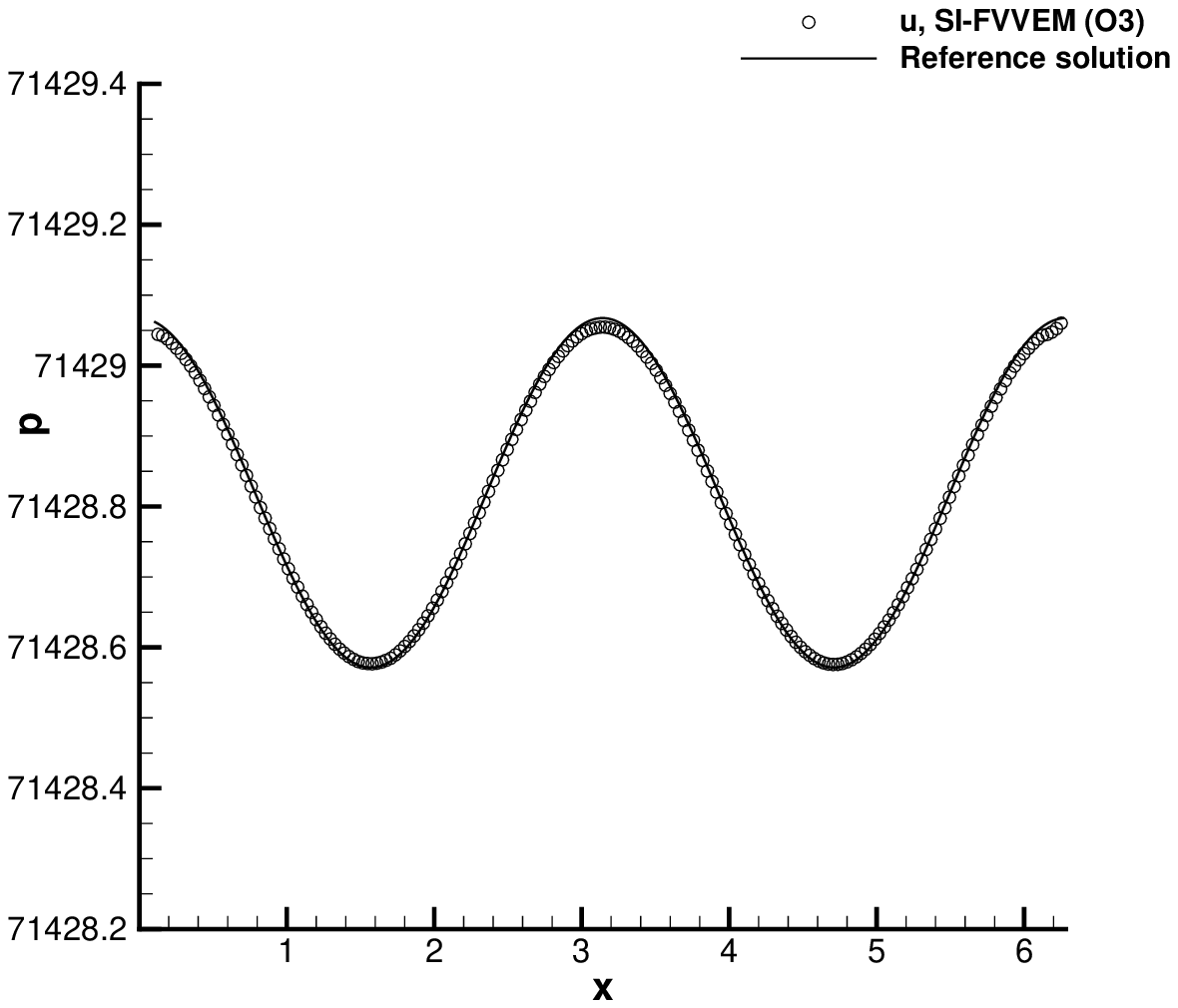} \\
		\end{tabular}
		\caption{Taylor-Green vortex at time $t_f=0.2$ with viscosity $\nu = 10^{-2}$ and pressure $p_{0}=10^{5}$. Top: pressure (left) and vorticity magnitude with the stream-traces (right). Bottom: one-dimensional cuts with 200 equidistant points along the $x-$axis and the $y-$axis for velocity components $u,v$ (left) and along the $x-$axis for the pressure $p$ with comparison against the exact solution.}
		\label{fig.TGV_mu1e-2}
	\end{center}
\end{figure}

The behavior of the SI-FVVEM approach in the framework of viscous flows is first analyzed using a modified Taylor-Green vortex problem whose initial condition is given by
\begin{equation}\label{eqn.tgv}
	\rho\left(\mathbf{x},t\right) = 1,\quad
	\vel \left(\mathbf{x},t\right) = \left( \begin{array}{r} 
		\sin(x)\cos(y) e^{-2\mu t} \\ 
		-\cos(x)\sin(y)  e^{-2\mu t}\end{array} \right), \quad 
	\press \left(\mathbf{x},0\right) = \frac{\press_{0}}{\gamma-1} + \frac{1}{4} \left(\cos(2x)+\cos(2y) \right) e^{-4\mu t}. 
\end{equation}
The simulations are run in the computational domain $\Omega=[0,2\pi]^2$ up to time $t_f=0.2$ using non-slip wall boundary conditions everywhere. As an example, one of the Voronoi grids employed, made of $N=3221$ elements, is depicted in Figure~\ref{fig.TGV_mu1e-2} together with the pressure contour plot for $\nu=10^{-2}$ and $p_0=10^{2}$. To ease the comparison with reference data in the bibliography, also the contour plot of the vorticity is shown. The 1D cuts of the the velocity components and pressure field match well with the reference solution, as observed in Figure~\ref{fig.TGV_mu1e-2}. 
Furthermore, the $L^{2}$ errors and convergence rates obtained using the second order scheme for $\nu=10^{-2}$ and $\nu=10^{-5}$ are reported in Table~\ref{tab.conv_rate_tgv}. The obtained results illustrate the asymptotic preserving property of the scheme with respect to the Reynolds number.

\begin{table}
	\begin{center} 
		\caption{Numerical convergence results of the SI-FVVEM scheme with second order of accuracy in space and time using the Taylor-Green vortex with $p_0=10^{2}$ on general polygonal meshes. The errors are measured in $L_2$ norm and refer to the pressure $p$ and velocity component $u$ at time $t_f=0.2$. The asymptotic preserving (AP) property of the scheme is studied by considering different Reynolds numbers $\Rey=\{ 10^{2}, 10^{5} \}$.
		}
		\begin{small}
			\renewcommand{\arraystretch}{1.1}	
			\begin{tabular}{ccccc}
				$\cvs(\Omega)$ & ${L_2}(p)$ & $\mathcal{O}(p)$ & ${L_2}(u)$ & $\mathcal{O}(u)$ \\
				\hline
				& \multicolumn{4}{c}{$\Rey=10^2$} \\
				\hline
				$\frac{1}{8}$ & $6.3271 \cdot 10^{-1}$ &  -      & $1.6742 \cdot 10^{0}$ & - \\
				$\frac{1}{15}$& $1.3686 \cdot 10^{-1}$ & $ 2.44$ & $3.5884 \cdot 10^{-1}$ & $ 2.45$ \\
				$\frac{1}{30}$& $2.7497 \cdot 10^{-2}$ & $ 2.32$ & $1.0891 \cdot 10^{-1}$ & $ 1.72$ \\
				$\frac{1}{45}$& $1.1556 \cdot 10^{-2}$ & $ 2.14$ & $4.7272 \cdot 10^{-2}$ & $ 2.06$ \\
				\hline
				& \multicolumn{4}{c}{$\Rey=10^5$} \\
				\hline
				$\frac{1}{8}$ & $6.3058 \cdot 10^{-1}$ & -       & $1.6859 \cdot 10^{0} $ & - \\
				$\frac{1}{15}$& $1.3649 \cdot 10^{-1}$ & $2.43 $ & $3.7171 \cdot 10^{-1}$ & $2.41 $ \\
				$\frac{1}{30}$& $2.7375 \cdot 10^{-2}$ & $2.32 $ & $1.1395 \cdot 10^{-1}$ & $1.71 $ \\
				$\frac{1}{45}$& $1.2662 \cdot 10^{-2}$ & $1.90 $ & $5.0073 \cdot 10^{-2}$ & $2.03 $ \\
				\hline
			\end{tabular}
		\end{small}
		\label{tab.conv_rate_tgv}
	\end{center}
\end{table}

\subsection{First problem of Stokes}
\begin{figure}[!htbp]
	\centering
	\includegraphics[width=0.32\linewidth]{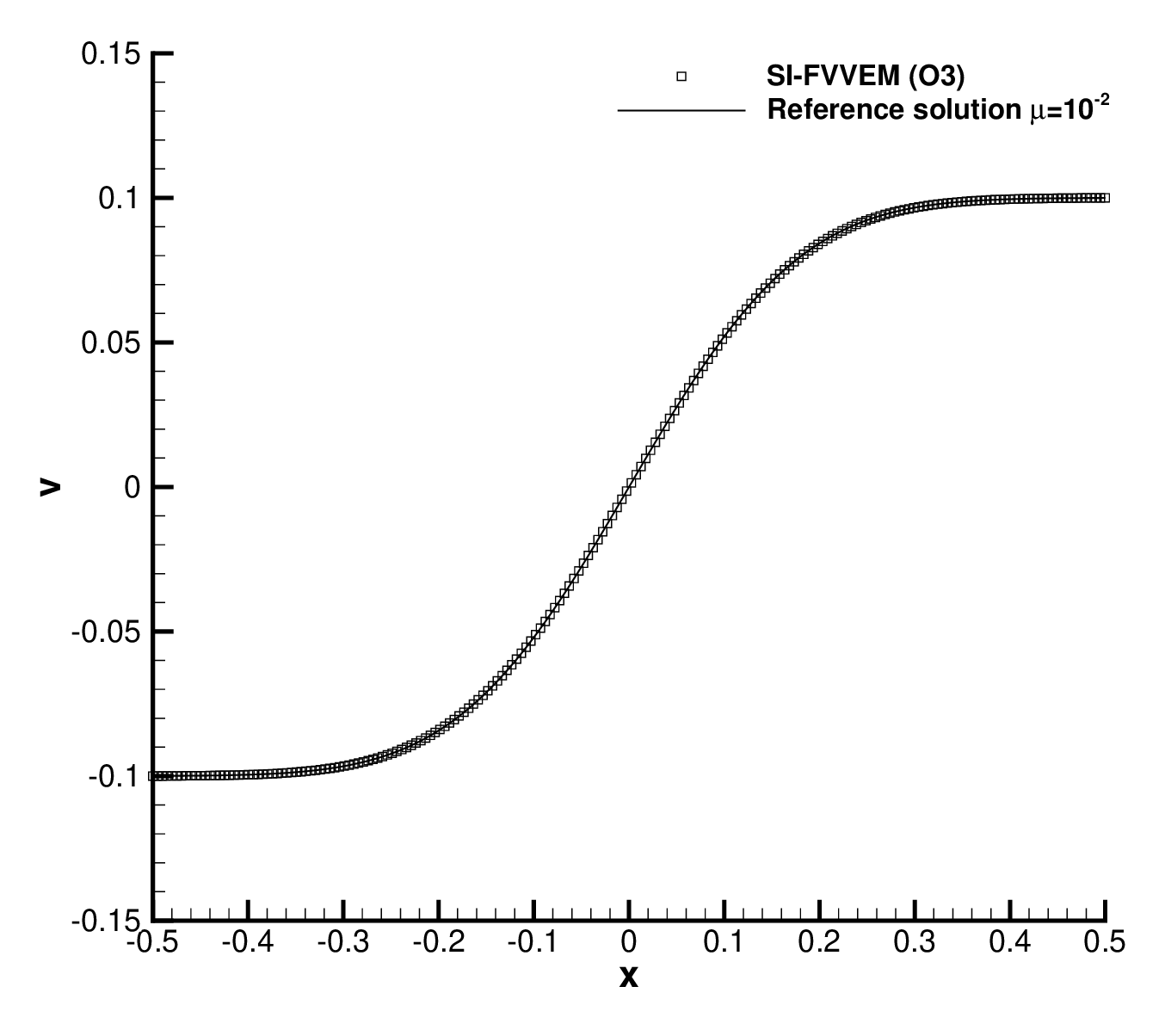}\hfill
	\includegraphics[width=0.32\linewidth]{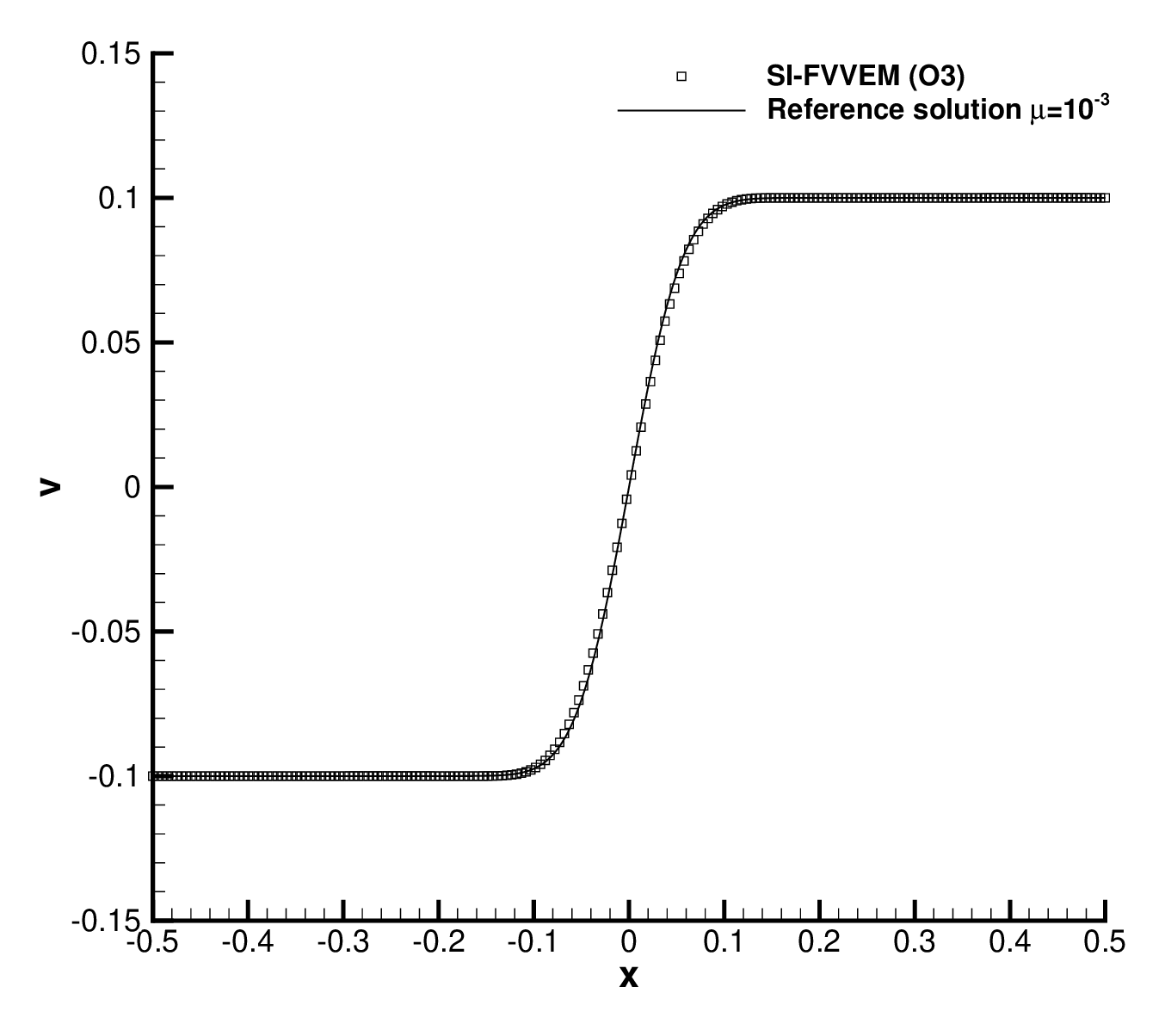}\hfill
	\includegraphics[width=0.32\linewidth]{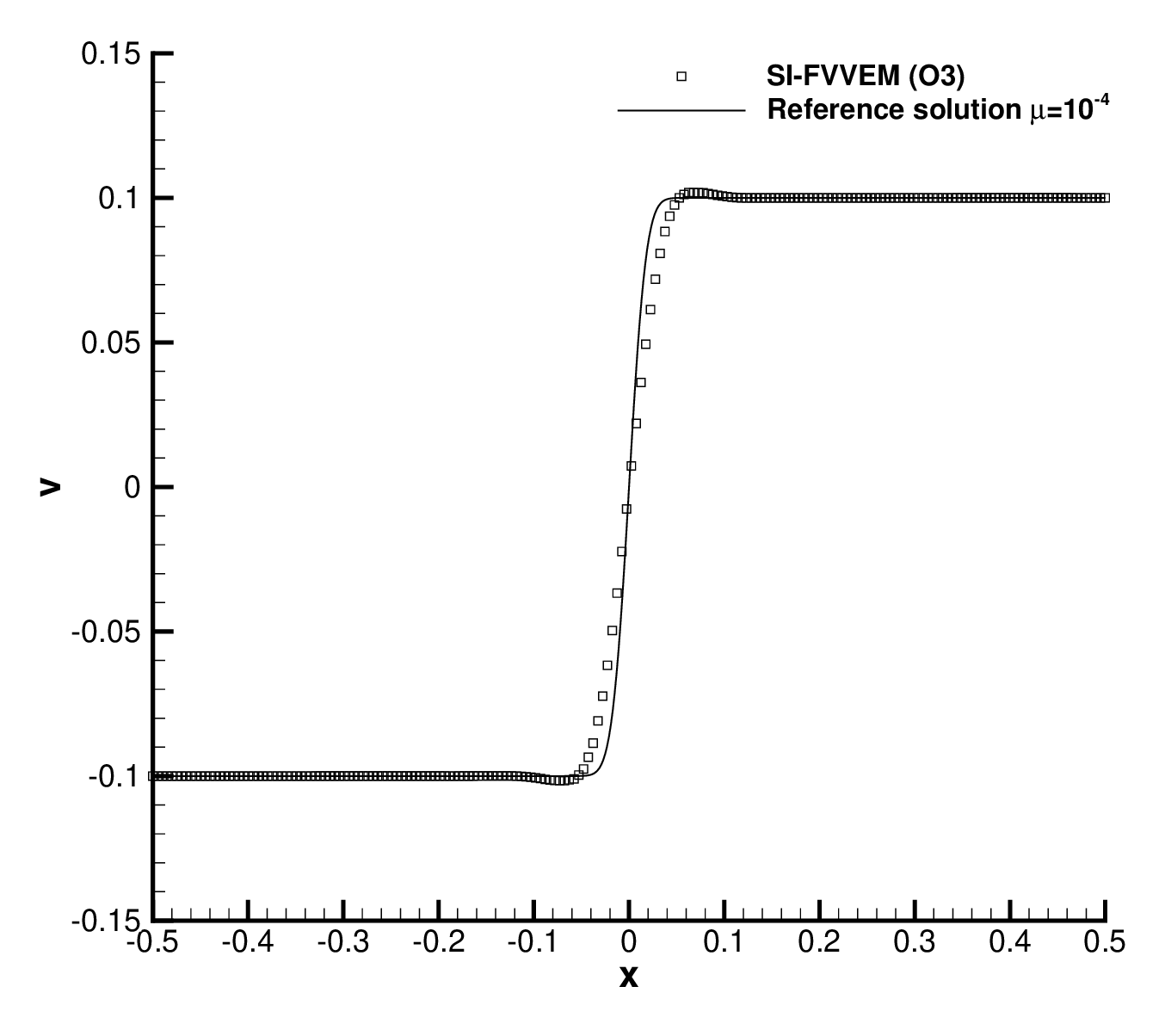}
	\caption{First problem of Stokes at time $t_f=1$. One-dimensional plot of the velocity component $v$  with $\mu=10^{-2}$ (left), $\mu=10^{-3}$ (middle), $\mu=10^{-4}$ (right).}
	\label{fig:FS2D}
\end{figure}

The first problem of Stokes is a classical benchmark for the validation of incompressible viscous flow solvers. It consists in an initial flow field of the form
\begin{equation}
	\rho\left(\mathbf{x},0\right) = 1,\qquad
	\press \left(\mathbf{x},0\right) = \frac{1}{\gamma}, \qquad
	u \left(\mathbf{x},0\right) = 0, \qquad
	v \left(\mathbf{x},0\right) = \left\lbrace \begin{array}{rr}
		-0.1 & \mathrm{ if } \; y \le 0,\\
		0.1 & \mathrm{ if } \; y > 0,
	\end{array}\right.
\end{equation}
which, in the incompressible limit, has an analytical solution for the vertical velocity component $v$ given by
\begin{equation}
	v \left(\mathbf{x},t\right) = \frac{1}{10} \mathrm{erf}\left( \frac{x}{2\sqrt{\mu t}}\right),
\end{equation}
with $\mathrm{erf}(z)$ the error function.
For the numerical simulations, we consider the computational domain \mbox{$\Omega=[-0.5,0.5]\times[-0.2,0.2]$} discretized using a polygonal grid generated employing $N_x=100$ and $N_y=10$ points in the $x-$ and $y-$ directions, respectively. Periodic boundary conditions are considered in the $y-$direction while the exact solution is imposed on the left and right boundaries. In Figure~\ref{fig:FS2D} we compare the results obtained using the second order SI-FVVEM scheme against the reference solution at time $t_f=1$.  A good agreement is observed for the three viscosity coefficients considered, namely $\nu=\{ 10^{-2}, 10^{-3}, 10^{-4}\}$.
 
\subsection{Viscous shock}
\begin{figure}[!htbp]
	\begin{center} 
		\includegraphics[trim=2 2 2 2,clip,width=0.9\textwidth]{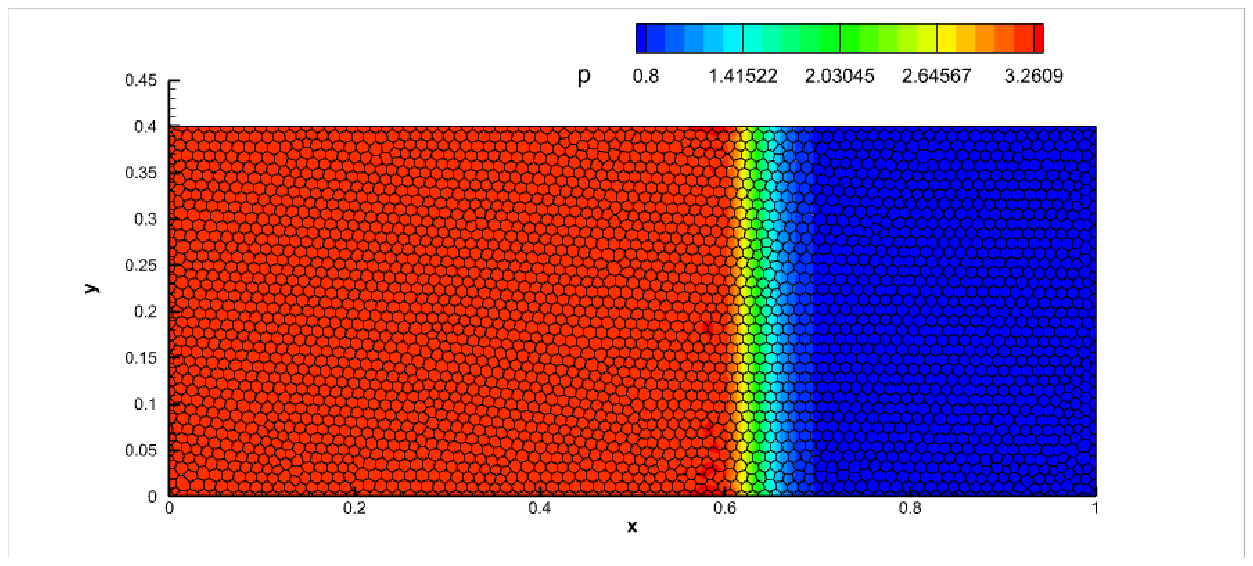}
		\caption{Viscous shock profile with shock Mach number $\mathrm{M_s}=2$ and Prandtl number $\Pr=0.75$ at time $t_f=0.2$. Voronoi tessellation and pressure contours.}
		\label{fig.ViscousShockcont}
	\end{center}
\end{figure}

\begin{figure}[!htbp]
	\begin{center}
		\begin{tabular}{cc} 
			\includegraphics[width=0.47\textwidth]{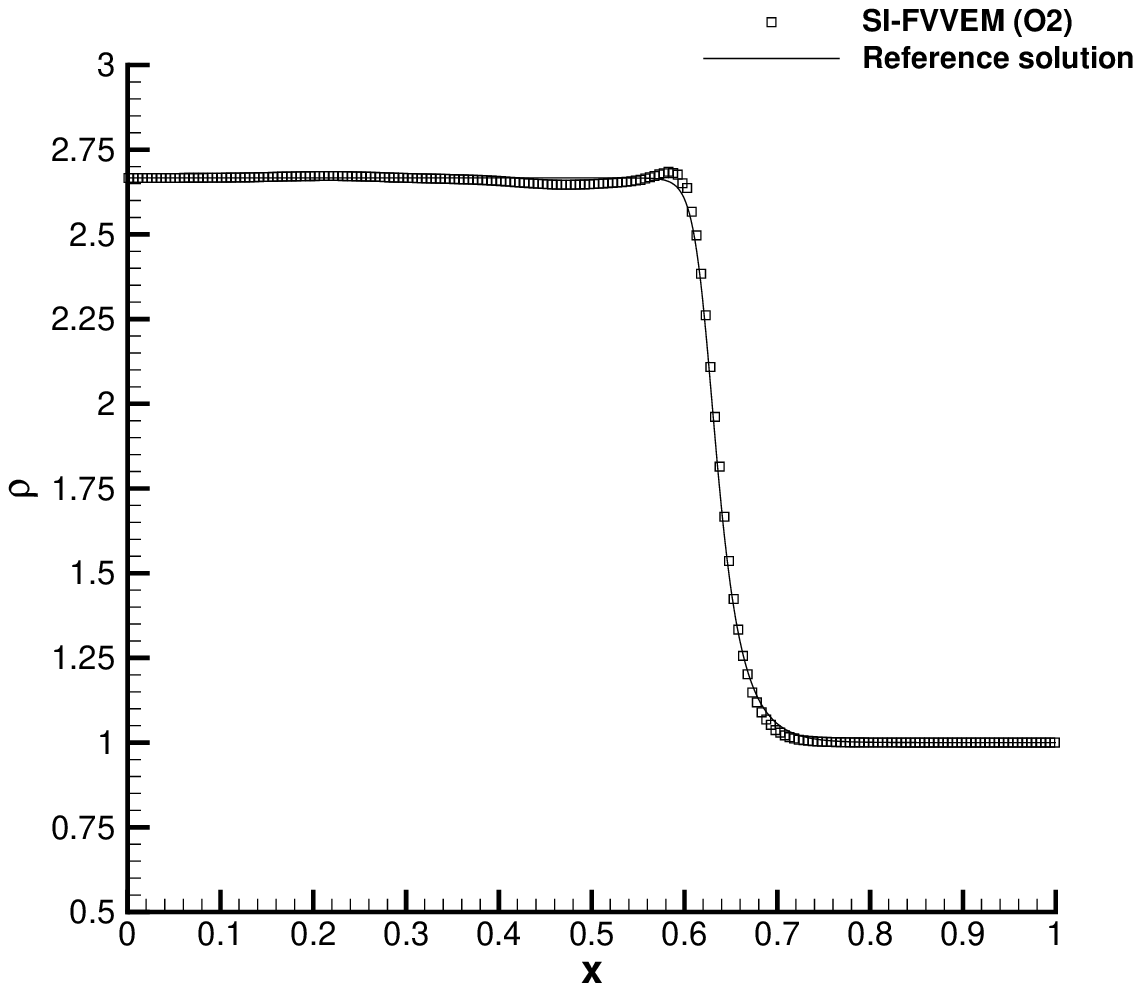} & 
			\includegraphics[width=0.47\textwidth]{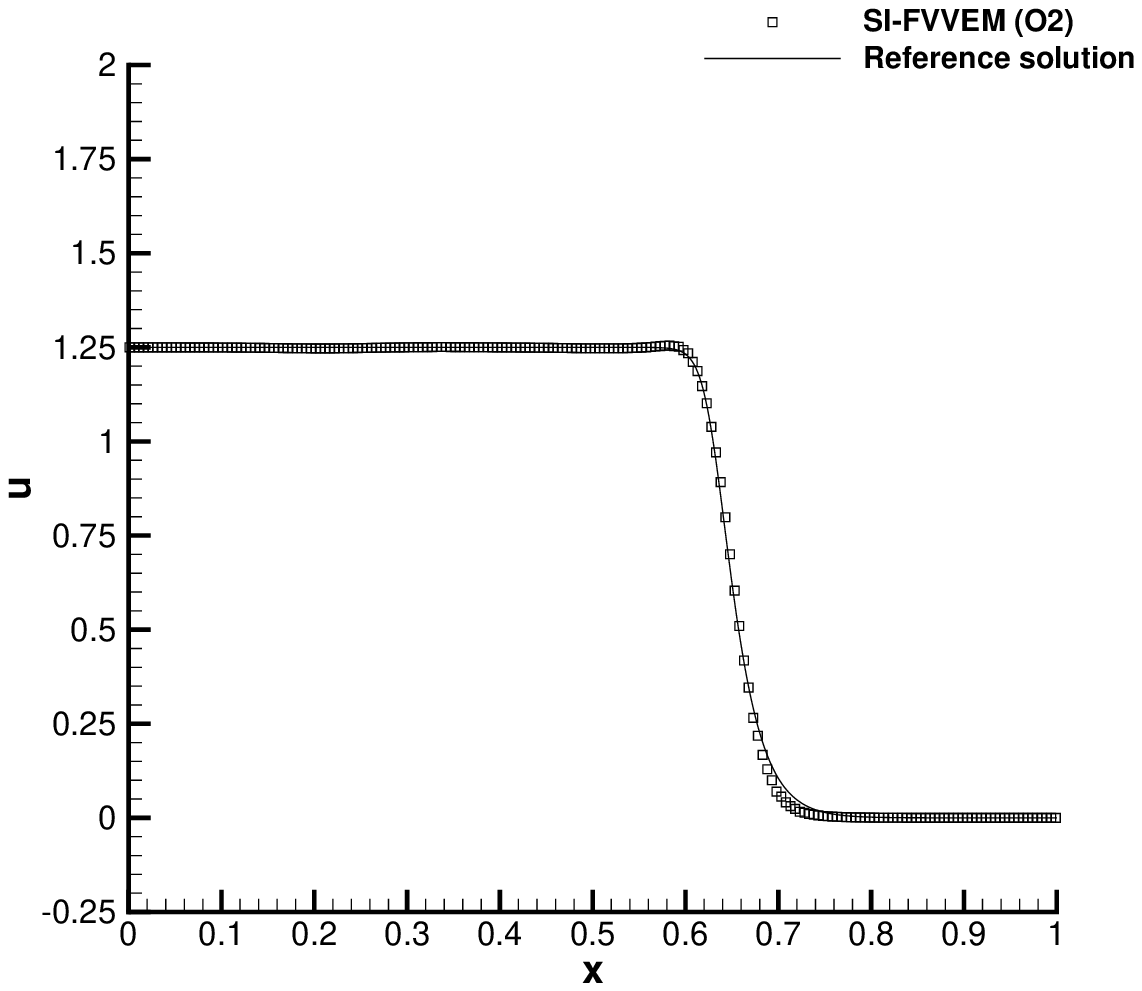} \\
			\includegraphics[width=0.47\textwidth]{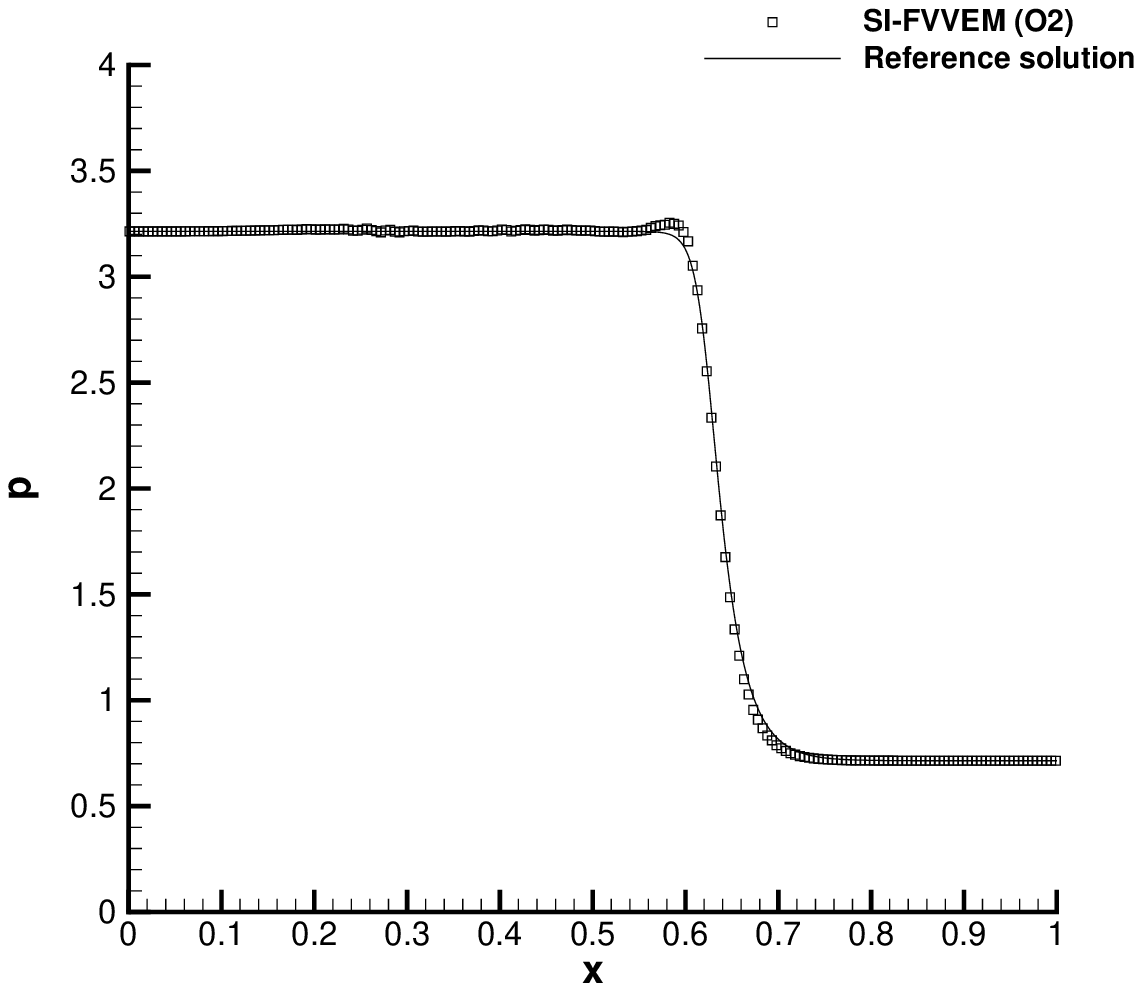} & 
			\includegraphics[width=0.47\textwidth]{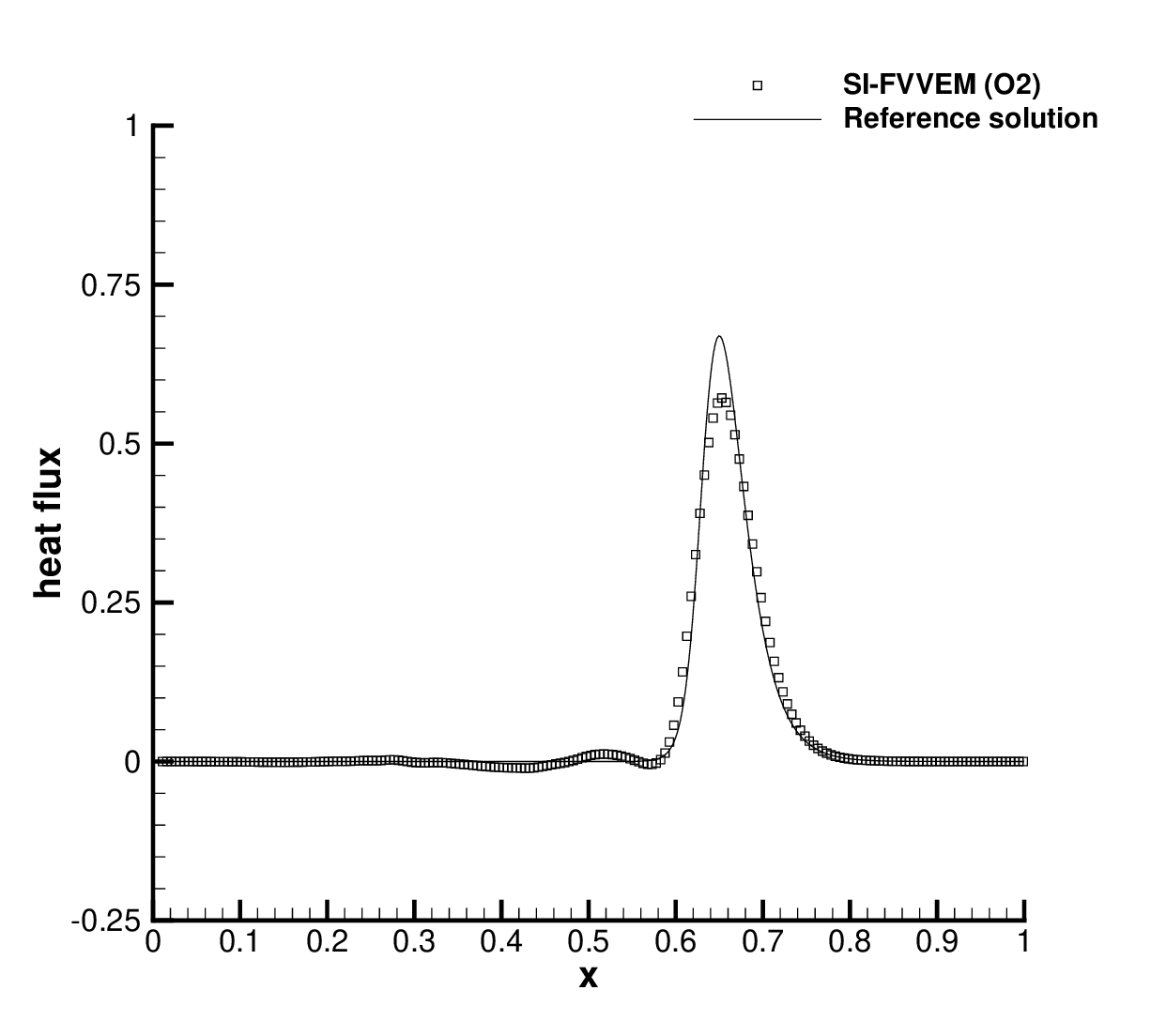} \\
		\end{tabular} 
		\caption{Viscous shock profile with shock Mach number $\mathrm{M_s}=2$ and Prandtl number $\Pr=0.75$ at time $t_f=0.2$. 
		Numerical solution compared against the reference solution for density $\rho$, horizontal velocity $u$, pressure $p$ and heat flux (from middle-left to bottom-right panels). We show a one-dimensional cut of 200 equidistant points along the $x-$direction at $y=0.1$.}
		\label{fig.ViscousShock}
	\end{center}
\end{figure}

To further study the behavior of the proposed scheme for high Mach number viscous flows, we consider the viscous shock benchmark whose exact solution for Pr$=0.75$ was derived by Becker \cite{Becker1923,BonnetLuneau,GPRmodel}. 
In particular, we set the fluid parameters as $c_v = 2.5$, $\mu=2 \cdot 10^{-2}$  and $\lambda = \left( 9+\frac{1}{3}\right)  \cdot 10^{-2}$ while the initial condition is detailed in \cite{ADERAFEDG}. The fluid is moving at shock Mach number $\text{M}_s=2$ with corresponding shock Reynolds number $\text{Re}_s=\frac{\rho_0 \, c_0 \, M_s \, L }{\mu}=100$.
The simulation is run on the computational domain $\Omega=[0,1] \times [0,0.4]$ discretized with $N=3217$ polygonal Voronoi elements 
of characteristic mesh spacing $h_x=1/200$, see Figure~\ref{fig.ViscousShockcont}. The initial shock wave is located at $x=0.25$. Figure~\ref{fig.ViscousShock} shows the comparison between the numerical solution and the exact solution at $t_{f}=0.2$. A very good agreement is observed for all reported quantities: density, velocity, pressure and heat flux. We remark that this
test case allows all terms contained in the Navier-Stokes system to be properly checked, since convection, thermal conduction and viscous stresses are present. Moreover, despite the underlying one-dimensional structure of the exact solution, this problem actually becomes multidimensional due to the use of unstructured meshes.

\subsection{Double shear layer}
\begin{figure}[!htbp]
	\begin{center}
		\begin{tabular}{cc}
			\includegraphics[trim=2 2 2 2,clip,width=0.47\textwidth]{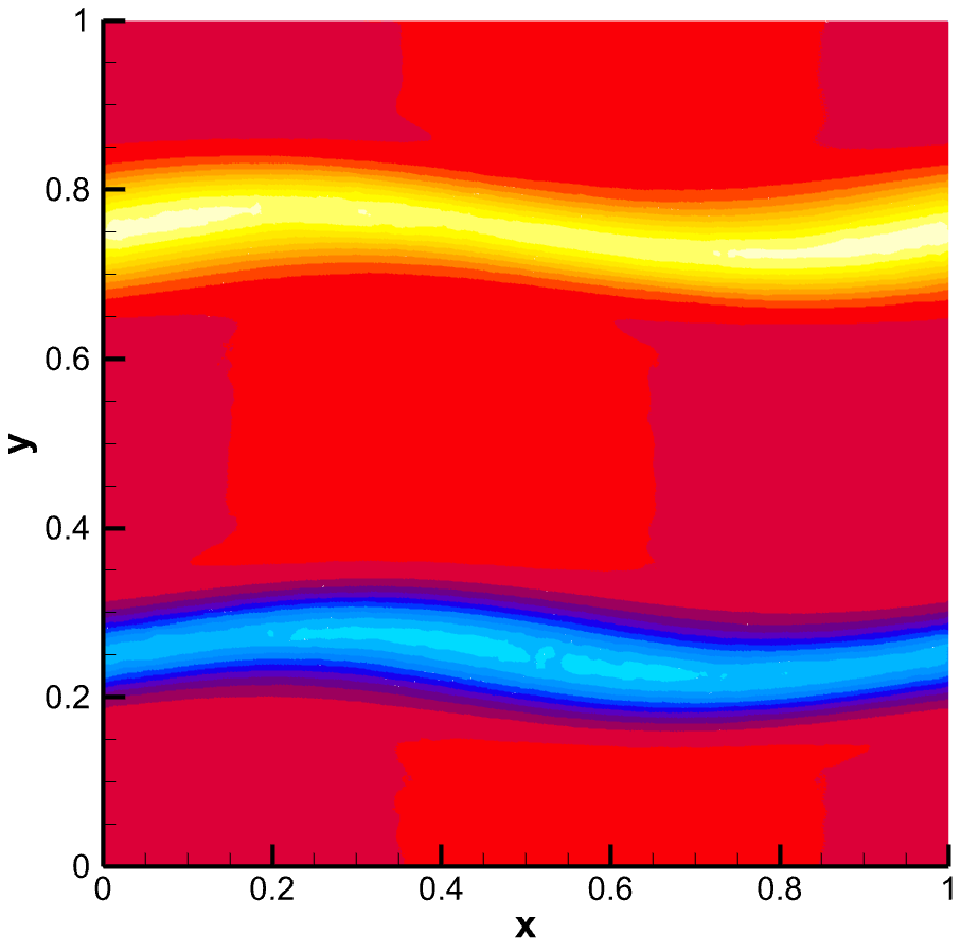}  &          
			\includegraphics[trim=2 2 2 2,clip,width=0.47\textwidth]{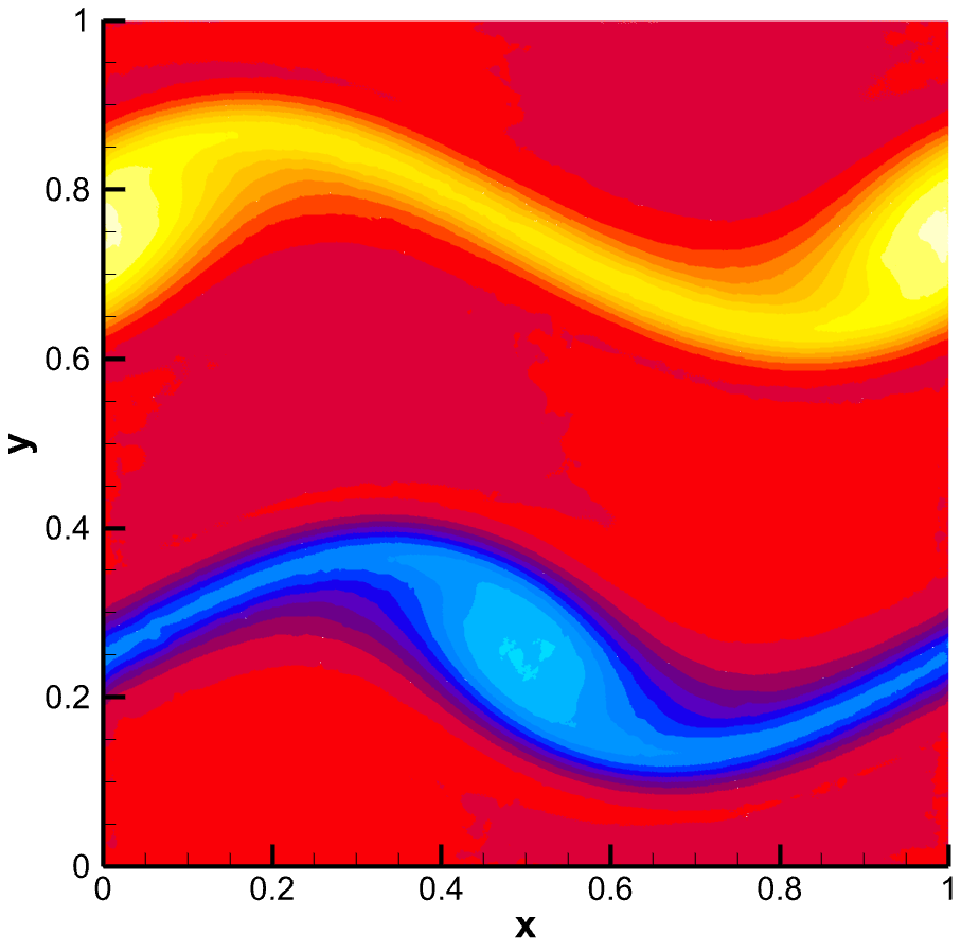}    \\
			\includegraphics[trim=2 2 2 2,clip,width=0.47\textwidth]{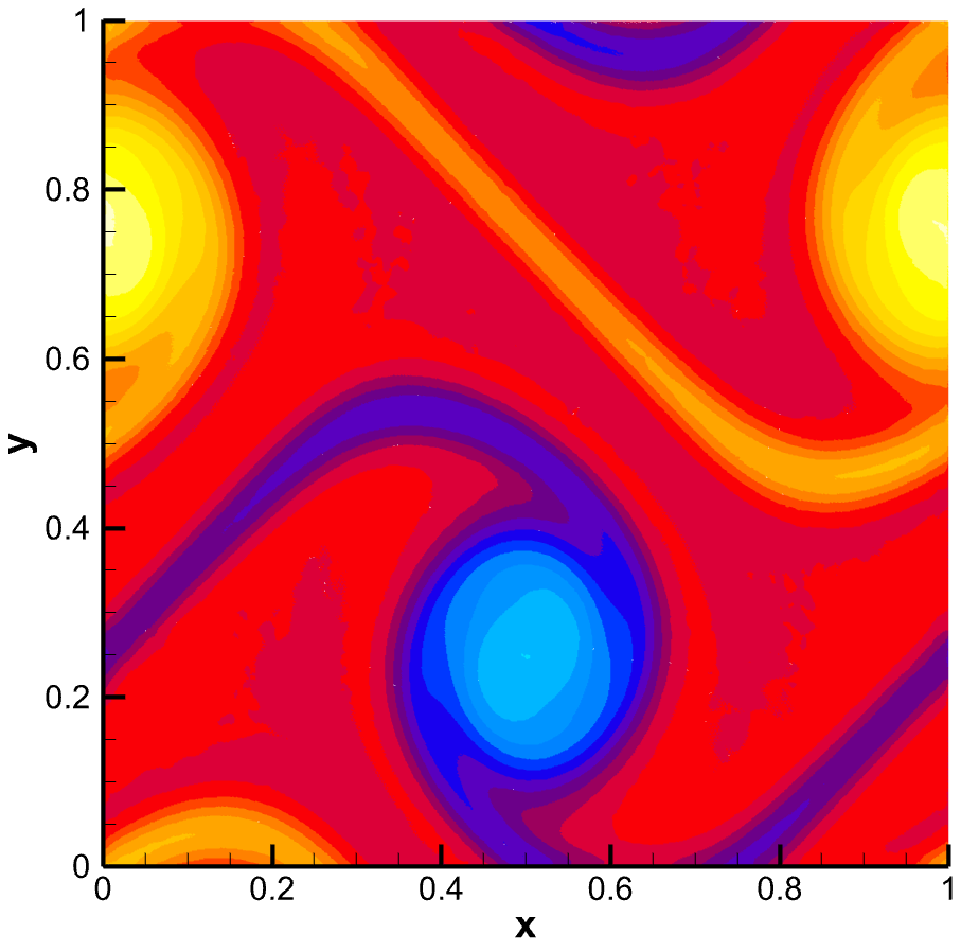}  &          
			\includegraphics[trim=2 2 2 2,clip,width=0.47\textwidth]{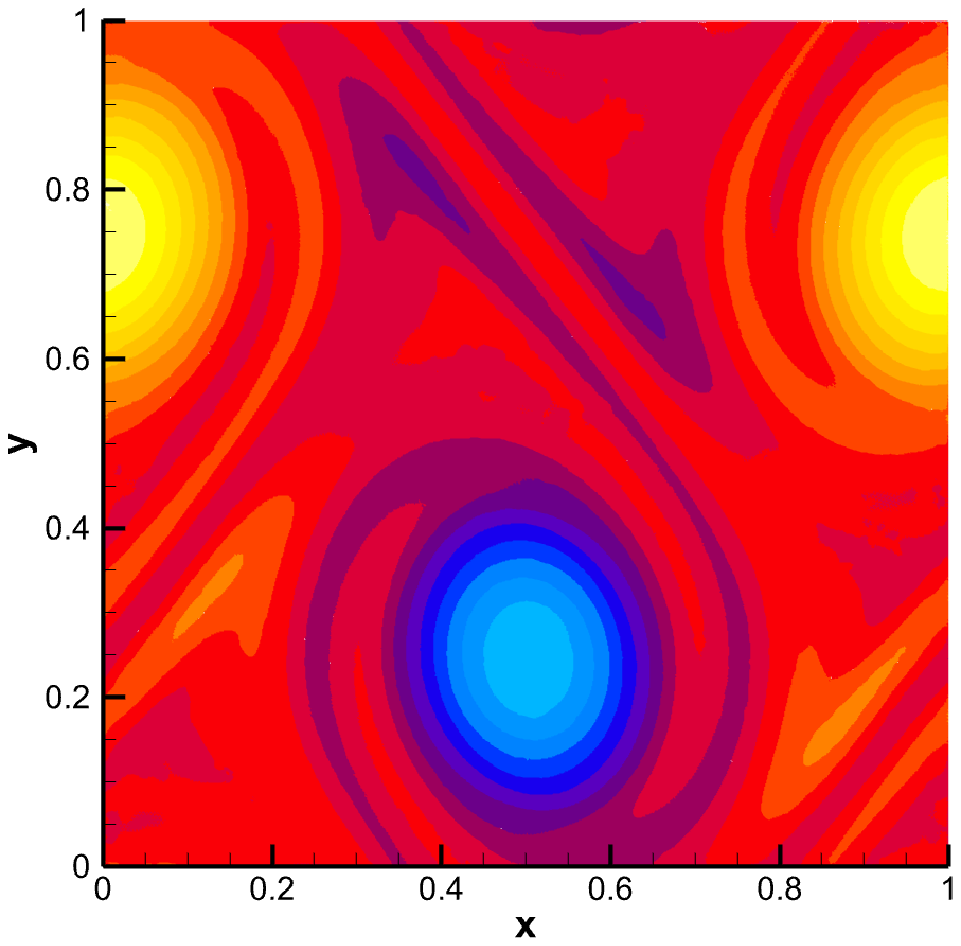} 
		\end{tabular}
		\caption{Double shear layer test. Vorticity at output times $t=0.4$, $t=0.8$, $t=1.2$ and $t=1.8$ (from top left ot bottom right). Plots use 21 equidistant contour lines in the range $[-26;26]$.}
		\label{fig.DSL}
	\end{center}
\end{figure}

This test involves high velocity gradients in the incompressible regime of the Navier-Stokes equations, and it has been proposed in \cite{BCG89}. The computational domain is $\Omega=[0,1]^2$ with periodic boundary conditions everywhere, and it is discretized with a mesh of size $h=1/200$, yielding a total number of $N=62355$ Voronoi cells. The fluid density is initially set to $\rho_0=1$, and the velocity field is prescribed by means of the following perturbed double shear layer profile:
\begin{equation}
	u=\left\{
	\begin{array}{l}
		\tanh\left( \theta \, (y-0.25) \right) \quad \textnormal{ if } y \leq 0.5, \\
		\tanh\left( \theta \, (0.75-y) \right) \quad \textnormal{ if } y > 0.5,
	\end{array}
	\right. \qquad v = \delta \sin(2\pi x), \label{eq:DSL_vel}
\end{equation} 
with $\theta=30$ and $\delta=0.05$. To approach the incompressible regime, the initial pressure is $p_0=10^4/\gamma$ so that a Mach number of $\text{M} \approx 10^{-2}$ is retrieved. The dynamic viscosity of the fluid is chosen to be $\mu=2\cdot 10^{-4}$ which corresponds to a Reynolds number of Re=$10000$. The simulation is run until the final time $t_f=1.8$ with the second order version of the SI-FVVEM scheme. Figure~\ref{fig.DSL} shows the distribution of the vorticity at different output times. Several vortices are generated by the dynamics of the flow, that are well captured by our numerical scheme, and which are in good agreement with results available in the literature \cite{Hybrid2,TD17,BT22}.

\subsection{Lid driven cavity}
\begin{figure}[!htbp]
	\begin{center}
		\begin{tabular}{cc} 
			\includegraphics[trim=2 2 2 2,clip,width=0.47\textwidth]{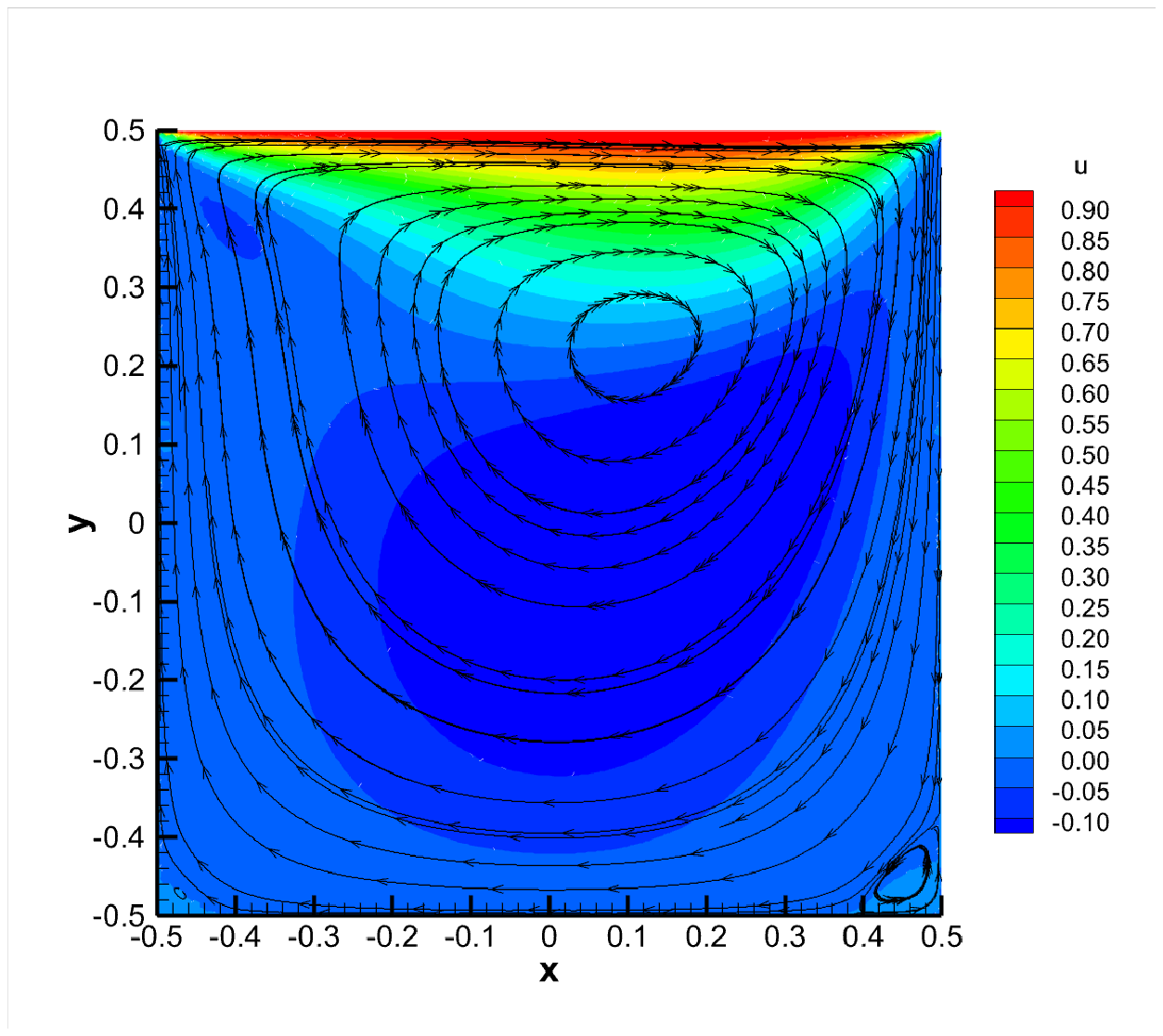} & 
			\includegraphics[trim=2 2 2 2,clip,width=0.47\textwidth]{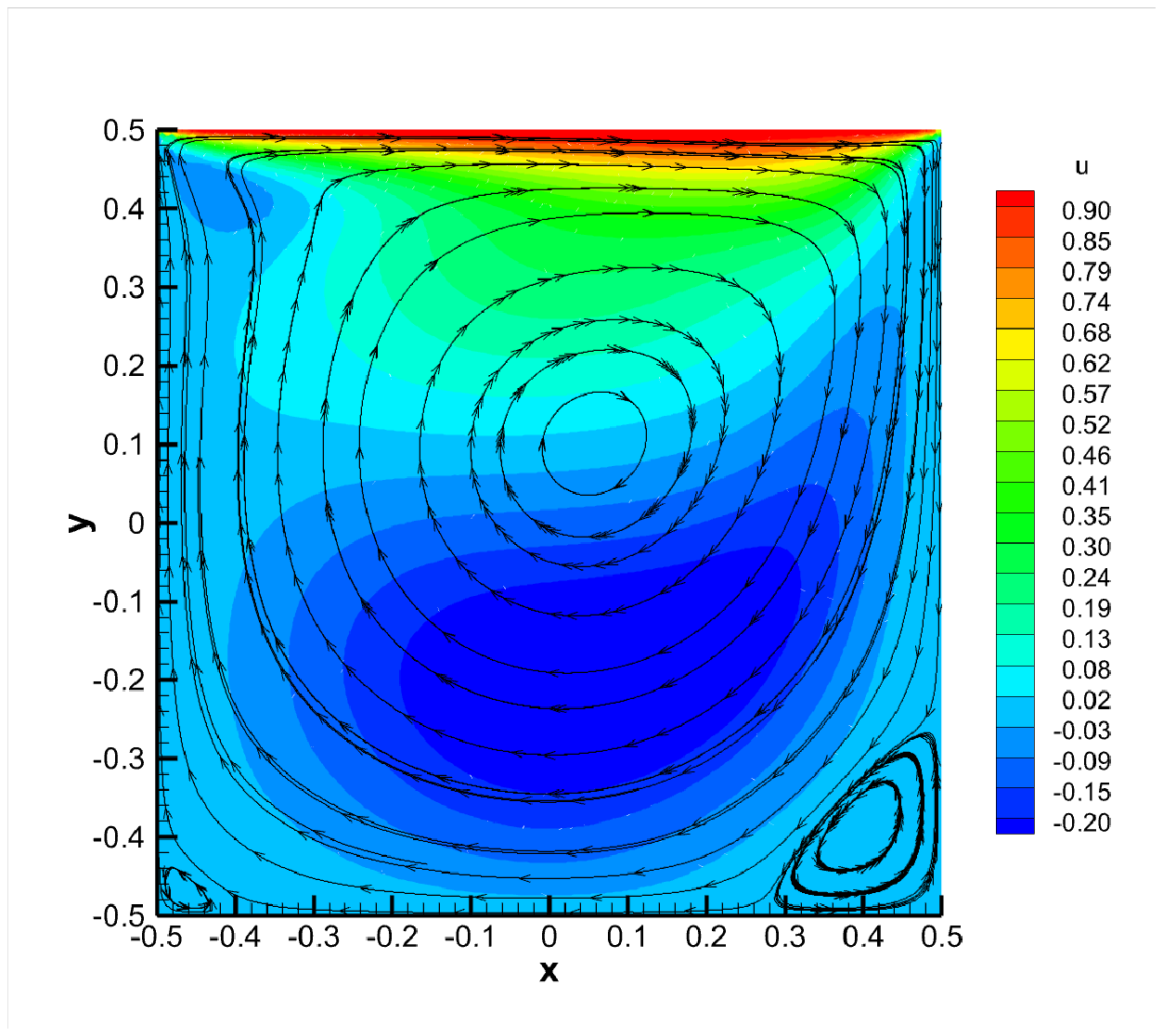} \\
			\includegraphics[width=0.47\textwidth]{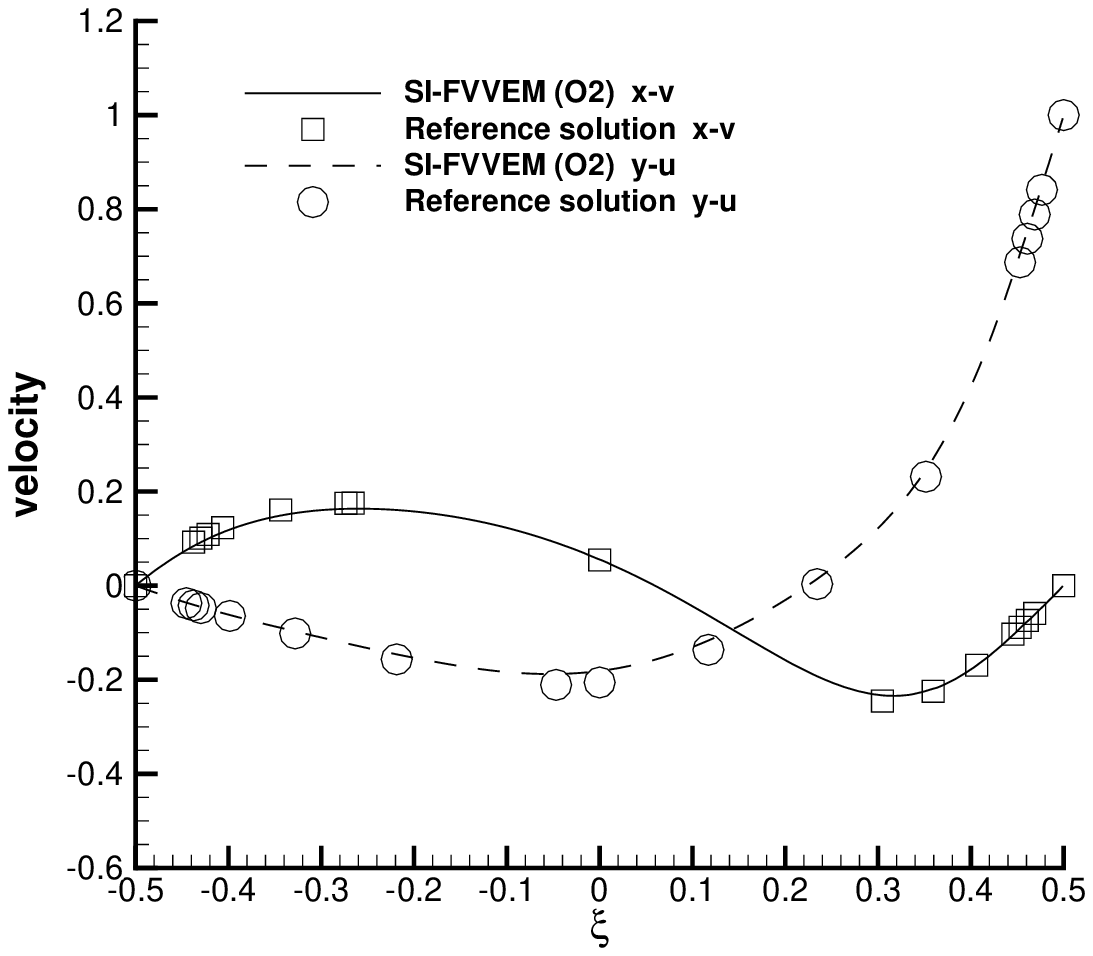} & 
			\includegraphics[width=0.47\textwidth]{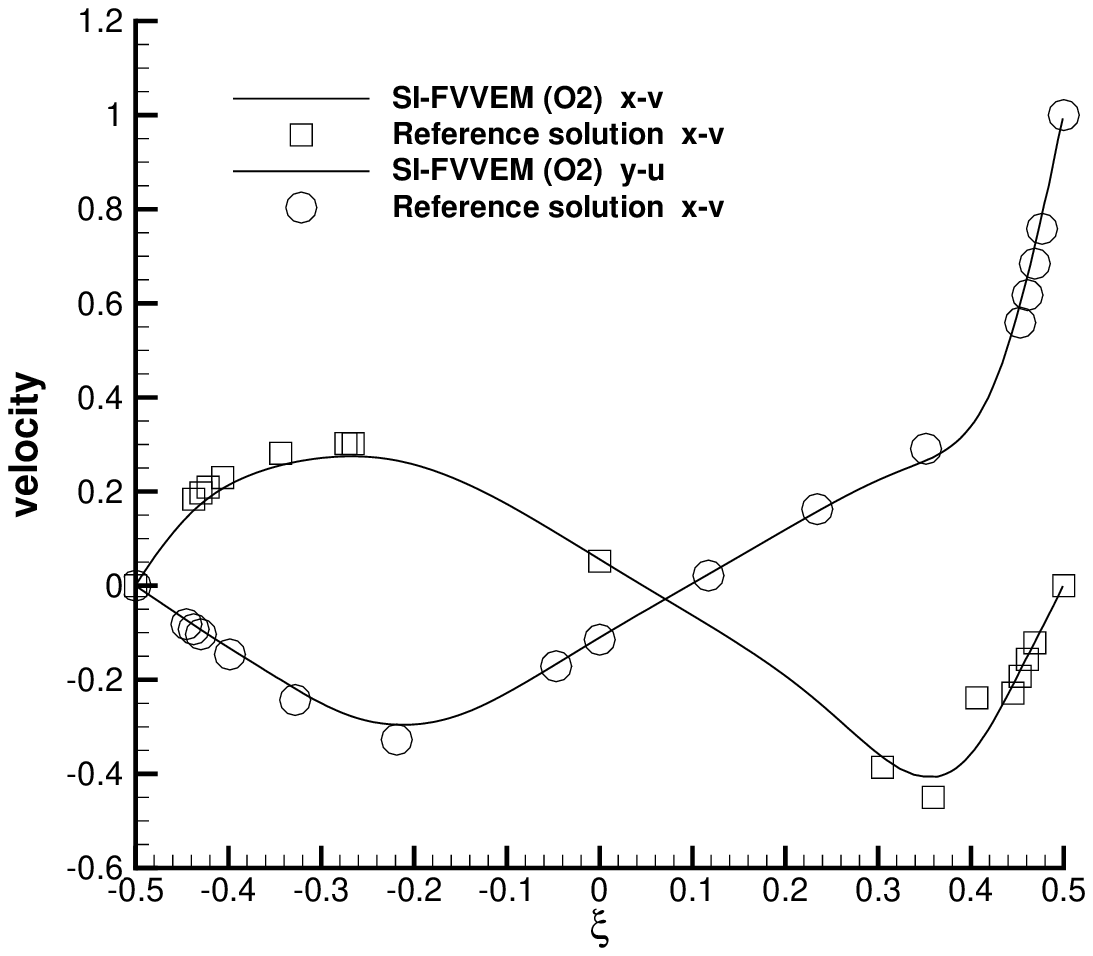} 
		\end{tabular} 
		\caption{Lid driven cavity flow at time $t_f=40$ for $\Rey=100$ (left) and $\Rey=400$ (right). Top: stream-traces of velocity on the computational domain $\Omega=[-0.5,0.5]^2$. Bottom: 1D cut comparison with the reference data in \cite{GGS82}.} 
		\label{fig.LidCavity1}
	\end{center}
\end{figure}

\begin{figure}[H]
	\begin{center}
		\begin{tabular}{cc} 
			\includegraphics[trim=2 2 2 2,clip,width=0.47\textwidth]{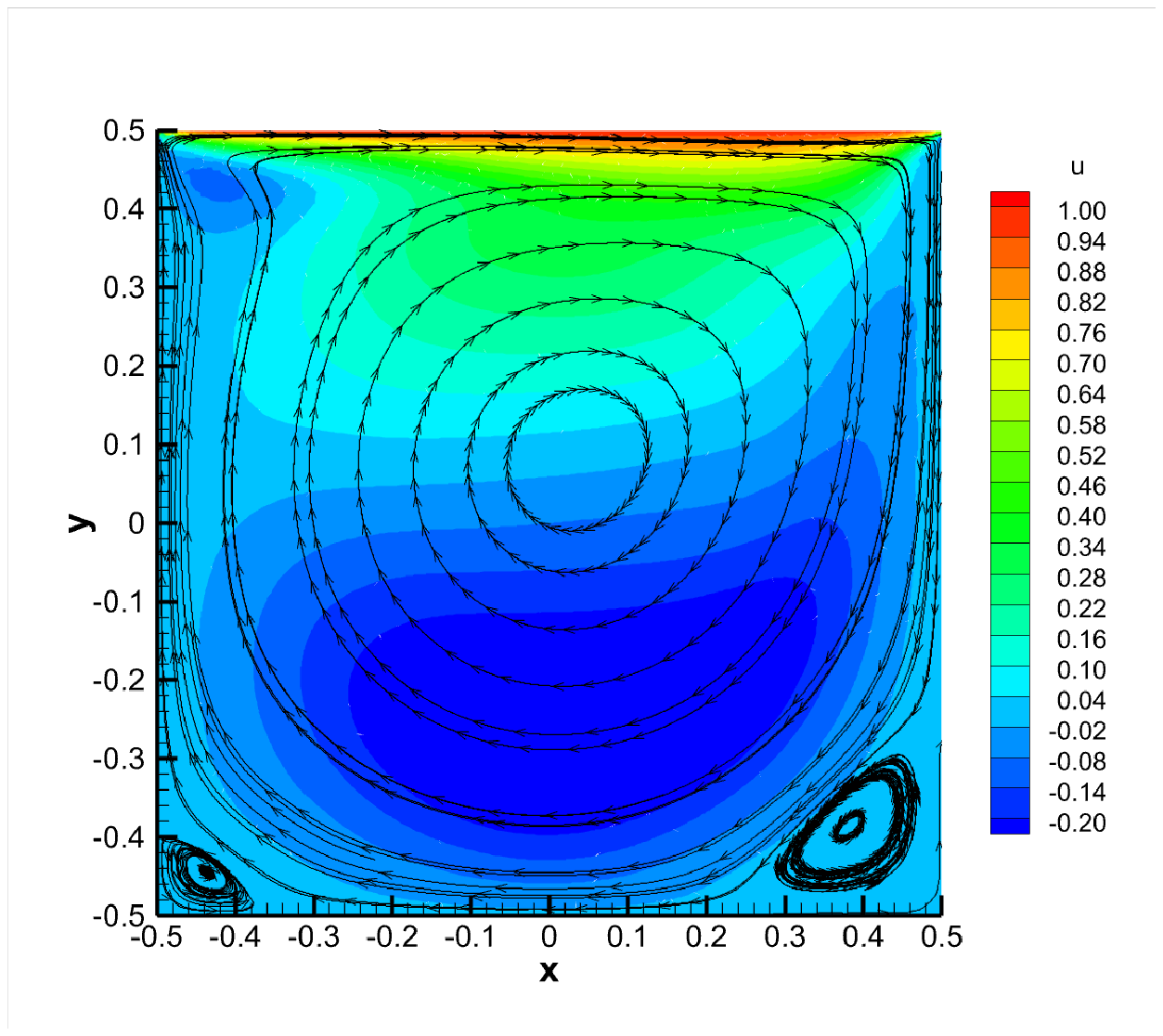} & 
			\includegraphics[trim=2 2 2 2,clip,width=0.47\textwidth]{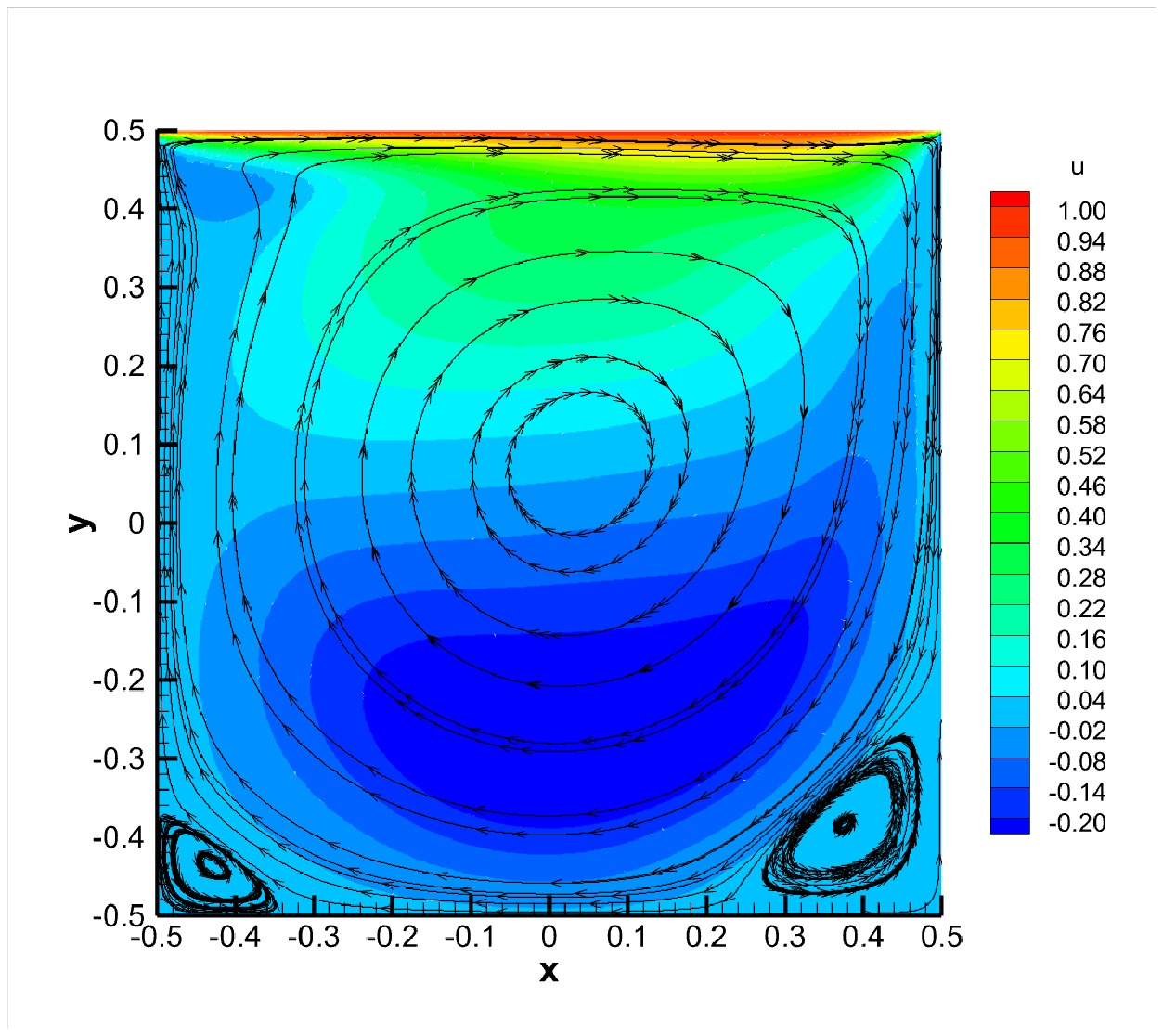} \\
			\includegraphics[width=0.47\textwidth]{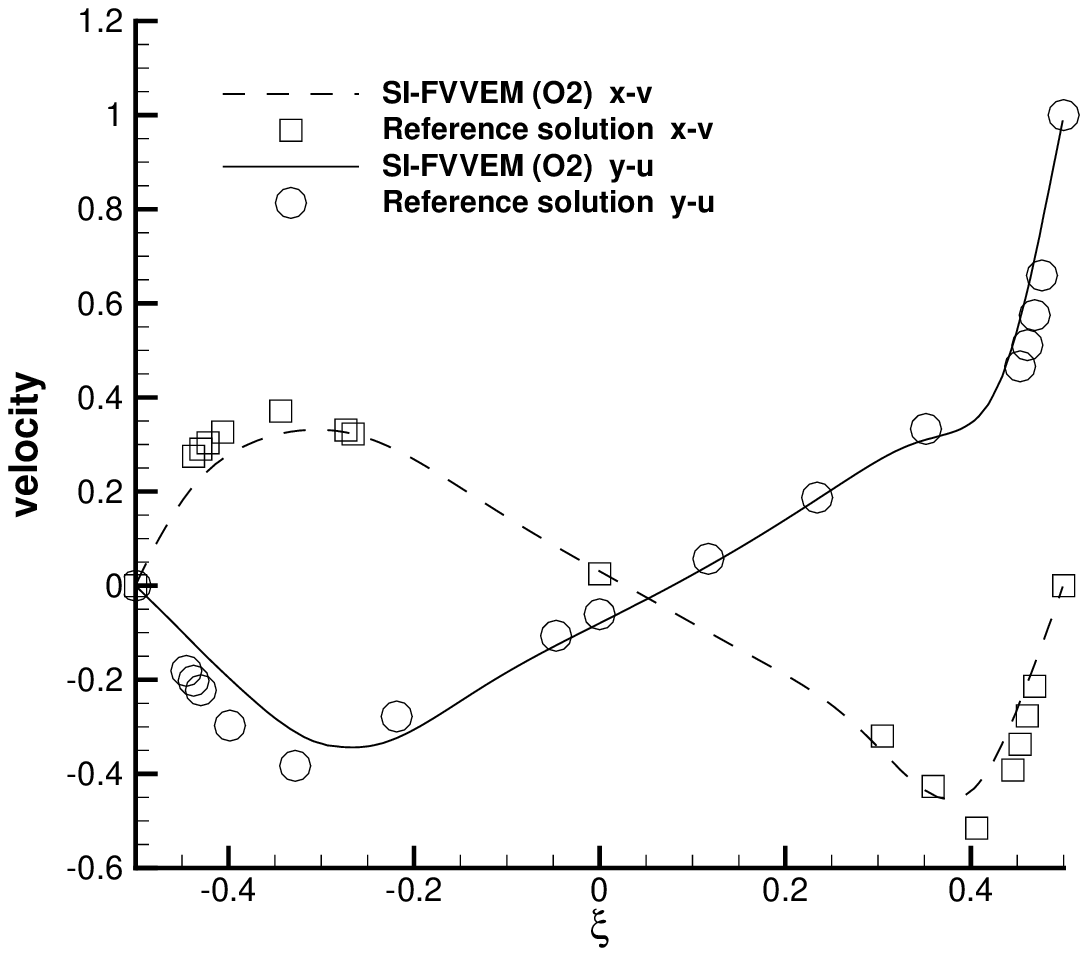} & 
			\includegraphics[width=0.47\textwidth]{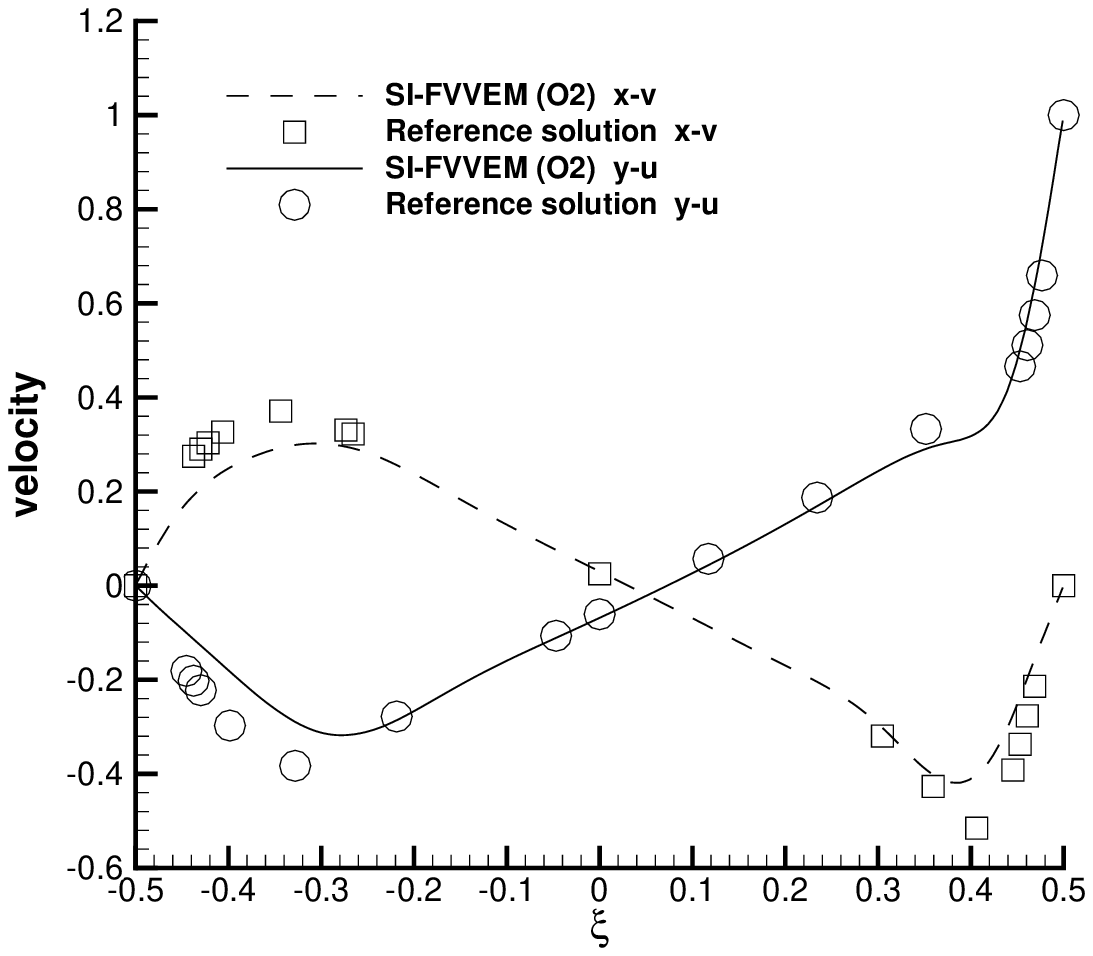} 
		\end{tabular} 
		\caption{Lid driven cavity flow at time $t_f=40$ for $\Rey=800$ (left) and $\Rey=1000$ (right). Top: stream-traces of velocity on the computational domain $\Omega=[-0.5,0.5]^2$. Bottom: 1D cut comparison with the reference data in \cite{GGS82}.} 
		\label{fig.LidCavity2}
	\end{center}
\end{figure}

As last test case, we analyze the lid driven cavity problem for a set of Reynolds numbers, namely $\Rey\in\left\lbrace 100, 400, 800, 1000 \right\rbrace$, see \cite{GGS82}. The computational domain $\Omega=[-0.5,0.5]^2$ is discretized using $N=6480$ polygons. No-slip wall boundary conditions are imposed everywhere apart from the upper boundary where we prescribe a moving velocity of $\vel(\xx,t)=(1,0)^{\top}$. We consider an initial fluid at rest with density $\rho(\xx,0)=1$ and pressure $\press(\xx,0)=10^{4}/\gamma$. The viscosity is chosen for each test case such that the desired Reynolds number is retrieved, while the Mach number is $\M \approx 8 \cdot 10^{-3}$ so we fall in the incompressible limit. The contour plots of the velocity field obtained using the second order SI-FVVEM scheme are depicted in Figure~\ref{fig.LidCavity1} for Re$=100$ and Re$=400$, and in Figure~\ref{fig.LidCavity2} for Re$=800$ and Re$=1000$. Moreover, the velocity profiles along the vertical and horizontal 1D cuts in the middle of the domain are compared with the corresponding reference solutions given in \cite{GGS82}, obtaining an overall very good agreement.

\section{Conclusions}\label{sec:conclusions}
We have presented a novel semi-implicit hybrid finite volume/virtual element method for all Mach number flows on general polygonal grids. 
The splitting of the governing equations into a convective sub-system, to be solved explicitly with FV, and a viscous and pressure sub-systems, discretized implicitly using VEM, allows the decoupling of the mean velocity waves from the fast pressure waves and the viscosity terms yielding an efficient scheme for both low Mach number flows and low Reynolds number flows. The finite volume method employed for the treatment of the non-linear convective terms leads to a robust scheme in the presence of large discontinuities which usually arise in compressible flows. On the other hand, the use of a virtual element approach to deal with the viscous and pressure Poisson-type problems permits to use arbitrary polygonal grids making the scheme suitable to deal with complex geometries. To attain high order of accuracy in time and space, we employ IMEX Runge-Kutta time integrators combined with high order CWENO reconstructions in the finite volume framework and high order virtual element schemes. An important property of the scheme is its asymptotic behavior in the incompressible limit. Asymptotic preserving properties are also ensured in the high Reynolds number limit \cite{HybridFVVEMinc}. Furthermore, the global energy conservation is proven to be guaranteed by solving a second pressure sub-system which accounts for the new kinetic energy. The proposed methodology has been successfully assessed for both supersonic flows and in the incompressible limit using a set of classical test cases. Viscous compressible flows in the high and low Mach number regime have also been successfully simulated by our novel SI-FVVEM schemes.

In the future, we foresee the extension of the proposed methodology to more complex systems of equations as the magnetohydrodynamics equations with its divergence-free condition which would require the development of a structure preserving scheme following the recent work presented in \cite{BTh24}. Moreover, moving mesh algorithms will be investigated in the framework of Arbitrary-Lagrangian-Eulerian schemes along the lines of \cite{LagrangeISO,SILVAins}, to account for moving curved control volumes even in the VEM context \cite{beirao2024curvedNCVEM}.

%
\section*{Acknowledgments}
WB received financial support by Fondazione Cariplo and Fondazione CDP (Italy) under the project No. 2022-1895 and by the Italian Ministry of University and Research (MUR) with the PRIN Project 2022 No. 2022N9BM3N. 
SB acknowledges the financial support from the Spanish Ministry of Science, Innovation and Universities (MCIN), the Spanish AEI (MCIN/AEI/10.13039/501100011033) and European Social Fund Plus under the project No. RYC2022-036355-I; from FEDER and the Spanish Ministry of Science, Innovation and Universities under project No. PID2021-122625OB-I00; and from the Xunta de Galicia (Spain) under project No. GI-1563 ED431C 2021/15.
MD was funded by the Italian Ministry of Education, University and Research (MIUR) in the frame of the Departments of Excellence  Initiative 2018--2027 attributed to DICAM of the University of Trento (grant L. 232/2016) and in the frame of the PRIN 2022 project \textit{High order structure-preserving semi-implicit schemes for hyperbolic equations}. 
MD was also co-funded by the European Union NextGenerationEU (PNRR, Spoke 7 CN HPC). Views and opinions expressed are however those of the author(s) only and do not necessarily reflect those of the European Union or the European Research Council. Neither the European Union nor the granting authorities can be held responsible for them.

WB and MD are members of the Gruppo Nazionale Calcolo Scientifico-Istituto Nazionale di Alta Matematica (GNCS-INdAM).

\bibliographystyle{plain}
\bibliography{./mibiblio}

\appendix
\section{IMEX schemes} \label{app.imex}
In what follows, we provide a detailed description of the two semi-implicit IMEX schemes introduced in Section~\ref{sec.imex} that have been employed for the simulation of the test cases presented in Section~\ref{sec:numericalresults}. Since our new semi-implicit scheme needs sub-stage values for the pressure, for which a conservative equation is not available, we provide the IMEX schemes in terms of \textit{states}, and not in terms of \textit{fluxes}, as presented in \cite{BFR16}. All the schemes employed in this work are stiffly accurate, therefore the final solution of the time integration coincides with the solution of the last stage of the Runge-Kutta method.

\paragraph{Second order semi-implicit IMEX scheme} The LSDIRK2(2,2,2) scheme is characterized by the following Butcher tableaux
\begin{equation}
	\begin{array}{c|cc}
		0 & 0 & 0 \\ \beta & \beta & 0 \\ \hline & 1-\gamma & \gamma
	\end{array} \qquad
	\begin{array}{c|cc}
		\gamma & \gamma & 0 \\ 1 & 1-\gamma & \gamma \\ \hline & 1-\gamma & \gamma
	\end{array}
	\label{eqn.IMEX2}
\end{equation}
with $\gamma=1-1/\sqrt{2}$ and $\beta=1/(2\gamma)$. Hence, to solve the autonomous system  \eqref{eqn.autonomous}, we perform the following steps.
\begin{description}
	\item[Step 1.] For $i=1$, we have
	\begin{eqnarray}
		&&\QE^{(1)} = \Q^{n},\label{eqn.LSDIRK2_step1a}\\
		&&\QI^{(1)} = \Q^{n} + \gamma \dt \mathcal{H}\left(\QE^{(1)},\QI^{(1)}\right) = \Q^{n} + \gamma \dt \mathcal{H}\left(\Q^{n},\QI^{(1)}\right). \label{eqn.LSDIRK2_step1b}
	\end{eqnarray}
	Thus, reordering the terms in \eqref{eqn.LSDIRK2_step1b}, we get 
	\begin{equation}
		\mathcal{H}\left(\QE^{(1)},\QI^{(1)}\right) = \frac{1}{\gamma \dt}\left( \QI^{(1)} - \Q^{n}\right).
	\end{equation}
	\item[Step 2.] For $i=2$, it results
	\begin{eqnarray}
		\QE^{(2)} &\!=&\! \Q^{n} + \dt \beta \mathcal{H}\left(\QE^{(1)},\QI^{(1)}\right) 
		= \frac{1-\beta}{\gamma }\Q^{n} +   \frac{\beta}{\gamma } \QI^{(1)},   \label{eqn.LSDIRK2_step2a}\\
		\QI^{(2)} 
		&\!=&\! \Q^{n} + \dt\left(1-\gamma\right)\mathcal{H}\left(\QE^{(1)},\QI^{(1)}\right)  + \gamma\dt \mathcal{H}\left(\QE^{(2)},\QI^{(2)}\right) \notag\\
		&\!=&\! \frac{1}{\gamma }\Q^{n} + \frac{1-\gamma}{\gamma}\QI^{(1)}  + \gamma\dt \mathcal{H}\left(\QE^{(2)},\QI^{(2)}\right). \label{eqn.LSDIRK2_step2b}
	\end{eqnarray}
	Hence,
	\begin{equation}		
		\mathcal{H}\left(\QE^{(2)},\QI^{(2)}\right) = - \frac{1}{\gamma^2 \dt }\Q^{n} + \frac{1}{\gamma\dt} \QI^{(2)}  - \frac{1-\gamma}{\gamma^2\dt}\QI^{(1)} .
	\end{equation}
	
	\item[Step 3.] Once we have computed $k_i\mathcal{H}\left(\QE^{(i)},\QI^{(i)}\right)$, $\forall i\in\left\lbrace 1,\dots,s\right\rbrace$, the numerical solution at the new time step is calculated from \eqref{eqn.QRKfinal}, thus we get $\Q^{n+1} = \QI^{(2)}$.
	
\end{description}

\paragraph{Third order semi-implicit IMEX scheme} The SA DIRK(3,4,3) scheme is described by the following Butcher tableaux with $\gamma=0.435866$: 
\begin{equation}
	\begin{array}{c|cccc}
		0 & 0 & 0 & 0 & 0 \\ \gamma & \gamma & 0 & 0 & 0 \\ 0.717933 & 1.437745 & -0.719812 & 0 & 0 \\ 1 & 0.916993 & 1/2 & -0.416993 & 0 \\ \hline  & 0 & 1.208496 & -0.644363 & \gamma
	\end{array} \qquad
	\begin{array}{c|cccc}
		\gamma & \gamma & 0 & 0 & 0 \\ \gamma & 0 & \gamma & 0 & 0 \\ 0.717933 & 0 & 0.282066 & \gamma & 0  \\ 1 & 0 & 1.208496 & -0.644363 & \gamma \\ \hline  & 0 & 1.208496 & -0.644363 & \gamma
	\end{array}
	\label{eqn.IMEX3}
\end{equation}
Consequently, the steps of the corresponding algorithm read as follows.
\begin{description}
	\item[Step 1.] $i=1$:
	\begin{eqnarray}
		\QE^{(1)} &\!=&\! \Q^{n},\\
		\QI^{(1)} &\!= &\!\Q^{n} + \gamma \dt \mathcal{H}\left(\QE^{(1)},\QI^{(1)}\right).
	\end{eqnarray}
	Then,
	\begin{equation}
		\mathcal{H}\left(\QE^{(1)},\QI^{(1)}\right) = \frac{1}{\gamma \dt}\left( \QI^{(1)} - \Q^{n}\right).
		\label{eqn.SADIRK_step1}
	\end{equation}
	
	\item[Step 2.] $i=2$:
	\begin{eqnarray}
		\QE^{(2)} &\!=&\! \Q^{n} + \dt \gamma \mathcal{H}\left(\QE^{(1)},\QI^{(1)}\right) 
		=  \QI^{(1)},   \label{eqn.SADIRK_step2a}\\
		\QI^{(2)} 
		&\!=&\! \Q^{n} + \dt\gamma\mathcal{H}\left(\QE^{(2)},\QI^{(2)}\right).
	\end{eqnarray}
	So,
	\begin{equation}
		\mathcal{H}\left(\QE^{(2)},\QI^{(2)}\right) = \frac{1}{\gamma\dt}\left(\QI^{(2)}  -\Q^{n}\right). 
		\label{eqn.SADIRK_step2b}
	\end{equation}
	
	\item[Step 3.] $i=3$:
	\begin{eqnarray}
		\QE^{(3)} &\!=&\! \Q^{n} +\tilde{a}_{31}\dt \mathcal{H}\left(\QE^{(1)},\QI^{(1)}\right) 
		+\tilde{a}_{32}\dt \mathcal{H}\left(\QE^{(2)},\QI^{(2)}\right)\notag \\
		&\!=&\! \left( 1 -  \frac{\tilde{a}_{31}}{\gamma } - \frac{\tilde{a}_{32} }{\gamma}\right) \Q^{n} +\frac{\tilde{a}_{31}}{\gamma } \QI^{(1)} 
		+  \frac{\tilde{a}_{32}}{\gamma}\QI^{(2)},\\
		\QI^{(3)} &\!=&\! \Q^{n} + a_{32} \dt \mathcal{H}\left(\QE^{(2)},\QI^{(2)}\right) + \gamma\dt \mathcal{H}\left(\QE^{(3)},\QI^{(3)}\right) \notag\\
		&\!=&\! \Q^{n} + \frac{a_{32}}{\gamma}\left(\QI^{(2)}  -\Q^{n}\right)
		+ \gamma\dt \mathcal{H}\left(\QE^{(3)},\QI^{(3)}\right).
	\end{eqnarray}
	Thus,
	\begin{equation}
		\mathcal{H}\left(\QE^{(3)},\QI^{(3)}\right)	 =  \left( - \frac{1}{\gamma\dt } + \frac{a_{32}}{\gamma^2\dt }\right) \Q^{n} - \frac{a_{32}}{\gamma^2\dt }\QI^{(2)} + \frac{1}{\gamma\dt }\QI^{(3)}.
	\end{equation}

	\item[Step 4.] $i=4$:
	\begin{eqnarray}
		\QE^{(4)} &\!\!\!=&\!\!\! \Q^{n} + \tilde{a}_{41} \dt \mathcal{H}\left(\QE^{(1)},\QI^{(1)}\right) 
		+ \tilde{a}_{42} \dt \mathcal{H}\left(\QE^{(2)},\QI^{(2)}\right)
		+ \tilde{a}_{43} \dt \mathcal{H}\left(\QE^{(3)},\QI^{(3)}\right) \notag\\
		&\!\!\!=&\!\!\! \left( 1 - \frac{ \tilde{a}_{41}}{\gamma } -\frac{\tilde{a}_{42}}{\gamma}  - \frac{\tilde{a}_{43} }{\gamma } + \frac{\tilde{a}_{43} a_{32}}{\gamma^2 }\right) \Q^{n} 
		+ \frac{ \tilde{a}_{41}}{\gamma } \QI^{(1)} 
		+ \left( \frac{\tilde{a}_{42}}{\gamma}  - \frac{\tilde{a}_{43} a_{32}}{\gamma^2 }\right) \QI^{(2)} + \frac{\tilde{a}_{43} }{\gamma }\QI^{(3)}, \\
		\QI^{(4)} &\!\!\!=&\!\!\! \Q^{n} +a_{42} \dt  \mathcal{H}\left(\QE^{(2)},\QI^{(2)}\right)  
		+ a_{43} \dt \mathcal{H}\left(\QE^{(3)},\QI^{(3)}\right)  +  \gamma\dt \mathcal{H}\left(\QE^{(4)},\QI^{(4)}\right) \notag\\
		&\!\!\!=&\!\!\! \left( 1  - \frac{a_{42}}{\gamma}  - \frac{a_{43} }{\gamma } 
		+ \frac{a_{43} a_{32}}{\gamma^2 }\right) \Q^{n}
		+ \left(   \frac{a_{42}}{\gamma} 
		- \frac{a_{43} a_{32}}{\gamma^2 }\right) \QI^{(2)} 
		+ \frac{a_{43} }{\gamma }\QI^{(3)}  
		+  \gamma\dt \mathcal{H}\left(\QE^{(4)},\QI^{(4)}\right). \notag\\
	\end{eqnarray}
	Hence,
	\begin{gather}
		\mathcal{H}\left(\QE^{(4)},\QI^{(4)}\right) 
		= 
		- \frac{1}{\gamma\dt}\left( 1  - \frac{a_{42}}{\gamma}  - \frac{a_{43} }{\gamma } 
		+ \frac{a_{43} a_{32}}{\gamma^2 }\right) \Q^{n}
		- \frac{1}{\gamma\dt}\left(   \frac{a_{42}}{\gamma} 
		- \frac{a_{43} a_{32}}{\gamma^2 }\right) \QI^{(2)} \notag\\
		- \frac{1}{\gamma\dt}\frac{a_{43} }{\gamma }\QI^{(3)}  
		+ \frac{1}{\gamma\dt} \QI^{(4)} .
	\end{gather}
	
	\item[Step 5.] Finally, $\Q^{n+1} = \QI^{(4)}$. 
\end{description}
\end{document}